\documentclass[11pt]{article}

\usepackage{color}
\usepackage{latexsym}
\usepackage{dsfont}
\usepackage{amssymb}
\usepackage{graphicx}
\usepackage{amsmath, amsfonts,amssymb,theorem,euscript,array,enumerate,amsfonts,mathrsfs}
\usepackage{hyperref}
\usepackage{appendix}

\usepackage{tikz}
\usetikzlibrary{fit,matrix,chains,positioning,decorations.pathreplacing,arrows}

\usepackage{mathtools}
\mathtoolsset{showonlyrefs}

\definecolor{darkgreen}{rgb}{0,0.7,0}

\newcommand{\independent}{\protect\mathpalette{\protect\independenT}{\perp}}
\def\independenT#1#2{\mathrel{\rlap{$#1#2$}\mkern2mu{#1#2}}}
\newcommand{\ba}{{\bf a}}

\def \b1{\bf{1}}

\def \N{\mathbb{N}}
\def \R{\mathbb{R}}

\def \Z{\mathbb{Z}}
\def \E{\mathbb{E}}
\def \F{\mathbb{F}}

\def \P{\mathbb{P}}

\def \T{\mathbb{T}}

\def \S{\mathbb{S}}

\def \bp{\boldsymbol{p}} 
\def \bpi{\boldsymbol{\pi}} 
 
\def \bF{\tilde F} 
 
\def \bd{\boldsymbol{d}} 
\def \bW{\boldsymbol{W}} 
 
\def \bnu{\boldsymbol{\nu}} 
\def \bA{\boldsymbol{A}} 
 
\def \balpha{\boldsymbol{\alpha}} 
\def \bAc{\boldsymbol{\Ac}} 
\def \mfa{\mathfrak{a}} 
 
\def \mra{\mathrm{a}} 
 
\def \mfd{\mathrm{d}} 
\def \d{\mathrm{d}} 
\def \bLi{\boldsymbol{\Lc}}

\def\esssup_#1{\underset{#1}{\mathrm{ess\,sup\, }}}

\def\argmin_#1{\underset{#1}{\mathrm{argmin\, }}}
\def\argmax_#1{\underset{#1}{\mathrm{argmax\, }}}

\def \Ac{{\cal A}}
\def \Bc{{\cal B}}
\def \Cc{{\cal C}}

\def \Fc{{\cal F}}
\def \Gc{{\cal G}}

\def \Lc{{\cal L}}
\def \Pc{{\cal P}}
\def \Rc{{\cal R}}

\def \Mc{{\cal M}}

\def \Oc{{\cal O}}
\def \Sc{{\cal S}}
\def \Tc{{\cal T}}
\def \Uc{{\cal U}}

\def \Wc{{\cal W}}
\def \Xc{{\cal X}}
\def \Yc{{\cal Y}}

\def \eps{\varepsilon}

\def \ep{\hbox{ }\hfill$\Box$}

\def\beqs{\begin{eqnarray*}}
\def\enqs{\end{eqnarray*}}
\def\beq{\begin{eqnarray}}
\def\enq{\end{eqnarray}}

\newtheorem{Theorem}{Theorem}[section]
\newtheorem{Definition}{Definition}[section]
\newtheorem{Proposition}{Proposition}[section]

\newtheorem{Lemma}{Lemma}[section]

\newtheorem{Remark}{Remark}[section]
\newtheorem{Example}{Example}[section]
\numberwithin{equation}{section}

\begin{document}

\title{
Mean-field Markov decision processes  with common noise and open-loop controls
}

\author{
M\'ed\'eric MOTTE
\footnote{LPSM, Universit\'e de Paris  \sf medericmotte at gmail.com
The author acknowledges support of the DIM MathInnov. 
}
\qquad\quad
Huy\^en PHAM
\footnote{LPSM, Universit\'e de  Paris, and CREST-ENSAE, \sf pham at lpsm.paris
This work was partially supported by the Chair Finance \& Sustainable Development / the FiME Lab (Institut Europlace de Finance)
}
}

\date{This version: june 2021 \\ to appear in {\it Annals of Applied Probability}}

\maketitle

\begin{abstract}
We develop an exhaustive study of Markov decision process (MDP) under mean field interaction both on states and actions  in the presence of common noise, and when optimization is performed over open-loop controls on infinite horizon. Such model, called CMKV-MDP for conditional McKean-Vlasov MDP, arises and is obtained here rigo\-rously  with a rate of convergence as the asymptotic problem of $N$-cooperative agents controlled by a social planner/influencer that observes the environment noises but not necessarily the individual states of the agents. 
We highlight  the crucial role of relaxed controls and randomization hypothesis for this class of models with respect to classical MDP theory. 
We prove the correspondence between CMKV-MDP and a general lifted MDP on the space of pro\-bability measures,  and establish the dynamic programming Bellman fixed point equation  satisfied by the value function, as well as the existence of $\epsilon$-optimal randomized feedback controls. The arguments of proof  involve an original measurable optimal coupling for the Wasserstein distance. 
This provides a procedure for learning strategies in a  large population of interacting collaborative agents.   
\end{abstract}

\vspace{5mm}


\vspace{5mm}


\maketitle

\newpage

\tableofcontents

\section{Introduction}


Optimal  control of McKean-Vlasov (MKV)  systems, also known as mean-field  control (MFC)   pro\-blems,  has sparked a great interest in the domain of applied probabilities during the last decade. 
In these optimization problems, the transition dynamics of the system and the reward/gain function depend not only on the state and action of the agent/controller, but also on their probability distributions.   
These problems are motivated from models of large population of interacting cooperative agents obeying to a social planner (center of decision), and are often justified heuristically 
as the asymptotic  regime with infinite number of agents under  Pareto efficiency. Such problems have found numerous applications in distributed energy, herd behavior, finance, etc.

A large literature has already emerged on continuous-time models for the optimal control of McKean-Vlasov dynamics, and  dynamic programming principle (in other words time consistency) has been established in this context in the papers  \cite{laupir16},  \cite{phawei17}, \cite{bayetal}, \cite{djepostan19}.  We point out  the work \cite{lac17}, which is the first paper to rigorously connect mean-field control to large systems of controlled processes, 
see also the recent paper  \cite{foretal18}, 
and  refer to the books \cite{benetal13}, \cite{cardelbook} for an overview of the subject.

\vspace{1mm}

\noindent {\bf Our work and main contributions.}  In this paper, we introduce  a general discrete time framework by providing an exhaustive study of Markov decision process (MDP) under mean-field interaction 
in the presence of {\it common noise}, and when optimization is performed over {\it open-loop controls} on infinite horizon. Such model is called conditional McKean-Vlasov MDP, shortly abbreviated in  the sequel as CMKV-MDP, and 
the set-up is the mathematical framework for a theory of reinforcement learning with mean-field interaction. 
Let us first briefly describe and motivate the main features of our framework:
\begin{enumerate}
    \item[(i)]  The controls are open-loop: it emulates the generally unobservable nature of the states in mean-field situations (political opinions, targeted advertising, 
    health conditions). Yet, the noises driving these states (people's characteristics, habits, surroundings, and navigation data) are nowadays observable.
    \item[(ii)]  The controls are {\em individualized}: the action applied to an individual only depends upon her own noise. Indeed, for the strategy to be implementable in practice, the action must be computed by a program running on the individual's own device (computer, phone), and, to limit data transfers, we impose that it only uses this device's data.
    \item[(iii)]  The dynamics of individuals depend upon a common noise, emulating the fact that they are influenced by common information (public data) which may vary over time. By considering an i.i.d. common noise sequence, we shall see that the value function is characterized as the unique fixed point of a Bellman operator, which is the basis for designing algorithms in reinforcement learning. 
\end{enumerate}


Compared to continuous-time models, discrete-time McKean-Vlasov control problems have been less studied in the literature. In  \cite{phawei16}, the authors consider a finite-horizon problem 
without common noise and state the dynamic programming (Bellman) equation for MFC with closed-loop (also called feedback) controls, that are restricted to depend on the state.  Very  recently, the works 
\cite{carlautan19}, \cite{guetal19} addressed Bellman equations for MFC problems in the context of reinforcement learning.  The paper \cite{guetal19} considers relaxed controls in their MFC formulation but without common noise, 
and derives the Bellman equation for the $Q$-value function as a deterministic control problem that we obtain here as a particular case (see our Remark \ref{remQ}).  The framework  in \cite{carlautan19} is closest to ours  by considering also common noise, however with the following differences: these authors restrict their attention to stationary feedback policies, 
and reformulate their MFC control problem as a MDP on the space of probability measures by deriving formally (leaving aside the measurability issues and assuming the existence of a stationary feedback control) 
the associated Bellman equation, which is then used for the development of $Q$-learning algorithms. Notice that \cite{carlautan19}, \cite{guetal19}  do not consider dependence upon the probability distribution of the control in 
 the state transition dynamics and reward function.

Besides the introduction of a general framework including a mean-field dependence on the pair state/action, our first main contribution is to rigorously connect CMKV-MDP to a large but finite system of MDP with interacting processes. We prove the almost sure and  $\mathbb{L}^1$ conditional propagation of chaos, i.e.,  the convergence, as the number of interacting agents  $N$ tends to infinity, of the state processes and gains of the $N$-individual population control problem towards the corresponding object in the CMKV-MDP. Furthermore, by relying on rate of convergence in Wasserstein distance of the empirical measure, we give a rate of convergence 
for the limiting CMKV-MDP under suitable Lipschitz assumptions on the state transition and reward functions, which is new to the best of our knowledge.

Our second contribution is to obtain the correspondence of  our CMKV-MDP  with a suitable lifted MDP on the space of probability measures. Starting from open-loop controls, this is achieved in general    
by introducing relaxed (i.e. measure-valued) controls in the enlarged state/action space, and by emphasizing  the measu\-rability issues arising in the presence of common noise and with continuous state space. In the special case without common noise or with finite state space, the relaxed control in the lifted MDP  is reduced to the usual notion in control theory,  
also known as mixed or randomized strategies in game theory. 
While it is known in standard MDP that an optimal control (when it exists) is in pure form,  relaxed control appears 
naturally in MFC where the social planner has to sample the distribution of actions  instead of simply  assigning  the same pure strategy among the population in order to perform the best possible collective gain. 

The reformulation of the original problem as a lifted MDP leads us to consider an associated dynamic programming equation written in terms of a Bellman fixed point equation  in the space of probability measures. 
Our third contribution is to establish rigorously  the Bellman equation satisfied by the state value function of the CMKV-MDP, and then by the state-action value function, called $Q$-function in the reinforcement learning termi\-nology. 
This is obtained under the crucial assumption that the initial information filtration is gene\-rated by an atomless random variable, i.e., that it is rich enough, and calls upon original measurable optimal  coupling results for the Wasserstein distance.  Moreover, and this is our fourth contribution,  the methodology of proof allows us to 
obtain as a by-product the existence of an $\epsilon$-optimal control, which is constructed from randomized feedback policies under a randomization hypothesis.   This shows in particular that the value function 
of CMKV-MDP over open-loop controls is equal to the value function over randomized feedback controls, and we highlight that it may be strictly larger than the value function of CMKV-MDP over ``pure"  feedback controls, i.e., without randomization. This is a notable difference with respect to  the classical (without mean-field dependence) theory of MDP as studied e.g. in  \cite{Bertsekas}, \cite{sutbar}. We discuss and illustrate with a set of simple examples the diffe\-rence of control strategies (open loop vs feedback), and the crucial role of the randomization hypothesis.

Finally, we discuss how to compute the value function and approximate optimal randomized feedback controls from the Bellman equation according to value or policy iteration methods and by discretization of the state space and of the space of probability measures.

\vspace{1mm}

\noindent {\bf Outline of the paper.}  The rest of the paper is organized  as follows. Section \ref{secchaos} carefully formulates both the $N$-individual model  and the CMKV-MDP, and show their connection 
by providing the rate of  convergence of the latter to the limiting MFC when $N$ goes to infinity.  In Section \ref{seclift}, we establish the correspondence  of the CMKV-MDP with a lifted MDP on the space of probability measures 
with usual relaxed controls when there is no common noise or when the state space is finite. In the general case considered in Section \ref{secmainlift}, we show how to lift the CMKV-MDP  by a suitable enlargement of the action space in order to get the correspondance with a MDP on the Wasserstein space. We then derive the associated Bellman fixed point equation satisfied by the value function, and obtain the existence of approximate randomized feedback controls. We also highlight the differences between open-loop vs feedback vs randomized controls.  
We conclude in Section \ref{secconclu} by indicating some questions for future research. Finally, we collect in the Appendix some useful and technical results including measurable coupling arguments 
used in the proofs of the paper.

\vspace{2mm}

\noindent {\bf Notations.}  Given two measurable spaces $(\Xc_1,\Sigma_1)$ and $(\Xc_2,\Sigma_2)$, we denote by $\text{pr}_{_1}$ (resp. $\text{pr}_{_2}$)  the projection function $(x_1,x_2)$ $\in$  
$\Xc_1\times\Xc_2$  $\mapsto$ $x_1$ $\in$ $\Xc_1$ (resp. $x_2$ $\in$ $\Xc_2$). For a measurable function $\Phi$ $:$ $\Xc_1$ $\rightarrow$ $\Xc_2$, and a positive measure $\mu_1$ on $(\Xc_1,\Sigma_1)$, the pushforward measure  $\Phi\star\mu_1$ is the measure on $(\Xc_2,\Sigma_2)$ defined by 
\begin{align*}
\Phi\star\mu_1 (B_2) &= \mu_1\big( \Phi^{-1}(B_2) \big), \quad \forall B_2 \in \Sigma_2. 
\end{align*}
We denote by $\Pc(\Xc_1)$ the set of probability measures on $\Xc_1$, and $\Cc(\Xc_1)$ the cylinder (or weak) $\sigma$-algebra on $\Pc(\Xc_1)$, that is the smallest 
$\sigma$-algebra making all the functions $\mu$ $\in$ $\Pc(\Xc_1)$ $\mapsto$ $\mu(B_1)$ $\in$ $[0,1]$, measurable for all $B_1$ $\in$ $\Sigma_1$. 

A probability kernel $\nu$ on $\Xc_1\times\Xc_2$, denoted $\nu$ $\in$ $\hat\Xc_2(\Xc_1)$, 
is a measurable mapping from $(\Xc_1,\Sigma_1)$ into $(\Pc(\Xc_2),\Cc(\Xc_2))$, and we shall write indifferently $\nu(x_1,B_2)$ $=$ $\nu(x_1)(B_2)$, for all 
$x_1$ $\in$ $\Xc_1$, $B_2$ $\in$ $\Sigma_2$. 
Given a probability measure $\mu_1$ on $(\Xc_1,\Sigma_1)$, and a probability kernel $\nu$ $\in$ $\hat\Xc_2(\Xc_1)$, we denote by $\mu_1\cdot\nu$ the probability measure on 
$(\Xc_1\times\Xc_2,\Sigma_1\otimes\Sigma_2)$ defined by 
\begin{align}
(\mu_1\cdot\nu)( B_1\times B_2) &=  \int_{B_1\times B_2}  \mu_1(\d x_1) \nu(x_1,\d x_2), \quad \forall B_1 \in \Sigma_1, \; B_2 \in \Sigma_2. 
\end{align} 
 Let $X_1$ and $X_2$ be two random variables valued respectively on $\Xc_1$ and $\Xc_2$, denoted $X_i$ $\in$ $L^0(\Omega;\Xc_i)$. 
 We denote by $\Lc(X_i)$ the probability distribution of $X_i$, and by $\Lc(X_2|X_1)$ the conditional probability distribution of $X_2$ given $X_1$. With these notations, when $X_2$ $=$ $\Phi(X_1)$, 
 then $\Lc(X_2)$ $=$ $\Phi\star\Lc(X_1)$. 

When $(\Yc,d)$ is a compact metric space, the set $\Pc(\Yc)$ of probability measures on $\Yc$ is equipped with the Wasserstein distance
\begin{align}
\Wc(\mu,\mu') &= \inf\Big\{ \int_{\Yc^2} d(y,y') \boldsymbol{\mu}(\d y,\d y'):  \boldsymbol{\mu} \in \boldsymbol{\Pi}(\mu,\mu') \Big\},
\end{align}
where $\boldsymbol{\Pi}(\mu,\mu')$ is the set of probability measures on $\Yc\times\Yc$ with marginals $\mu$ and $\mu'$, i.e., $\text{pr}_{_1}\star\boldsymbol{\mu}$ $=$ $\mu$, and $\text{pr}_{_2}\star\boldsymbol{\mu}$ $=$ $\mu'$. 
Since $(\Yc,d)$ is compact, it is known (see e.g. Corollary  6.13 in \cite{Villani}) that the Borel $\sigma$-algebra generated by the Wasserstein metric coincides with the cylinder $\sigma$-algebra on $\Pc(\Yc)$, i.e., 
Wasserstein distances metrize weak convergence.  
We also recall the dual Kantorovich-Rubinstein representation of the Wasserstein distance
\begin{align}
\Wc(\mu,\mu') &=  \sup\Big\{ \int_\Yc \phi \; \mathrm{d}(\mu-\mu'):  \phi \in L_{lip}(\Yc;\R), [\phi]_{lip} \leq 1 \Big\}, 
\end{align}
where $L_{lip}(\Yc;\R)$ is the set of Lipschitz continuous functions $\phi$ from $\Yc$ into $\R$, and $[\phi]_{lip}$ $=$ $\sup\{ |\phi(y)-\phi(y')|/d(y,y'): y,y' \in \Yc, y\neq y'\}$.

\section{The $N$-agent and the limiting McKean-Vlasov MDP} \label{secchaos}

We formulate the mean-field Markov Decision Process (MDP) in a large population model with indistinguishable agents $i$ $\in$ $\N^*$ $=$ $\N\setminus\{0\}$.

Let $\Xc$ (the state space) and $A$ (the action space) be two compact Polish spaces equipped respectively with their metric $d$ and $d_A$.  We denote by $\Pc(\Xc)$ (resp. $\Pc(A)$) the 
space of probability measures on $\Xc$ (resp. $A$) equipped respectively with their Wasserstein distance  $\Wc$ and $\Wc_A$.  We also consider the product space $\Xc\times A$, equipped with the metric $\bd((x,a),(x',a'))$ $=$ $d(x,x')$ $+$ $d_A(a,a')$, $x,x'$ $\in$ $\Xc$, $a,a'$ $\in$ $A$,  and the associated space of probability measure $\Pc(\Xc\times A)$, equipped with its Wasserstein distance $\bW$. 
Let $G$, $E$, and $E^0$ be three measurable spaces, representing respectively  the initial information, idiosyncratic noise, and common noise spaces. 

We denote by $\Pi_{OL}$ the set of sequences $(\pi_t)_{t\in\N}$ (called {\it open-loop policies}) where $\pi_t$ is a measurable function from $G\times E^t\times (E^0)^t$ into $A$ for $t\in\N$.  

Let  $(\Omega, \Fc,\P)$ be a probability space on which are defined the following family of mutually i.i.d. random variables 
\begin{itemize}
\item $(\Gamma^i,\xi^i)_{i\in\N^\star}$ (initial informations and initial states) valued in  $G\times \Xc$
\item $(\eps^i_t)_{i\in\N^\star, t\in\N}$ (idiosyncratic noises) valued in $E$ with probability distribution $\lambda_\eps$ 
\item  $\eps^0:=(\eps^0_t)_{t\in\N}$ (common noise) valued in $E^0$. 
\end{itemize}
We assume that $\Fc$ contains an atomless random variable, i.e., $\Fc$ is rich enough, so that any probability measure $\nu$ on $\Xc$ (resp. $A$ or $\Xc\times A$) can be represented by the law of some 
random variable $Y$ on $\Xc$ (resp. $A$ or $\Xc\times A$), and we write $Y$ $\sim$ $\nu$, i.e., $\Lc(Y)$ $=$ $\nu$. 
Given an open-loop policy $\pi$, we associate an {\it open-loop control} for individual $i$ $\in$ $\N^*$  as the process $\alpha^{i,\pi}$ defined by   
\begin{align}
\alpha^{i,\pi} &= \pi_t(\Gamma^i,(\eps_s^i)_{s\leq t},(\eps_s^0)_{s\leq t}), \quad t \in \N. 
\end{align}
In other words, an open-loop control is a non-anticipative  process that depends on the initial information, the past idiosyncratic and common noises, but not on the states of the agent in contrast with closed-loop control.

Given $N\in\N^*$, and $\pi$ $\in$ $\Pi_{OL}$, the state process of  agent $i$ $=$ $1,\ldots,N$ in an $N$-agent MDP is given by the dynamical system
\begin{equation} \label{dynXi}
\left\{
\begin{array}{rcl}
	X^{i,N,\pi}_0 & = & \xi^i  \\
	X^{i,N,\pi}_{t+1} & = &  F(X^{i,N,\pi}_t, \alpha^{i,\pi}_t, \frac{1}{N}\sum_{j=1}^N\delta_{(X^{j,N,\pi}_t,\alpha^{j,\pi}_t)},\eps^i_{t+1}, \eps^0_{t+1}),  \quad 
	t\in\N,
\end{array}
\right.
\end{equation}
where $F$ is  a measurable function from $\Xc \times A \times \Pc(\Xc\times A)\times E\times E^0$ into $\Xc$, called state transition function.    
The $i$-th individual contribution to the influencer's gain over an infinite horizon  is defined  by
\begin{align}
J_i^{N,\pi} &:=   \sum_{t=0}^\infty \beta^tf\Big(X^{i,N,\pi}_t, \alpha^{i,\pi}_t, \frac{1}{N}\sum_{j=1}^N\delta_{(X^{j,N,\pi}_t,\alpha^{j,\pi}_t)}\Big), \quad i=1,\ldots,N, 
\end{align}
where the reward $f$ is a mesurable real-valued function on $\Xc \times A \times \Pc(\Xc\times A)$, assumed to be bounded (recall that $\Xc$ and $A$ are compact spaces), and $\beta$  is a positive  discount factor in $[0,1)$. 
The influencer's renormalized and expected  gains are 
\begin{align}
J^{N,\pi}  \; := \; \frac{1}{N}\sum_{i=1}^N J_i^{N,\pi}, & \quad \quad  V^{N,\pi} \; := \; \E \big[ J^{N,\pi}\big], 
\end{align} 
and the optimal value of the influencer is $V^N :=  \sup_{\pi \in \Pi_{OL}}V^{N,\pi}$.  Observe that the agents are indistinguishable in the sense that the initial pair of information/state $(\Gamma^i,\xi^i)_i$, and 
idiosyncratic noises are i.i.d., and the state transition function $F$, reward function $f$, and discount factor $\beta$ do not depend on $i$.

\vspace{1mm}

Let us now consider the asymptotic problem when the number of agents $N$ goes to infinity. In view of the propagation of chaos argument, we expect that the state process of agent $i$ $\in$ 
$\N^*$ in the infinite population model to be governed by the conditional  McKean-Vlasov dynamics
\begin{equation} \label{dynXiMKV}
\left\{
\begin{array}{rcl}
	X^{i,\pi}_0 & = & \xi^i  \\
	X^{i,\pi}_{t+1} & = &  F(X^{i,\pi}_t, \alpha^{i,\pi}_t, \P^0_{(X^{i,\pi}_t,\alpha^{i,\pi}_t)},\eps^i_{t+1}, \eps^0_{t+1}),  \quad 
	t\in\N. 
\end{array}
\right.
\end{equation}
Here, we denote by $\P^0$ and $\E^0$ the conditional probability and expectation knowing the common noise $\eps^0$, and then, given a random variable $Y$ valued in $\Yc$, we denote by 
$\P^0_Y$ or $\Lc^0(Y)$ its conditional law knowing $\eps^0$, which is a random variable valued in $\Pc(\Yc)$ (see Lemma \ref{condlaw}). 
The $i$-th individual contribution to the influencer's gain in the infinite population model is
\beqs
J_i^{\pi} &:=&  \sum_{t=0}^\infty \beta^tf\big(X^{i,\pi}_t, \alpha^{i,\pi}_t, \P^0_{(X^{i,\pi}_t,\alpha^{i,\pi}_t)}\big), \quad i \in \N^*, 
\enqs
and we define the conditional gain, expected gain, and optimal value, respectively by
\begin{align}
J^\pi \; := \;  \E^0 \big[ J_i^{\pi} \big] \; =  \; \E^0 \big[ J_1^{\pi} \big], &  \quad i \in \N^*  \quad (\mbox{by indistinguishability of the agents}), \nonumber  \\
V^\pi \; := \; \E\big[ J^\pi \big], & \quad\quad  V := \sup_{\pi \in \Pi_{OL}}V^\pi.  \label{defVMKV}
\end{align}
Problem \eqref{dynXiMKV}-\eqref{defVMKV} is called conditional McKean-Vlasov Markov decision process, CMKV-MDP in short.

 \vspace{1mm}

The main goal of this Section is to rigorously connect the finite-agent model to the infinite population model with mean-field interaction by proving the convergence of the $N$-agent MDP to the CMKV-MDP. 
First, we prove the almost sure convergence of the state processes under some continuity assumptions on the state transition function. 

\begin{Proposition}
\label{propaschaos} 
Assume that for all $(x_0,a,\nu_0,e^0)\in \Xc\times A\times\Pc(\Xc \times A)\times E^0$, the (random) function 
\begin{align*}
(x,\nu) \in \Xc\times\Pc(\Xc\times A) & \longmapsto \;  F(x,a,\nu,\eps^1_1, e^0) \in  \Xc  \quad \mbox{is continuous in } (x_0,\nu_0) \; a.s.
\end{align*}
i.e. $\P\big[\underset{(x,\nu)\rightarrow (x_0,\nu_0)}{\lim}F(x,a,\nu,\eps^1_1, e^0)=F(x_0,a,\nu_0,\eps^1_1, e^0)\big]$ $=$ $1$.
Then, for any $\pi\in\Pi_{OL}$, $X^{i,N,\pi}_t \underset{N\rightarrow \infty}{\overset{a.s.}{\rightarrow}}X^{i,\pi}_t$, $i$ $\in$ $\N^*$, $t$ $\in$ $\N$, and 
\beqs
\bW\Big(\frac{1}{N}\sum_{i=1}^N \delta_{(X^{i,N,\pi}_t,\alpha_t^{i,\pi})}, \P^0_{(X^{i,\pi}_t,\alpha_t^{i,\pi})}\Big) \underset{N\rightarrow \infty}{\overset{a.s.}{\longrightarrow}} 0, \; \quad \frac{1}{N}\sum_{i=1}^N d(X^{i,N,\pi}_t,X^{i,\pi}_t)\underset{N\rightarrow \infty}{\overset{}{\rightarrow}} 0\quad a.s.
\enqs
Furthermore, if we assume  that  for all $a\in A$, the function $(x,\nu) \in \Xc\times\Pc(\Xc\times A)$ $\mapsto$  
$f(x,a,\nu)$ is continuous, then  
\beqs
J_i^{N,\pi} \underset{N\rightarrow \infty}{\overset{a.s.}{\longrightarrow}} J_i^\pi, \quad \; J^{N,\pi}  \underset{N\rightarrow \infty}{\overset{a.s.}{\longrightarrow}} J^\pi, \quad \;\; 
V^{N,\pi} \underset{N\rightarrow \infty}{\overset{}{\longrightarrow}} V^\pi, \quad \; \text{and} \; \quad\underset{N\rightarrow \infty}{\liminf} V^N  \; \geq \;  V.
\enqs
\end{Proposition}
{\begin{Remark}\label{remregexp}
\rm{
We stress the importance of making the continuity assumption for $F$ only a.s. w.r.t. the noise $\eps^1_1$. Indeed, most of the practical applications we have in mind concerns finite spaces $\Xc$, for which $F$ cannot be, strictly speaking, 
continuous, since it maps from a continuous space $\Pc(\Xc\times A)$ to a finite space $\Xc$. However, in this case, the discontinuities will essentially be placed on the boundaries between each region where $F$ is constant. Our assumption then simply formalizes the idea that these boundaries depend upon $\eps^1_1$ in a way that no point $(x_0,a,\nu_0,e^0)\in \Xc\times A\times\Pc(\Xc \times A)\times E^0$ has a positive probability (w.r.t. $\eps^1_1$) to be on a boundary. In other words, we assume that the discontinuity sets always vary with $\eps^1_1$, which is a natural phenomenon. 
}
\ep
\end{Remark}
{\bf Proof.}  
Fix $\pi$ $\in$ $\Pi_{OL}$. We omit the dependence  of the state processes and  control on $\pi$, and denote by 
$\nu^{N}_t$ $:=$ $\frac{1}{N}\sum_{i=1}^N \delta_{(X^{i,N}_t,\alpha^i_t)}$,  $\nu^{N,\infty}_t$ $:=$ 
$\frac{1}{N}\sum_{i=1}^N \delta_{(X^{i}_t,\alpha^i_t)}$,  and $\nu_t$ $:=$ $\P^0_{(X^i_t,\alpha^i_t)}$. 

\vspace{1mm}

\noindent (1) We first  prove the convergence of trajectories, for all $i$ $\in$ $\N^\star$, 
\beqs
\P \big[\underset{N\rightarrow\infty}{\text{lim}}(X^{i,N}_t,\nu^N_t)=(X^i_t,\nu_t)\big] = 1, &  & 
\P\Big[  \underset{N\rightarrow\infty}{\text{lim}} \frac{1}{N}\sum_{j=1}^N d(X^{j,N}_t,X^{j}_t)  =  0 \Big] = 1,
\enqs
by induction on $t\in\N$.

\vspace{1mm}

\noindent{\it  - Initialization.}  For $t=0$, we have $X^{i,N}_0=\xi^i=X^{i}_0$, $\alpha_0^i$ $=$ $\pi_0(\Gamma_i)$, 
for all $N\in\N_\star$ and $i\in\N^\star$, which obviously implies that   $X^{i,N}_0 \underset{N\rightarrow \infty}{\overset{a.s.}{\rightarrow}}X^{i}_0$ and 
$\frac{1}{N} \sum_{j=1}^N d(X^{j,N}_0,X^{j}_0) \underset{N\rightarrow \infty}{\overset{}{\rightarrow}} 0$. 
Moreover, $\Wc(\nu_0^N,\nu_0)$ $=$ $\Wc( \nu^{N,\infty}_0, \nu_0)$, which converges to zero a.s., when $N$ goes to infinity, by the 
weak convergence of empirical measures (see \cite{vanvaart}), and the fact that Wasserstein distance metrizes weak convergence. 
 
\vspace{1mm}

\noindent{\it - Hereditary Property.}  We have 
\begin{equation}\label{xni}
X^{i,N}_{t+1} = F(X^{i,N}_t, \alpha^i_t, \nu^N_t, \eps^i_{t+1},\eps^0_{t+1})	\quad \text{and}\quad X^{i}_{t+1} = F(X^{i}_t,  \alpha^i_t, \nu_t, \eps^i_{t+1},\eps^0_{t+1}).
\end{equation}
By a simple conditioning, we notice that $\P\big[\underset{N\rightarrow\infty}{\text{lim}}X^{i,N}_{t+1}=X^i_{t+1}\big]$ $=$ 
$\E\big[P\big((X^{i,N}_t,\nu_t^N)_N,X_t^i,\nu_t,\alpha_t^i\big)\big]$, where  
\beqs
P\big( (x_N,\nu_N)_N,x,\nu,a) \big) &=&  \P\Big[ \underset{N\rightarrow \infty}{\text{lim}}F(x_N, a,\nu_N, \eps^i_{t+1},\eps^0_{t+1})=F(x, a,\nu, \eps^i_{t+1},\eps^0_{t+1})\Big].
\enqs
By the continuity assumption on $F$, we see that $P\big( (x_N,\nu_N)_N,x,\nu,a) \big)$ is bounded from below by  
${\bf 1}_{\underset{N\rightarrow\infty}{\text{lim}}(x_N,\nu_N)=(x,\nu)}$, and thus
\begin{align}\label{eqprop-as}
\P\big[ \underset{N\rightarrow\infty}{\text{lim}}X^{i,N}_{t+1}=X^i_{t+1} \big] &\geq \;  \P \big[\underset{N\rightarrow\infty}{\text{lim}}(X^{i,N}_t,\nu^N_t)=(X^i_t,\nu_t)\big]. 
\end{align}
This proves by induction hypothesis that  $\P\big[ \underset{N\rightarrow\infty}{\text{lim}}X^{i,N}_{t+1}=X^i_{t+1} \big]$ $=$ $1$.  

Let us now show that ${\frac{1}{N}\sum_{i=1}^N d(X^{i,N}_{t+1}, X^{i}_{t+1}) \underset{N\rightarrow \infty}{\overset{a.s.}{\rightarrow}} 0}$. From \eqref{xni}, we have
\begin{eqnarray}
d(X^{i,N}_{t+1}, X^{ i}_{t+1}) & \leq & \Delta_{D_N} F(X^{i}_t,\alpha^i_t,\nu_t,\eps^i_{t+1}, \eps^0_{t+1}) \label{eq_theorem_chaos_propagation_1}
\end{eqnarray}
with  $D_N:= \max[d(X^{i,N}_t,X^{i}_t), \bW(\nu^N_t, \nu_t)]$, and 
$$
\Delta_y F(x,a,\nu,e,e^0):=\sup_{(x',\nu')\in D} \{d(F(x',a,\nu',e,e^0),F(x,a,\nu,e,e^0)){\bf 1}_{\max[d(x',x), \bW(\nu',\nu)]\leq y}\},
$$
where $D$ is a fixed countable dense set of the separable space $\Xc\times\Pc(\Xc\times A)$, which implies that $(y,x, a,\nu,e,e^0 )\mapsto \Delta_y F(x,a,\nu,e,e^0) $ is a measurable function. Fix $\epsilon>0$. 
Let $\Delta_\Xc$ denote the diameter of the compact metric space $\Xc$. We thus have, for any $\eta >0$, 
\begin{eqnarray}
d(X^{i,N}_{t+1}, X^{i}_{t+1}) & \leq & d(X^{i,N}_{t+1}, X^{i}_{t+1}){\bf 1}_{D_N\leq \eta}+ d(X^{i,N}_{t+1}, X^{i}_{t+1}){\bf 1}_{D_N> \eta}\\
& \leq & \Delta_\eta F(X^i_t, \alpha^i_t, \nu_t, \eps^1_{t+1},\eps^0_{t+1})+ \Delta_\Xc{\bf 1}_{d(X^{N,i}_t,X^{i}_t)> \eta} 
+ \Delta_{\Xc}{\bf 1}_{\bW(\nu^N_t,\nu_t)> \eta},  \label{eq_theorem_chaos_propagation_11}
\end{eqnarray}
and thus  
\begin{eqnarray}
& & \frac{1}{N}\sum_{i=1}^N d(X^{i,N}_{t+1}, X^{i}_{t+1}) \\
& \leq & 
\frac{1}{N}\sum_{i=1}^N\Delta_\eta F(X^i_t, \alpha^i_t, \nu_t, \eps^1_{t+1},\eps^0_{t+1})+ \frac{\Delta_\Xc}{\eta N}\sum_{i=1}^N d(X^{i,N}_t,X^{i}_t)+ 
\Delta_\Xc{\bf 1}_{\bW(\nu^N_t,\nu_t) > \eta}.  \label{eq_theorem_chaos_propagation_12}
\end{eqnarray}
The second and third terms in the right hand side converge to $0$ by induction hypothesis, and by Proposition \ref{annexe_conditional_wassertein_convergence}, the first term converges to 
\beqs
\E^0 \big[ \Delta_\eta F(X^1_t,\alpha^1_t, \nu_t, \eps^{1}_{t+1},\eps^0_{t+1}) \big] &=&  \E^0 \left[\E\left[\Delta_\eta F(x^1_t,a^1_t, \nu_t, \eps^{1}_{t+1},e^0_{t+1})\right]_{(x^1_t, a^1_t,e^0_{t+1}):=(X^1_t,\alpha^1_t, \eps^0_{t+1})}\right]. 
\enqs
As $\eta\rightarrow 0$, the inner expectation tends to zero by continuity assumption on $F$ and by dominated convergence. Then, the outer expectation converges to zero by conditional dominated convergence, and will thus be smaller than $\frac{\epsilon}{2}$ for $\eta$ small enough, which implies that 
$\frac{1}{N}\sum_{i=1}^N d(X^{i,N}_{t+1}, X^{i}_{t+1}) $ will be smaller than $\epsilon$ for $N$ large enough. 

Let us finally prove that $\bW(\nu^N_{t+1},\nu_{t+1})\overset{a.s.}{\underset{N\rightarrow \infty}{\rightarrow}}0$. We have $\bW(\nu^N_{t+1},\nu_{t+1})$ $\leq$ 
$\bW(\nu^N_{t+1},\nu^{N,\infty}_{t+1})$ $+$ $\bW(\nu^{N,\infty}_{t+1},\nu_{t+1})$. To dominate the first term $\bW(\nu^N_{t+1},\nu^{N,\infty}_{t+1})$, notice that, given a variable $U_N\sim \Uc(\{ 1,...,N\})$, the random measure 
$\nu^N_{t+1}$ (resp. $\nu^{N,\infty}_{t+1}$) is, at fixed $\omega\in\Omega$, the law of the pair of random variable ($X^{U_N,N}_{t+1}(\omega),\alpha_{t+1}^{U_N}(\omega))$ 
(resp. $(X^{U_N}_{t+1}(\omega),\alpha_{t+1}^{U_N}(\omega)$) where we stress that only $U_N$ is random here, essentially selecting each sample of these empirical measures with probability $\frac{1}{N}$. Thus, by definiton of the Wasserstein distance, $\bW(\nu^N_{t+1},\nu^{N,\infty}_{t+1}) $ is dominated by $\E[d(X^{U_N,N}_{t+1}(\omega), X^{U_N}_{t+1}(\omega))]$ $=$ 
$\frac{1}{N}\sum_{i=1}^N d(X^{i,N}_{t+1}(\omega), X^{i}_{t+1} (\omega))$,  which has been proved to converge to zero. For the second term, observe that 
$\alpha^{i}_{t+1}$ $=$ $\pi_{t+1}(\Gamma^i, (\eps^i_s)_{s\leq t+1}, (\eps^0_s)_{s\leq t+1})$, and by  Proposition \ref{condlaw-adapted},  
there exists a measurable function $f_{t+1}:\Xc\times G\times E^{t+1} \times (E^0)^{t+1}\rightarrow \Xc$ such that 
$X^{i,N}_{t+1}$ $=$ $f_{t+1}(\Gamma^i,(\eps^i_s)_{s\leq t+1}, (\eps^0_s)_{s\leq t+1}))$ .  From  Proposition \ref{annexe_conditional_wassertein_convergence}, we then deduce that $\bW(\nu^{N,\infty}_{t+1},\nu_{t+1})$ converges to zero as $N$ goes to infinity. This concludes the induction. 

\vspace{2mm}

\noindent (2) Let us now study the convergence of gains. By the continuity assumption on $f$, we have 
$f(X^{i,N}_t,\alpha^i_t,\nu^N_t) \overset{a.s.}{\underset{N\rightarrow \infty}{\rightarrow}} f(X^{i}_t,\alpha^i_t,\nu_t) $ for all $t$ $\in$ $\N$. Thus, as $f$ is bounded, we get  by dominated convergence that $J_i^{N,\pi} \overset{a.s.}{\underset{N\rightarrow \infty}{\rightarrow}} J_i^\pi$. Let us now study the convergence of $J^{N,\pi}$ to $J^\pi$. We write 
\beqs
\vert J^{N,\pi}-J^\pi \vert &\leq & \frac{1}{N}\sum_{i=1}^N \vert J^{N,\pi}_i- J_i^\pi \vert +\Big\vert \frac{1}{N}\sum_{i=1}^N J_i^\pi - J^\pi \Big\vert \: =: \; 
S_N^1 + S_N^2. 
\enqs
The second term $S_N^2$ converges  a.s. to zero by Propositions \ref{condlaw-adapted} and 
\ref{annexe_conditional_wassertein_convergence}, as $N$ goes to infinity. On the other hand,  
\beqs
S_N^1 & \leq & \sum_{t=0}^\infty \beta^t \Delta_N(f), \quad \mbox{ with }  \; 
\Delta_N(f) := \frac{1}{N}\sum_{i=1}^N \big\vert f(X^{i,N}_t,\alpha^i_t, \nu^N_t)-f(X^i_t,\alpha^i_t, \nu_t) \big\vert. 
\enqs
By the same argument as above in (1) for showing  that  ${\frac{1}{N}\sum_{i=1}^N d(X^{i,N}_{t+1}, X^{i}_{t+1}) \underset{N\rightarrow \infty}{\overset{a.s.}{\rightarrow}} 0}$, we 
prove that $\Delta_N(f)$  tends a.s.  to zero  as $N\rightarrow \infty$. Since $f$ is bounded, we deduce by the dominated convergence theorem that $S_N^1$ converges a.s. to zero 
as $N$ goes to infinity, and thus $J^{N,\pi} \overset{a.s.}{\underset{N\rightarrow \infty}{\rightarrow}} J^\pi$. By dominated convergence, we then also obtain that 
$V^{N,\pi}$ $=$  $\E[J^N] \overset{a.s.}{\underset{N\rightarrow \infty}{\rightarrow}} \E[J^\pi] =V^\pi$.   
Finally, by considering an $\epsilon$-optimal policy $\pi_\epsilon$ for $V$, we have  
\beqs
\underset{N\rightarrow \infty}{\liminf} V^N \geq  \underset{N\rightarrow \infty}{\text{lim}} V^{N,\pi_\epsilon} \; = \;  V^{\pi_\epsilon} \;  \geq \;  V-\epsilon,
\enqs
which implies, by sending $\epsilon$ to zero, that $\underset{N\rightarrow \infty}{\text{liminf}} \; V^N \geq V$.
\ep

\vspace{3mm}

Next, under Lipschitz assumptions on the state transition and reward functions, we prove the corresponding convergence in $\mathbb{L}^1$, which implies the convergence of the optimal value,  
and also a rate of convergence in terms of the rate of convergence in Wasserstein distance of the empirical measure.

\vspace{3mm} 

\noindent $({\bf HF_{lip}})$ There exists $K_F$ $>$ $0$,  such that   for all $a\in A$, $e^0$ $\in$ $E^0$, $x,x'\in\Xc$, $\nu,\nu'\in\Pc(\Xc\times A)$,
\begin{align*}
\E \big[ d\big(F(x,a,\nu,\eps^1_1, e^0), F(x',a,\nu',\eps^1_1, e^0)\big)\big] & \leq  \; K_F \big( d(x,x') + \bW(\nu,\nu')  \big)). 
\end{align*}

\vspace{1mm} 

\noindent $({\bf Hf_{lip}})$ There exists $K_f$ $>$ $0$,  such that   for all $a\in A$,  $x,x'\in\Xc$, $\nu,\nu'\in\Pc(\Xc\times A)$,
\begin{align*}
d(f(x,a,\nu), f(x',a,\nu'))  & \leq  \; K_f  \big( d(x,x') + \bW(\nu,\nu')  \big)). 
\end{align*}

\vspace{2mm}

{\begin{Remark}
\rm{
Here again, we stress the importance of making the regularity assumptions for $F$ in {\em expectation} only. For the same argument as in Remark \ref{remregexp}, when $\Xc$ is finite, $F$ cannot be, strictly speaking, 
Lipschitz. However, $F$ can be Lipschitz {\em in expectation}, e.g. once integrated w.r.t. the idiosyncratic noise, and, again, it is a very natural assumption.} 
\ep
\end{Remark}

\vspace{3mm}

In the sequel, we shall denote by $\Delta_\Xc$ the diameter of the metric space $\Xc$, and define 
\begin{align} \label{defMN}
M_N &:= \underset{\nu\in\Pc(\Xc\times A)}{\sup}\E [\bW(\nu_N,\nu)],
\end{align}
 where $\nu_N$ is the empirical measure $\nu_N$ $=$ $\frac{1}{N}\sum_{n=1}^N\delta_{Y_n}$,  $(Y_n)_{1\leq n\leq N}$ are i.i.d. random variables with law $\nu$. 
 We recall in the next Lemma recent results about non asymptotic bounds for  the mean rate of convergence in Wasserstein distance of the empirical measure. 

\vspace{3mm}

\begin{Lemma}\label{lemWass}
We have $M_N\underset{N\rightarrow \infty}{\rightarrow}0$. Furthermore,
\begin{itemize}
\item If $\Xc\times A\subset \R^d$ for some $d\in\N^\star$, then:
$M_N=\Oc(N^{-\frac{1}{2}})$	for $d=1$, $M_N=\Oc(N^{-\frac{1}{2}}\log (1+N))$	for $d=2$, and $M_N=\Oc(N^{-\frac{1}{d}})$	for $d\geq 3$.
\item If for all $\delta>0$, the smallest number of balls with radius $\delta$ covering the compact metric set $\Xc\times A$ with diameter $\Delta_{\Xc\times A}$  
is smaller than $\Oc\Big( \big(\frac{\Delta_{\Xc\times A}}{\delta}\big)^\theta\Big)$ for $\theta>2$, then $M_N= \Oc(N^{-1/\theta})$.
\end{itemize}
\end{Lemma}
{\bf Proof.} The second point is proved in \cite{fougui15}, and the third one in \cite{boileg14}. 
\ep

 \begin{Remark} \label{remdiscret}
 {\rm In the case where the state and action spaces $\Xc$ and $A$ are finite, let 
 $\phi:\Xc\times A\xhookrightarrow{} \R$ be any injective function. Then $\phi^{-1}$ is necessarily Lipschitz. From the dual Kantorovich representation of the Wasserstein distance, we have
 \beqs
 M_N&\leq& \underset{\nu\in\Pc(\Xc\times A)}{\sup} \E \left[\underset{g\in L_{lip}(\Xc\times A;\R),[g]_{lip}\leq 1}{\sup}\int_{\Xc\times A} g \text{d}(\nu_N-\nu)\right]\\
 &\leq& \underset{\nu\in\Pc(\Xc\times A)}{\sup} \E \left[\underset{g\in L_{lip}(\Xc\times A;\R),[g]_{lip}\leq 1}{\sup}\int_{\Xc\times A} g\circ \phi^{-1} \text{d}\big((\phi\star\nu)_N-\phi\star\nu)\big)\right]\\
  &\leq& C\underset{\nu\in\Pc(\Xc\times A)}{\sup} \E \left[\Wc((\phi\star\nu)_N,\phi\star\nu)\right]\leq C\underset{\nu\in\Pc(\phi(\Xc\times A))}{\sup} \E \left[\Wc(\nu_N,\nu)\right],
 \enqs
 where $C$ is the Lipschitz constant of $\phi^{-1}$. Thus, by the second point in Lemma \ref{lemWass}, $M_N=\Oc(N^{-\frac{1}{2}})$. 
 }
 \ep
 \end{Remark}

\begin{Theorem} \label{theoL1chaos} 
\label{theolaw} 
Assume $({\bf HF_{lip}})$. For all $i$ $\in$ $\N^*$,  $t\in\N$, 
\begin{align}
\underset{\pi\in\Pi_{OL}}{\sup}\E \big[ d(X^{i,N,\pi}_t, X^{i,\pi}_t) \big] & = \;  \Oc(M_N), \label{convX} \\
\underset{\pi\in\Pi_{OL}}{\sup} \E \Big[ \bW\Big(\frac{1}{N}\sum_{i=1}^N \delta_{(X^{i,N,\pi}_t,\alpha^{i,\pi})}, \P^0_{(X^{i,\pi}_t,\alpha^{i,\pi})}\Big)  \Big] & = \; \Oc(M_N).  \label{convnu}
\end{align}
Furthermore, if we assume $({\bf Hf_{lip}})$, and set  $\gamma$ $=$ $\min\big[1,\frac{\vert\ln\beta\vert}{(\ln 2K_F)_+}\big]$,  there exists a constant $C$ $=$ 
$C(K_F,K_f,\beta, \gamma)$ (explicit in the proof) such that for $N$ large enough, 
\beq \label{rateV}
\underset{\pi\in\Pi_{OL}}{\sup}\vert V^{N,\pi}-V^\pi\vert & \leq & C M_N^\gamma,
\enq
and thus $\vert V^{N}-V\vert =\Oc(M_N^\gamma) $. Consequently, any $\eps-optimal$ policy for the CMKV-MDP is $\Oc(\eps)$-optimal for the $N$-agent MDP problem,  
and conversely, any $\eps-optimal$ policy for the $N$-agent MDP problem is $\Oc(\eps)$-optimal for the CMKV-MDP,  for $N$ large enough, namely 
$M_N^\gamma$ $=$ $\Oc(\epsilon)$. 
\end{Theorem}
{\bf Proof.}  
Given $\pi\in\Pi_{OL}$,  denote by  $\nu^{N,\pi}_t:= \frac{1}{N}\sum_{i=1}^N \delta_{(X^{i,N,\pi}_t,\alpha^{i,\pi}_t)}$, 
${\nu^{N,\infty,\pi}_t := \frac{1}{N}\sum_{i=1}^N \delta_{(X^{i,\pi}_t,\alpha^{i,\pi}_t)}} $ and ${\nu_t^\pi:=\P^0_{(X^{i,\pi}_t,\alpha^i_t)}}$.  
By definition, $\alpha^{i,\pi}_t=\pi_t(\Gamma^i, (\eps^i_s)_{s\leq t}, (\eps^0_s)_{s\leq t})$, and by Lemma \ref{condlaw-adapted}, 
we have $X^{i,\pi}_t$ $=$ $f^\pi_t(\xi^i,\Gamma^i, (\eps^i_s)_{s\leq t},(\eps^0_s)_{s\leq t})$ for some measurable function ${f^\pi_t\in L^0(\Xc\times G\times E^t\times (E^0)^t, \Xc)}$. 
By definition of $M_N$ in \eqref{defMN}, we have 
\begin{align} \label{inegnu} 
\E \Big[ \bW(\nu_t^{N,\infty,\pi},\nu_t^\pi) \Big] & \leq \; M_N, \quad \forall N \in \N, \;  \forall \pi \in \Pi_{OL}. 
\end{align}
Let us now prove \eqref{convX} by induction on $t$ $\in$ $\N$.  At $t=0$, $X^{i,N,\pi}_0=X^{i,\pi}_0$ $=$ $\xi^i$,  and the result is obvious. Now assume that it holds true at time $t\in\N$ and let us show that 
it then holds  true at time $t+1$. 
By a simple conditioning argument, $\E [d(X^{i,N,\pi}_{t+1},X^{i,\pi}_{t+1})] $ $=$ $\E\big[ \Delta\big( X^{i,N,\pi}_{t}, X^{i,\pi}_{t},\alpha^i_t,\nu^{N,\pi}_t,\nu_t^\pi,\eps^0_{t+1} \big)\big]$,
where
\begin{align}
\Delta(x,x',a,\nu,\nu',e^0) &= \;  \E [d(F(x,a,\nu, \eps^i_{t+1},e^0),F(x',a, \nu', \eps^i_{t+1},e^0))] \nonumber \\
& \leq \;   K_F \big( d(x,x') +  \bW(\nu,\nu') \big), \label{Delta}
\end{align}
by $({\bf HF_{lip}})$. On the other hand, we have 
\begin{align}
\E \big[ \bW (\nu^{N,\pi}_t,\nu_t^\pi) \big] & \leq \;  \E\big[ \bW(\nu^{N,\pi}_t,\nu^{N,\infty,\pi}_t)\big] + \E\big[\bW(\nu^{N,\infty,\pi}_t,\nu_t^\pi)\big] \nonumber \\
& \leq \; \E [d(X^{i,N,\pi}_{t},X^{i,\pi}_{t})]  +  M_N,  \label{inegnunun}
\end{align}
where we used the fact that  $\bW(\nu^{N,\pi}_{t},\nu^{N,\infty,\pi}_{t})$ $\leq$ $\frac{1}{N}\sum_{i=1}^N d(X^{i,N,\pi}_{t}, X^{i,\pi}_{t})$, and \eqref{inegnu}. 
It follows from \eqref{Delta}  that 
\begin{align}
\E \big[ d(X^{i,N,\pi}_{t+1},X^{i,\pi}_{t+1})\big] 
& \leq \;  K_F \Big(  2 \E [d(X^{i,N,\pi}_{t},X^{i,\pi}_{t})]  +  M_N \Big), \quad \forall \pi \in \Pi_{OL},  \label{estid} 
\end{align}
 which proves that ${\underset{\pi\in\Pi_{OL}}{\sup}\E [d(X^{i,N,\pi}_{t+1},X^{i,\pi}_{t+1})] }$ $=$ $\Oc(M_N)$ by induction hypothesis, and thus   \eqref{convX}. 
 Plugging \eqref{convX} into \eqref{inegnunun} then yields \eqref{convnu}.

\vspace{2mm}

Let us now prove the convergence of gains. From $({\bf Hf_{lip}})$, and \eqref{inegnunun},  we have
\begin{align}
 \vert V^{N,\pi}-V^\pi\vert 
& \leq \;  K_f  \sum_{t=0}^\infty \beta^t  \E \Big[ d(X^{i,N,\pi}_t,X^{i,\pi}_t) +  \bW(\nu^{N,\pi}_t,\nu_t^\pi) \Big] \nonumber \\
& \leq \;  K_f \Big( 2  \sum_{t=0}^\infty \beta^t \delta_t^N + \frac{M_N}{1-\beta} \Big), \quad \forall \pi \in \Pi_{OL}, \label{eqrate1}
\end{align}
where we set $\delta_t^N$ $:=$ ${\underset{\pi\in\Pi_{OL}}{\sup}\E [d(X^{i,N,\pi}_{t},X^{i,\pi}_{t})] }$.  From \eqref{estid}, we have: 
$\delta_{t+1}^N$ $\leq$ $2K_F \delta_t^N + K_F M_N$, $t$ $\in$ $\N$, with $\delta_0^N$ $=$ $0$, and so by induction:
\begin{align}
\delta_t^N & \leq \; \frac{K_F}{1-2K_F} M_N  + s_t\big(\frac{K_F}{|2K_F-1|} M_N\big), \;\;\;  s_t(m) \; := \; m (2K_F)^t, \; m \geq 0. 
\end{align}
 where we may assume w.l.o.g. that $2K_F$ $\neq$ $1$.  Observing that we obviously have  $\delta_t^N$ $\leq$ $\Delta_\Xc$ (the diameter of $\Xc$), we deduce that 
\begin{align}
\sum_{t=0}^\infty \beta^t \delta_t^N & \leq \; \frac{K_F}{(1-2K_F)(1-\beta)} M_N  +  S\big(\frac{K_F}{|2K_F-1|} M_N \big) \label{interdelta} \\
S(m) & := \; \sum_{t=0}^\infty \beta^t   \min \big[ s_t(m);\Delta_\Xc \big], \quad m \geq 0.   
\end{align}
If $2\beta K_F<1$, we clearly have
\beqs
S(m)\leq \sum_{t=0}^\infty (m\beta 2K_F)^t=\frac{m}{1-\beta 2K_F}, 
\enqs
and so by \eqref{eqrate1}, $\vert V^{N,\pi}-V^\pi\vert=\Oc(M_N)$. Let us now study the case $2\beta K_F>1$. In this case, in particular, $2K_F>1$, thus  
$t$ $\mapsto$ $s_t(m)$ is nondecreasing, and so
\begin{align}
S(m) &\leq\; \sum_{t=0}^\infty \int_t^{t+1}\beta^t   \min \big[ s_t(m);\Delta_\Xc \big]ds \\ 
&\leq\; \frac{1}{\beta}\sum_{t=0}^\infty \int_t^{t+1}\beta^s   \min \big[ m (2K_F)^s;\Delta_\Xc \big]ds \\
&\leq\; \frac{1}{\beta} \int_0^{\infty}e^{-\vert\ln\beta\vert s}   \min \big[ m e^{\ln(2K_F)s};\Delta_\Xc \big]ds. \label{serieint}
\end{align}
Let $t_\star$ be such that $me^{\ln (2K_F)t_\star}=\Delta_\Xc$, i.e.
$t_\star$ $=$ $\frac{\ln(\Delta_\Xc/m)}{\ln(2 K_F)}$. 
Then, 
\beqs
\int_0^{\infty}e^{-\vert\ln\beta\vert s}   \min \big[ m e^{\ln(2K_F)s};\Delta_\Xc \big]ds &\leq& m\int_0^{t^\star} e^{\ln (2K_F \beta)s} ds+\Delta_\Xc\int_{t^\star}^\infty e^{\ln (\beta) s}  ds\\
&\leq& \frac{m}{\ln (2K_F \beta)}\left[ e^{\ln (2K_F \beta) t^\star}-1 \right]-\frac{\Delta_\Xc}{\ln \beta}e^{\ln (\beta)  t^\star}.
\enqs
After substituting $t_\star$ by its explicit value, we then obtain
\beqs
\int_0^{\infty}e^{-\vert\ln\beta\vert s}   \min \big[ m e^{\ln(2K_F)s};\Delta_\Xc \big]ds &\leq& \frac{m}{\ln (2K_F \beta)}\Big[ \Big(\frac{\Delta_\Xc}{m}\Big)^{\frac{\ln (2K_F\beta)}{\ln(2 K_F)}} -1 \Big]-\frac{\Delta_\Xc}{\ln \beta}\left(\frac{\Delta_\Xc}{m}\right)^{ \frac{\ln (\beta)}{\ln(2 K_F)}}\\
&\leq& \Delta_\Xc\Big(\frac{1}{\ln (2K_F \beta)}-\frac{1}{\ln \beta} \Big)\Big(\frac{\Delta_\Xc}{m}\Big)^{ \frac{\ln (\beta)}{\ln(2 K_F)}} -\frac{m}{\ln (2K_F \beta)}\\
&\leq& \Oc\Big(m^{\min\big[1,\frac{\vert \ln \beta\vert}{\ln (2 K_F)}\big]}\Big)=\Oc\Big(M_N^{\min\big[1,\frac{\vert \ln \beta\vert}{\ln (2K_F)}\big]}\Big). 
\enqs
This, combined with \eqref{eqrate1}, \eqref{interdelta} and \eqref{serieint}, concludes the proof.
\ep

\begin{Remark}
{\rm If  the Lipschitz constant  in $({\bf HF_{lip}})$ satisfies $\beta 2K_F$ $<$ $1$, then we can take $\gamma$ $=$ $1$ in the rate of convergence \eqref{rateV} of the optimal value.

In the particular case when $F$ and $f$ depend on the joint distribution $\nu$ $\in$ $\Pc(\Xc\times A)$ only through its marginals on $\Pc(\Xc)$ and $\Pc(A)$, which is 
the usual framework considered in controlled mean-field dynamics, then a careful look in the above proof shows that  the rate of convergence of the CMKV-MDP will be expressed in terms of 
\begin{align} \label{defMN2}
\tilde M_N &:=  \max \Big\{  \underset{\mu\in\Pc(\Xc)}{\sup}\E [\Wc(\mu_N,\mu)],  \underset{\upsilon\in\Pc(A)}{\sup}\E [\Wc_A(\upsilon_N,\upsilon)] \Big\}, 
\end{align}
instead of $M_N$ in \eqref{defMN}, where here $\mu_N$ (resp. $\upsilon_N$) is the empirical measure associated to $\mu$ (resp. $\upsilon$) $\in$ $\Pc(\Xc)$ (resp. $\Pc(A)$). From Lemma \ref{lemWass}, the speed of convergence of $\tilde M_N$ is faster 
than the one of $M_N$. For instance when $\Xc$ $\subset$ $\R$, $A$ $\subset$ $\R$, then $\tilde M_N$ $=$ $\Oc(N^{-1/2})$, while $M_N$ $=$ $\Oc\big(N^{-1/2}\log(1+N)\big)$.

The propagation of chaos in mean-field control problems has been addressed in several papers, see \cite{lac17}, \cite{cardelbook}, \cite{lautan20}, \cite{dje20}.  
For instance, the recent paper \cite{dje20} dealing with continuous time mean-field diffusion processes, shows that the propagation of chaos holds when the optimization is performed over asymmetric open-loop controls where each individual control is allowed to depend upon all the past noises (including the noises from other individuals), in contrast to our model where the optimization is made over 
symmetric ``individualized'' open-loop controls (each control $\alpha^i$ only depends upon the $i$-th individual's noise and the common noise). However,  rate of convergence is not provided in \cite{dje20}. To the best of our knowledge, convergence rate has been obtained only in the continuous time case with feedback controls and finite horizon in Chapter 6, vol 2 of \cite{cardelbook}. The originality of Theorem \ref{theoL1chaos}, besides providing a convergence rate in the individualized open-loop case, is that we are dealing with an infinite horizon. The above proof illustrates how two phenomenons are involved in this result: for $t\leq t_\star$, it is the regularity of $F$ which ensures the closeness of the $N$-agent processes and the McKean-Vlasov processes, and thus the closeness of their rewards.  As $t$ becomes larger, because the horizon is infinite, the $N$-agent and McKean-Vlasov processes drift away from each other due to errors accumulation, but the compactness of $\Xc$ prevents them from getting too far away from each other, and for $t>t_\star$, the discounting factor $\beta
^t$, by being smaller and smaller, makes the discounted rewards geometrically tend to $0$ and be neglectable. What is remarkable is that these combined phenomenons work together well enough to still ensure a fast convergence, even in infinite horizon, in $\Oc(M_N^\gamma)$, with an explicit order $\gamma\leq 1$.} 
\ep
\end{Remark}

\section{Lifted MDP on $\Pc(\Xc)$} \label{seclift}

Theorem \ref{theolaw} justifies the CMKV-MDP as a macroscopic approximation of  the $N$-agent MDP problem with mean-field interaction.   
Observe that the computation of the conditional gain, expected gain and optimal value of the CMKV-MDP in \eqref{defVMKV}, only requires the variables associated to one agent, called representative agent. Therefore,  we place ourselves in a reduced universe, dismissing other individuals variables, and rena\-ming the representative agent's initial information by $\Gamma$, initial state 
by $\xi$, idiosyncratic noise by $(\eps_t)_{t\in\N}$.  In the sequel, we shall denote by $\Gc$ $=$ $\sigma(\Gamma)$ the $\sigma$-algebra generated by the random variable 
$\Gamma$, hence representing the initial information filtration, and by $L^0(\Gc;\Xc)$  the set of $\Gc$-measurable random variables valued in $\Xc$.  We shall assume that 
the initial state $\xi$ $\in$ $L^0(\Gc;\Xc)$, which means that  the policy has access to the agent's initial state through the initial information filtration $\Gc$. 

An open-loop control for the representative agent is a process $\alpha$, which is adapted to the filtration generated by 
$(\Gamma, (\eps_s)_{s\leq t},(\eps^0_s)_{s\leq t})_{t\in\N}$, and associated to an open-loop policy by: $\alpha_t$ $=$ $\alpha_t^\pi$ $:=$ 
$\pi_t(\Gamma, (\eps_s)_{s\leq t},(\eps^0_s)_{s\leq t})$ for some $\pi$ $\in$ $\Pi_{OL}$. We denote by $\Ac$ the set of open-loop controls, and given $\alpha$ $\in$ $\Ac$, 
$\xi$ $\in$ $L^0(\Gc;\Xc)$,  the state process $X$ $=$ $X^{\xi,\alpha}$ of the representative agent is governed by  
\begin{align} \label{dynXMKV2} 
 X^{}_{t+1} & =   F(X^{}_t, \alpha^{}_t, \P^0_{(X^{}_t,\alpha_t)},\eps_{t+1}, \eps^0_{t+1}),  \quad 	t\in\N, \;\; X_0 = \xi.
\end{align}
For $\alpha$ $=$ $\alpha^\pi$, $\pi$ $\in$ $\Pi_{OL}$, we write indifferently $X^{\xi,\pi}$ $=$ $X^{\xi,\alpha}$, and the expected gain $V^\alpha$ $=$ $V^\pi$ equal to
\begin{align} \label{defValpha}
V^\alpha(\xi) &= \E \Big[ \sum_{t\in\N} \beta^t f(X_t,\alpha_t,\P^0_{(X^{}_t,\alpha_t)}) \Big], 
\end{align}
where we stress the dependence upon the initial state $\xi$. 
 The value function to the CMKV-MDP is then defined by 
 \begin{align}
 V(\xi) &= \sup_{\alpha\in\Ac} V^\alpha(\xi), \quad \xi \in L^0(\Gc;\Xc).  
 \end{align}

 Let us now show how one can lift the CMKV-MDP to a (classical)  MDP on the space of probability measures $\Pc(\Xc)$.  
 We  set $\F^0$ as the filtration generated by the common noise  $\eps^0$. 
Given an open-loop control $\alpha$ $\in$ $\Ac$, and its state process $X$ $=$ $X^{\xi,\alpha}$, denote by $\{\mu_t$ $=$ $\P^0_{X_t}$, $t$ $\in$ $\N\}$, the random $\Pc(\Xc)$-valued process, and notice from  Proposition \ref{condlaw-adapted} that $(\mu_t)_t$  is $\F^0$-adapted.  
From \eqref{dynXMKV2}, and recalling the pushforward measure notation, we have
 \begin{align} \label{dynmu}
 \mu_{t+1} & = F\big(\cdot, \cdot, \P^0_{(X_t,\alpha_t)}, \cdot, \eps^0_{t+1}\big)\star \big(\P^0_{(X_t,\alpha_t)}\otimes\lambda_{\eps_{}}\big), \quad a.s.
 \end{align}
 As the probability distribution $\lambda_\eps$ of the idiosyncratic noise  is a fixed parameter, the above relation means that  
 $\mu_{t+1}$ only depends on $\P^0_{(X_t,\alpha_t)} $ and $\eps^0_{t+1}$.  Moreover, by introducing the so-called {\it relaxed control} associated to the open-loop control $\alpha$ as 
 \begin{align} \label{defrelax}
 \hat\alpha_t(x) & = \; \Lc^0\big(\alpha_t   | X_t = x \big), \quad t \in \N, 
 \end{align} 
 which is  valued in $\hat A(\Xc)$, the set of probability kernels on $\Xc\times A$ (see Lemma \ref{condlaw}), 
 we see from Bayes formula that $\P^0_{(X_t,\alpha_t)}$ $=$ $\mu_t \cdot \hat{\alpha}_t$. 
 The dynamics relation \eqref{dynmu} is then written as 
 \begin{align}
  \mu_{t+1} & = 
  \hat{F}(\mu_t,\hat{\alpha}_t,\eps^0_{t+1}), \quad t \in \N, \;  
 \end{align}
 where the function $\hat F$ $:$ $\Pc(\Xc)\times\hat A(\Xc)\times E^0$ $\rightarrow$ $\Pc(\Xc)$ is defined by
\begin{align} \label{defhatF}
\hat F(\mu,\hat a,e^0) &= F\big(\cdot, \cdot,\mu \cdot \hat{a}, \cdot,e^0\big)\star 
  \big((\mu\cdot \hat{a})\otimes\lambda_{\eps}\big). 
\end{align} 

 On the other hand, by the law of iterated conditional expectation, the expected gain can be written as
\begin{align}
V^\alpha(\xi)  &= 
\E \Big[  \sum_{t\in\N} \beta^t \E^0\big[ f(X_t,\alpha_t,\P^0_{(X_t,\alpha_t)})\big]  \Big],
\end{align}
with the conditional expectation term  equal to
\begin{align}
\E^0\big[f(X_t,\alpha_t,\P^0_{(X_t,\alpha_t)}) \big] &= 
\hat{f}(\mu_t, \hat{\alpha}_t),
\end{align}
where  the function $\hat f$ $:$ $\Pc(\Xc)\times\hat A(\Xc)$ $\rightarrow$ $\R$ is defined by
 \begin{align} \label{defhatf} 
 \hat f(\mu,\hat a) &=   \int_{\Xc\times A}  f(x,a,\mu\cdot \hat{a})(\mu\cdot \hat{a})(\d x,\d a).  
 \end{align}


The above derivation suggests to consider a MDP with state space $\Pc(\Xc)$, action space $\hat A(\Xc)$, a  state transition function 
$\hat F$ as in \eqref{defhatF}, a discount factor $\beta$ $\in$ $[0,1)$, and a reward function $\hat f$ as in \eqref{defhatf}.  A  key point is to endow $\hat A(\Xc)$ with a suitable $\sigma$-algebra in order to have measurable  functions $\hat F$, $\hat f$, and $\F^0$-adapted process $\hat\alpha$ valued in $\hat A(\Xc)$, so that the MDP with characteristics  
$(\Pc(\Xc),\hat A(\Xc),\hat F,\hat f,\beta)$ is well-posed. This issue is investigated in the next sections, first in special cases, and  then in general case by a suitable enlargement of the action space.

\subsection{Case without common noise}
  
When there is no common noise, the original state transition function $F$ is defined from  $\Xc\times A \times \Pc(\Xc\times A)\times E$ into $\Xc$, and the associated function $\hat F$ is then 
defined from $\Pc(\Xc)\times\hat A(\Xc)$ into $\Pc(\Xc)$ by 
\begin{align}
\hat F(\mu,\hat a) &= F\big(\cdot, \cdot, \mu \cdot \hat{a}, \cdot\big)\star \big((\mu\cdot \hat{a})\otimes\lambda_{\eps}\big). 
\end{align}
In this case, we are simply reduced to a deterministic control problem on the state space $\Pc(\Xc)$ with dynamics
\begin{align}
 \mu_{t+1} &= \hat F(\mu_t,\kappa_t), \quad t \in \N, \;\;  \mu_0 \; = \; \mu \in \Pc(\Xc), 
\end{align}
controlled by  $\kappa$ $=$ $(\kappa_t)_{t\in\N}$ $\in$ $\widehat\Ac$, the set of deterministic sequences valued in $\hat A(\Xc)$, and cumulated gain/value function:
 \begin{align}
 \widehat V^{\kappa}(\mu)  &= \sum_{t=0}^\infty \beta^t \hat f(\mu_t,\kappa_t), \quad 
 \widehat V(\mu) \; = \;  \sup_{\kappa\in\widehat\Ac} \widehat V^{\kappa}(\mu), \quad \mu \in \Pc(\Xc),
 \end{align}
 where the bounded function $\hat f$ $:$ $\Pc(\Xc)\times\hat A(\Xc)$ $\rightarrow$ $\R$ is defined as in \eqref{defhatf}. 
 Notice that there are no measurability issues for $\hat F$, $\hat f$, as the problem is deterministic and all the quantities defined above are well-defined.

 \vspace{1mm}

 We aim to prove  the correspondence and equivalence between the MKV-MDP and the above deterministic control problem. From similar derivation as in \eqref{dynmu}-\eqref{defhatf} 
 (by taking directly law under $\P$ instead of $\P^0$),  we clearly see that for any  
 $\alpha$ $\in$ $\Ac$, $V^\alpha(\xi)$ $=$ $\widehat V^{\hat\alpha}(\mu)$, with $\mu$ $=$ $\Lc(\xi)$, and $\hat\alpha$ $=$ $\Rc_\xi(\alpha)$ where $\Rc_\xi$ is the relaxed operator 
 \begin{equation}
 \begin{array}{rcl}
 \Rc_\xi :  \Ac & \longrightarrow & \widehat\Ac \\
 \alpha = (\alpha_t)_t  & \longmapsto & \hat\alpha  = (\hat\alpha_t)_t: \; \hat\alpha_t(x) = \Lc\big( \alpha_t  | X_t^{\xi,\alpha} = x \big), \;\; t \in \N, \; x \in \Xc.  
 \end{array}
 \end{equation}
 It follows that $V(\xi)$ $\leq$ $\widehat V(\mu)$. In order to get the reverse inequality, we have to show that $\Rc_\xi$ is surjective. Notice that this property is not always satisfied:  for instance, when 
 the $\sigma$-algebra generated by $\xi$ is equal to $\Gc$, then for any $\alpha$ $\in$ $\Ac$, $\alpha_0$ is $\sigma(\xi)$-measurable at time $t$ $=$ $0$, 
 and thus  $\Lc(\alpha_0 | \xi)$ is a Dirac distribution, hence cannot  be equal to an arbitrary probability kernel $\kappa_0$ $=$ $\hat a$ $\in$ $\hat A(\Xc)$. 
 We shall then make  the following randomization hypothesis.

\vspace{2mm}

\noindent  {\bf Rand$(\xi,\Gc)$:} There exists a uniform random variable $U$ $\sim$ $\Uc([0,1])$, which is  $\Gc$-measurable and independent of $\xi$ $\in$ $L^0(\Gc;\Xc)$.

\begin{Remark} \label{remrand} 
{\rm The randomization hypothesis  {\bf Rand$(\xi,\Gc)$} implies in particular that $\Gamma$ is atomless, i.e., $\Gc$ is rich enough, and thus $\Pc(\Xc)$ $=$ $\{ \Lc(\zeta): \zeta \in L^0(\Gc;\Xc)\}$.  
Furthermore,  it means that  there is extra randomness in  $\Gc$ besides $\xi$, so that one can freely randomize via the uniform random variable $U$ 
the first action given $\xi$ according to any probability kernel $\hat{a}$. Moreover,  one can extract from $U$, by standard separation  
of the decimals of $U$ (see Lemma 2.21 in \cite{Kallenberg}), 
an i.i.d. sequence of uniform variables $(U_t)_{t\in\N}$, which are $\Gc$-measurable, independent of $\xi$, and  can then be used to randomize the subsequent actions. 
}
\ep
\end{Remark}

\begin{Theorem} {\bf (Correspondence in the no common noise case)}\label{theonocommon}
 
\noindent Assume that {\bf Rand$(\xi,\Gc)$} holds true. Then $\Rc_\xi$ is surjective from $\Ac$ into $\widehat\Ac$, and we have $V(\xi)$ $=$ $\widehat V(\mu)$, for $\mu$ $=$ $\Lc(\xi)$. 
Moreover, for $\epsilon\geq 0$, if $\alpha^\epsilon$ $\in$ $\Ac$ is an $\epsilon$-optimal control for $V(\xi)$, then  $\Rc_\xi(\alpha^\epsilon)$ $\in$ $\widehat\Ac$ is an $\epsilon$-optimal control for 
$\widehat V(\mu)$, and conversely, if $\hat\alpha^\epsilon$ $\in$ $\widehat\Ac$ is an $\epsilon$-optimal control for 
$\widehat V(\mu)$, then any $\alpha^\epsilon$ $\in$ $\Rc_\xi^{-1}(\hat\alpha^\epsilon)$ is an $\epsilon$-optimal control for $V(\xi)$. Consequently, an optimal control for $V(\xi)$ exists iff an optimal control for $\widehat V(\mu)$ exists.
\end{Theorem}
{\bf Proof.}
In view of the above discussion, we only need to prove the surjectivity of $\Rc_\xi$. Fix a control $\kappa\in\widehat{\Ac}$ for the MDP on $\Pc(\Xc)$. By Lemma 2.22 in \cite{Kallenberg}, for all $t\in\N$, there exists  
a measurable function $\mra_t:\Xc\times [0,1]\rightarrow A$ such that $\P_{\mra_t(x,U)}$ $=$  $\kappa_t(x)$, for all $x$ $\in$ $\Xc$. It is then clear that the control $\alpha$ defined recursively 
by $\alpha_t:= \mra_t(X^{\xi,\alpha}_t,U_t)$, where $(U_t)_t$ is an i.i.d. sequence of $\Gc$-measurable uniform variables independent of $\xi$ under  {\bf Rand$(\xi,\Gamma)$}, 
satisfies $\Lc(\alpha_t\mid X^{\xi,\alpha} =x) = \kappa_t(x)$ (observing that $U_t$ is independent of $X_t^{\xi,\alpha}$), 
and thus $\hat{\alpha}=\kappa$, which proves the surjectivity of $\Rc_\xi$.
\ep

\begin{Remark}
{\rm The above correspondence  result  shows in particular that the value function $V$ of the MKV-MDP is law invariant, in the sense that it depends on its initial state $\xi$ only via its probability law $\mu$ $=$ $\Lc(\xi)$,  
for $\xi$ satisfying the randomization hypothesis.
}
\ep
\end{Remark}

\subsection{Case with finite state space $\Xc$ and with common noise}

We consider the case with common noise but when the state space $\Xc$ is finite, i.e.,  its cardinal  $\#\Xc$ is finite, equal to $n$. 
   
 In this case, one can identify $\Pc(\Xc)$ with the simplex $\S^{n-1}$ $=$ $\{ p = (p_i)_{i=1,\ldots,n} \in [0,1]^n: \sum_{i=1}^n p_i = 1\}$, by associating any probability distribution 
 $\mu$ $\in$ $\Pc(\Xc)$ to its weights $(\mu(\{x\}))_{x\in\Xc}$ $\in$ $\S^{n-1}$.  
 We also identify the action space $\hat A(\Xc)$ with $\Pc(A)^n$ by associating any probability kernel $\hat a$ $\in$ $\hat A(\Xc)$ to 
 $(\hat a(x))_{x\in\Xc}$ $\in$ $\Pc(A)^n$, and thus $\hat A(\Xc)$ is naturally endowed with the product $\sigma$-algebra of the Wasserstein metric space $\Pc(A)$.

\begin{Lemma} \label{lemmesur}
Suppose that $\#\Xc$ $=$ $n$ $<$ $\infty$. Then, $\hat F$ in \eqref{defhatF} is a measurable function from $\S^{n-1}\times\Pc(A)^n\times E^0$ into $\S^{n-1}$, $\hat f$ in \eqref{defhatf} is a real-valued measurable function on $\S^{n-1}\times\Pc(A)^n$. Moreover, for any $\xi$ $\in$ $L^0(\Gc;\Xc)$, and  $\alpha$ $\in$ $\Ac$,  the $\Pc(A)^n$-valued process $\hat\alpha$ defined by 
$\hat\alpha_t(x)$ $=$ $\Lc^0( \alpha_t | X_t^{\xi,\alpha} = x)$, $t$ $\in$ $\N$, $x$ $\in$ $\Xc$,  is $\F^0$-adapted. 
\end{Lemma}
{\bf Proof.}
By Lemma  \ref{wassmes}, it is clear, by measurable composition, that we only need to prove that $\Psi$ $:$ $(\mu,\hat{a})\in (\Pc(\Xc),\hat A(\Xc))\mapsto \mu \cdot \hat{a}\in \Pc(\Xc\times A)$ is measurable. However, in this finite state space  case, $\mu \cdot \hat{a}$ is here simply equal to $\sum_{x\in\Xc} \mu(x) \hat{a}(x)$ and, thus $\Psi$ is  clearly measurable. 
\ep

\vspace{2mm}

 In view of Lemma \ref{lemmesur}, the MDP with characteristics  $(\Pc(\Xc)\equiv \S^{n-1},\hat A(\Xc)  \equiv \Pc(A)^n,\hat F,\hat f,\beta)$ is well-posed.  Let us then 
 denote by $\widehat\Ac$ the set of $\F^0$-adapted processes valued in $\Pc(A)^n$, and given $\kappa$ $\in$ $\widehat\Ac$,  consider the controlled dynamics in $\S^{n-1}$ 
 \begin{align} \label{MDPmu}
 \mu_{t+1} &=  \hat F(\mu_t,\kappa_t,\eps_{t+1}^0), \quad t \in \N,  \; \mu_0 = \mu \in \S^{n-1}, 
 \end{align}
 the associated expected gain and value function 
 \begin{align} \label{defvalueMDPmu}
 \widehat V^{\kappa}(\mu) \; = \; \E \Big[ \sum_{t=0}^\infty \beta^t \hat f(\mu_t,\kappa_t) \Big], & \quad \widehat V(\mu) \; = \; \sup_{\kappa\in\hat\Ac}  \widehat V^{\kappa}(\mu).  
 \end{align}
 
  We aim to prove  the correspondence and equivalence between the CMKV-MDP and the MDP \eqref{MDPmu}-\eqref{defvalueMDPmu}. From the derivation  in \eqref{dynmu}-\eqref{defhatf} and by 
  Lemma \ref{lemmesur}, we see that for any  $\alpha$ $\in$ $\Ac$,  $V^\alpha(\xi)$ $=$ $\widehat V^{\hat\alpha}(\mu)$, where $\mu$ $=$ $\Lc(\xi)$, and $\hat\alpha$ $=$ $\Rc^0_\xi(\alpha)$ where 
  $\Rc^0_\xi$ is the relaxed operator 
 \begin{equation} \label{defrelax0}
 \begin{array}{rcl}
 \Rc^0_\xi :  \Ac & \longrightarrow & \widehat\Ac \\
 \alpha = (\alpha_t)_t  & \longmapsto & \hat\alpha  = (\hat\alpha_t)_t: \; \hat\alpha_t(x) = \Lc^0\big( \alpha_t  | X_t^{\xi,\alpha} = x \big), \;\; t \in \N, \; x \in \Xc.  
 \end{array}
 \end{equation}
 It follows that $V(\xi)$ $\leq$ $\widehat V(\mu)$. In order to get the reverse inequality from the surjectivity of  $\Rc^0_\xi$, we need again as in the no common noise case to make some randomization hypothesis. It turns out that when $\Xc$ is finite, this randomization hypothesis is simply reduced  to the atomless property of $\Gamma$.

\begin{Lemma} \label{lemrandom} 
Assume that  $\Gamma$ is  atomless, i.e.,  $\Gc$ is rich enough. Then, any  $\xi$ $\in$ $L^0(\Gc;\Xc)$ taking  a countable number of values,  satisfies {\bf Rand$(\xi,\Gamma)$}.
\end{Lemma}
{\bf Proof.}
Let $S$ be a countable set s.t. $\xi\in S$ a.s., and $\P[\xi=x]>0$ for all $x\in S$. Fix $x\in S$ and denote by $\P_x $ the probability ``knowing $\xi=x$", i.e.,  
$\P_x[B]$ $:=$ $\frac{\P[B, \xi=x]}{\P[\xi=x]}$, for all $B\in\Fc$. It is clear that, endowing $\Omega$ with this probability, $\Gamma$ is still atomless, and so  there exists a $\Gc$-measurable random variable $U_x$ that is uniform under $\P_x$. Then, the random variable $U:=\sum_{x\in S}U_x{\bf 1}_{\xi=x}$ is a $\Gc$-measurable uniform random variable under $\P_x$ for all $x\in S$, which implies that it is a uniform variable under $\P$, independent of $\xi$. 
\ep

\vspace{1mm}

\begin{Theorem} \label{theodiscret} 
{\bf (Correspondance with the MDP on $\Pc(\Xc)$ in the $\Xc$ finite case)}

\noindent 	Assume that  $\Gc$ is rich enough.  Then $\Rc^0_\xi$ is surjective from $\Ac$ into $\widehat\Ac$, and  $V(\xi)$ $=$ $\widehat V(\mu)$, for any $\mu$ $\in$ $\Pc(\Xc)$, 
$\xi$  $\in$ $L^0(\Gc;\Xc)$ s.t.  $\mu$ $=$ $\Lc(\xi)$.   
Moreover, for $\epsilon\geq 0$, if $\alpha^\epsilon$ $\in$ $\Ac$ is an $\epsilon$-optimal control for $V(\xi)$, then  $\Rc^0_\xi(\alpha^\epsilon)$ $\in$ $\widehat\Ac$ is an $\epsilon$-optimal control for 
$\widehat V(\mu)$. Conversely, if $\hat\alpha^\epsilon$ $\in$ $\widehat\Ac$ is an $\epsilon$-optimal control for 
$\widehat V(\mu)$, then any $\alpha^\epsilon$ $\in$ $(\Rc_\xi^0)^{-1}(\hat\alpha^\epsilon)$ is an $\epsilon$-optimal control for $V(\xi)$. 
Consequently, an optimal control for $V(\xi)$ exists iff an optimal control for $\widehat V(\mu)$ exists.
\end{Theorem}
{\bf Proof.}
From the derivation in \eqref{MDPmu}-\eqref{defrelax0},  we only need to prove the surjectivity of $\Rc_\xi^0$. Fix $\kappa\in \widehat{\Ac}$ and let $\bpi_t\in L^0((E^0)^t;\hat{A}(\Xc))$ be such that
$\kappa_t=\bpi_t((\eps^0_s)_{s\leq t})$. As $\Xc$ is finite, by definition of the $\sigma$-algebra on $\hat A(\Xc)$, $\bpi_t$ can be seen as a measurable function in $L^0((E^0)^t\times \Xc;\Pc(A))$. Let $\phi\in L^0(A,\R)$ be an embedding as in Lemma \ref{annexe_embedding}. By Lemma  \ref{wassmes}, we know that $\phi \star \bpi_t$ is in $L^0((E^0)^t\times \Xc; \Pc(\R))$. Given $m\in \Pc(\R)$ we denote by $F^{-1}_m$ the generalized inverse of its distribution function, and it is known that the mapping $m\in (\Pc(\R),\Wc)\mapsto F^{-1}_m \in (L^1_{caglad}(\R),\Vert \cdot\Vert_1)$ is an isometry and is thus measu\-rable. Therefore, $F_{\phi\star \bpi_t}^{-1}$ is in $L^0((E^0)^t\times \Xc;(L^1_{caglad}(\R),\Vert \cdot\Vert_1 ))$. Finally, the mapping $(f,u)\in(L^1_{caglad}(\R),\Vert \cdot\Vert_1 ) \times ([0,1],\Bc([0,1]))\mapsto f(u)\in(\R,\Bc(\R))$ is measurable, since it is the limit of the sequence $n\sum_{i\in \Z} {\bf 1}_{[\frac{i+1}{n}, \frac{i+2}{n} )}(u)\int_{\frac{i}{n}}^{\frac{i+1}{n} } f(y)dy$ when $n\rightarrow \infty$. Therefore, the mapping
\beqs
\mra_t: (E^0)^t\times \Xc\times [0,1] &\longrightarrow & A\\
((e^0_s)_{s\leq t},x,u) &\longmapsto & \phi^{-1} \circ F_{\phi\star \bpi_t((e^0_s)_{s\leq t},x)}^{-1}(u)
\enqs
is measurable. 
We thus define, by induction, $\alpha_t$ $:=$ $\mra_t((\eps^0_s)_{s\leq t}, X^{\xi,\alpha}_t, U_t)$. By construction and by the generalized inverse simulation method, it is clear that  
$\hat{\alpha}_t=\kappa_t$. 
\ep

\vspace{1mm}

\begin{Remark} \label{lipdiscex}
{\rm   
We point out  that when both state space $\Xc$ and action space $A$ are finite,  equipped with the metrics $d(x,x')$ $:=$ ${\bf 1}_{x\not =x'}$, $x,x'$ $\in$ $\Xc$ and $d_A(a,a')$ $:=$ ${\bf 1}_{a\not =a'}$, $a,a'$ $\in$ $A$, the transition function $\hat F$ and reward function $\hat f$ of the lifted MDP on $\Pc(\Xc)$  inherits the Lipschitz condition $({\bf HF_{lip}})$ and $({\bf Hf_{lip}})$ used for the propagation of chaos. 
Indeed, it is known that the Wasserstein distance obtained from $d$ (resp. $d_A$) coincides with twice the total variation distance, and thus  to the $\mathbb{L}^1$ distance when naturally embedding $\Pc(\Xc)$ (resp. $\Pc(A)$) in $[0,1]^{\# \Xc}$ (resp. $[0,1]^{\# A}$). Thus, embedding $\hat A (\Xc)$ in $\Mc_{\#\Xc, \# A}([0,1])$, the set of $\#\Xc \times \# A$ matrices with coefficients valued in $[0,1]$,  we have
\beqs
\Vert \hat{F}(\mu, \hat{a}, e^0), \hat{F}(\nu, \hat{a}', e^0)\Vert_1 \leq (1+K_F)(2\Vert \mu-\mu'\Vert_1 + \underset{x\in \Xc}{\sup}\Vert \hat{a}_{x,\cdot}-\hat{a}'_{x,\cdot}\Vert_1).
\enqs
We obtain a similar property for $f$. In other words, lifting the CMKV-MDP not only turns it into an MDP, but also its state and action spaces $[0,1]^{\# \Xc}$ and $[0,1]^{\# \Xc\times\# A }$ are very standard, and its dynamic and reward are Lipschitz functions with factors of the order of $K_F$ and $K_f$ according to the norm $\Vert \cdot\Vert_1$. Thus, due to the standard nature of this MDP, most MDP algorithms can be applied and their speed will be simply expressed in terms of the original parameters of the CMKV-MDP, $K_F$ and $K_f$.}
\ep	
\end{Remark}

\begin{Remark}
{\rm As in the no common noise case, the correspondence result in the  finite state space case for $\Xc$  shows notably that the value function of the CMKV-MDP is law-invariant. 

The general case (common noise and continuous state space $\Xc$) raises multiple issues for establishing the equivalence  between CMKV-MDP and the lifted MDP on 
$\Pc(\Xc)$.  First, we have to endow the action space $\hat A(\Xc)$ with a suitable $\sigma$-algebra for the lifted MDP to be well-posed: on the one hand,  this $\sigma$-algebra has to be large enough to make the functions  $\hat{F}$ $:$ $\Pc(\Xc)\times \hat{A}(\Xc)\times E^0$ $\rightarrow$ $\Pc(\Xc)$ and $\hat{f}$ $:$ 
$\Pc(\Xc)\times \hat{A}(\Xc)$ $\rightarrow$ $\R$ measurable,  and on the other hand, it should be small enough to make the process $\hat\alpha$ $=$ 
$\Rc_\xi^0(\alpha)$ $\F^0$-adapted for any  control $\alpha$ $\in$ $\Ac$ in the CMKV-MDP.  
Beyond the well-posedness issue of the lifted MDP, the second important concern is the surjectivity of the relaxed operator $\Rc_\xi^0$ from $\Ac$ into $\widehat\Ac$. 
Indeed, if we try to adapt the proof of Theorem \ref{theodiscret} to the case of a continuous state space $\Xc$, the issue is that we cannot in general equip $\hat{A}(\Xc)$ with a $\sigma$-algebra such that 
$L^0((E^0)^t;\hat{A}(\Xc))$ $=$ $L^0((E^0)^t\times \Xc;\Pc(A))$, and thus we cannot see $\bpi_t\in L^0((E^0)^t;\hat{A}(\Xc))$ as an element of $L^0((E^0)^t\times \Xc;\Pc(A))$, which is crucial because the control $\alpha$ (such that $\hat{\alpha}=\kappa$) is defined with $\alpha_t$ explicitly depending upon $\bpi_t((\eps^0_s)_{s\leq t},X_t)$.

In the next section, we shall fix these measurability issues in the general case, and prove the correspondence between the CMKV-MDP and a general lifted MDP on 
$\Pc(\Xc)$.  
}
\ep
\end{Remark}

\section{General case and Bellman fixed point equation in $\Pc(\Xc)$} \label{secmainlift}

We address  the general case with common noise and possibly continuous state space $\Xc$, and our aim is to state the correspondence of the CMKV-MDP with a suitable lifted MDP on $\Pc(\Xc)$  associated to a  Bellman fixed point equation, characterizing the value function, and obtain as a by-product an $\epsilon$-optimal control. 
We proceed as follows: 
\begin{itemize}
\item[(i)]  We first introduce a well-posed lifted MDP on $\Pc(\Xc)$ by enlarging the action space to $\Pc(\Xc\times A)$, and call 
$\tilde V$ the corresponding value function, which satisfies: $V(\xi)$ $\leq$ $\tilde V(\mu)$, for $\mu$ $=$ $\Lc(\xi)$.     
\item[(ii)] We then consider the Bellman equation associated to this well-posed  lifted MDP on $\Pc(\Xc)$, which admits a unique fixed point, called $V^\star$. 
\item[(iii)] Under the randomization hypothesis for $\xi$, we show the existence of an $\epsilon$-randomized feedback policy, which yields both an $\epsilon$-randomized feedback control  for the CMKV-MDP and an $\epsilon$-optimal  feedback control for $\tilde V$. This proves that $V(\xi)$ $=$ $\tilde V(\mu)$ $=$ $V^*(\mu)$, for $\mu$ $=$ 
$\Lc(\xi)$. 
\item[(iv)] Under the condition that $\Gc$ is rich enough, we conclude that $V$ is law-invariant and is equal to $\tilde V$ $=$ $V^\star$, hence satisfies the Bellman equation.   
\end{itemize}
Finally, we show how to compute from the Bellman equation by value or policy iteration  approximate optimal strategy and value function.

\subsection{A general lifted  MDP on $\Pc(\Xc)$}

We start again from the relation \eqref{dynmu} describing the evolution of $\mu_t$ $=$ $\P^0_{X_t}$, $t$ $\in$ $\N$, for a state process $X_t$ $=$ $X_t^{\xi,\alpha}$ controlled by  $\alpha$ $\in$ 
$\Ac$:  
\begin{align} \label{dynmu2} 
\mu_{t+1} &= \;  F(\cdot, \cdot,\P^0_{(X_t,\alpha_t)}, \cdot, \eps^0_{t+1})\star (\P^0_{(X_t,\alpha_t)}\otimes\lambda_{\eps_{}}), \quad a.s.
\end{align}
Now, instead of disintegrating as in Section \ref{seclift},  the conditional law of the pair $(X_t,\alpha_t)$, as  $\P^0_{(X_t,\alpha_t)}$ $=$ $\mu_t\cdot\hat\alpha_t$ 
where $\hat\alpha$ $=$ $\Rc_\xi^0(\alpha)$ is the relaxed control in \eqref{defrelax0}, we directly consider the control process $\balpha_t$ $=$ $\P^0_{(X_t,\alpha_t)}$, $t$ $\in$ $\N$, which is $\F^0$-adapted (see Proposition \ref{condlaw-adapted}), and valued in  the  space of proba\-bility measures $\bA$ $:=$ $\Pc(\Xc\times A)$, naturally endowed with the $\sigma$-algebra of its Wasserstein metric. 
Notice that this $\bA$-valued control $\balpha$ obtained from the CMKV-MDP  has to satisfy by definition  
the marginal constraint  $\text{pr}_{_1}\star\balpha_t$ $=$ $\mu_t$ at any time $t$. In order to tackle this marginal constraint, we shall  rely on the following coupling results.

\begin{Lemma} \label{lemcoupling1} 
{\bf (Measurable coupling)} 
 
\noindent There exists a  measurable function $\zeta$ $\in$ $L^0(\Pc(\Xc)^2\times \Xc\times [0,1];\Xc)$ s.t.  for any $(\mu,\mu')$ $\in$ 
$\Pc(\Xc)$, and if $\xi$ $\sim$ $\mu$, then 
\begin{itemize}
\item $\zeta(\mu,\mu',\xi,U)$ $\sim$ $\mu'$, where  $U$ is an uniform random variable independent of   $\xi$.
\item 
\begin{itemize}
\item[(i)] When $\Xc\subset \R$:  
\beqs
\E \big[ d(\xi,\zeta(\mu,\mu',\xi,U)) \big] &=& \Wc(\mu,\mu'). 
\enqs
\item[(ii)]  In general when $\Xc$ Polish:  $\forall$  $\eps>0$, $\exists \eta>0$ s.t. 
\beqs
\Wc(\mu,\mu') \; < \; \eta  &\Rightarrow&  \E \big[ d(\xi,\zeta(\mu,\mu',\xi,U))\big] \; < \; \eps. 
\enqs
\end{itemize}
\end{itemize}
\end{Lemma}
{\bf Proof.}
See Appendix \ref{coupling}. 
\ep

\begin{Remark}
{\rm
Lemma \ref{lemcoupling1}  can be seen as a measurable version of the well-known coupling result in optimal transport, which states that given $\mu$, $\mu'$ $\in$ $\Pc(\Xc)$, there exists 
$\xi$ and $\xi'$ random variables with $\Lc(\xi)$ $=$ $\mu$, $\Lc(\xi')$ $=$ $\mu'$ such that $\Wc(\mu,\mu')$ $=$ $\E\big[d(\xi,\xi')]$. A similar measurable optimal coupling is proved in 
\cite{fontbana} under the assumption that there exists a transfer function realizing an optimal coupling between $\mu$ and $\mu'$. However, such transfer function does not always exist, for instance when $\mu$ has atoms 
but not $\mu'$. 
Lemma \ref{lemcoupling1} builds a measurable coupling without making such assumption (essentially using the uniform variable $U$ to randomize when $\mu$ has atoms).
}
\ep
\end{Remark}

\vspace{1mm}

From the measurable coupling function $\zeta$ as in Lemma \ref{lemcoupling1}, we define the coupling projection 
$\bp$ $:$ $\Pc(\Xc)\times\bA\rightarrow\bA$ by 
\begin{align}
\bp(\mu,\ba) &= \; \Lc\big( \zeta(\text{pr}_1\star \ba,\mu,\xi',U),\alpha_0 \big), \quad \mu \in \Pc(\Xc), \ba \in \bA, 
\end{align}
where $(\xi',\alpha_0)$ $\sim$ $\ba$, and $U$ is a uniform random variable independent of $\xi'$. 

\begin{Lemma} \label{lemcoupling2} 
{\bf (Measurable coupling projection)} 

\noindent The  coupling projection $\bp$ is a measurable function  from $\Pc(\Xc)\times\bA$ into $\bA$, and for all $(\mu,\ba)$ $\in$ $\Pc(\Xc)\times\bA$:   
\begin{align} \label{trivial}
\text{pr}_{_1}\star \boldsymbol{p}(\mu,\ba) \; = \; \mu,   & \hspace{2mm}    
\mbox{ and if }  \hspace{2mm}  \text{pr}_1\star \ba=\mu, \text{ then } \bp(\mu,\ba)=\ba.
\end{align}
\end{Lemma}
{\bf Proof.} 
By construction, it is clear that $\zeta(\mu,\mu,\xi,U)$ $=$ $\xi$, and so relation \eqref{trivial} is obvious. 
The only result that is not trivial  is the measurability of $\bp$. Observe  that $\bp (\mu,\ba)=g(\mu,\ba,\cdot,\cdot,\cdot)\star (\ba\otimes \Uc([0,1]))$ where $g$ is the measurable function  
\beqs
g : \Pc(\Xc)\times \Pc(\Xc\times A)\times \Xc\times A\times [0,1]&\longrightarrow &\Xc \times A\\
(\mu,\ba, x, a, u) & \longmapsto & (\zeta (\text{pr}_1 \star \ba, \mu, x, u),a)
\enqs
We thus conclude by Lemma \ref{wassmes}. 
\ep


\vspace{3mm}

By using this coupling  projection $\bp$, we see that the dynamics \eqref{dynmu2} can be written as 
\begin{align} \label{dynmu3} 
\mu_{t+1} &= \; \bF(\mu_t,\balpha_t,\eps_{t+1}^0),  \quad t \in \N, 
\end{align} 
where the function $\bF$ $:$ $\Pc(\Xc)\times\bA\times E^0$ $\rightarrow$ $\Pc(\Xc)$ defined by 
\begin{align}
\bF(\mu,\ba,e^0) &= \; F(\cdot, \cdot,\bp(\mu,\ba),\cdot,e^0)\star\big(\bp(\mu,\ba)\otimes\lambda_{\eps_{}}\big), 
\end{align}
is clearly measurable. Let us also define the measurable  function $\tilde f$ $:$ $\Pc(\Xc)\times\bA$  $\rightarrow$ $\R$ by 
\begin{align}
\tilde f(\mu,\ba) & = \;  \int_{\Xc\times A}  f(x,a,\bp(\mu,\ba))\bp(\mu,\ba)(\d x,\d a).   
\end{align}

The MDP with characteristics $(\Pc(\Xc),\bA=\Pc(\Xc\times A),\tilde F,\tilde f,\beta)$ is then well-posed. Let us then denote by $\bAc$ the set of 
$\F^0$-adapted processes  valued in $\bA$, and given an open-loop control $\bnu$ $\in$ $\bAc$, consider the controlled dynamics
\begin{align} \label{dynmugen}
\mu_{t+1} &= \; \bF(\mu_t,\bnu_t,\eps_{t+1}^0),  \quad t \in \N, \; \mu_0 = \mu \in \Pc(\Xc),   
\end{align} 
with associated expected gain/value function
\begin{align} \label{deftildeV}
\widetilde V^{\bnu}(\mu) & = \; 
\E \Big[ \sum_{t\in\N} \beta^t \tilde f(\mu_t,\bnu_t) \Big],   \quad\quad \widetilde V(\mu) \; = \; \sup_{\bnu\in\bAc} \widetilde V^{\bnu}(\mu). 
\end{align}

Given $\xi$ $\in$ $L^0(\Gc;\Xc)$, and $\alpha$ $\in$ $\Ac$, we set $\balpha$ $=$ $\bLi^0_\xi(\alpha)$, where $\bLi_\xi^0$ is the lifted operator 
\begin{equation} \label{defR0}
 \begin{array}{rcl}
 \bLi^0_\xi :  \Ac & \longrightarrow & \bAc \\
 \alpha = (\alpha_t)_t  & \longmapsto & \balpha  = (\balpha_t)_t: \; \balpha_t = \P^0_{(X_t^{\xi,\alpha},\alpha_t)}, \;\; t \in \N.  
 \end{array}
 \end{equation}
By construction from \eqref{dynmu3}, we see that $\mu_t$ $=$ $\P^0_{X_t^{\xi,\alpha}}$, $t$ $\in$ $\N$, follows the dynamics \eqref{dynmugen} with the control $\bnu$ $=$ $\bLi^0_\xi(\alpha)$ 
$\in$ $\bAc$. Moreover, by the law of iterated conditional expectation, and the definition of $\tilde f$, 
the expected gain of the CMKV-MDP can be written as
\begin{align}
V^\alpha(\xi)  &= \; 
\E \Big[  \sum_{t\in\N} \beta^t \E^0\big[ f(X_t^{\xi,\alpha},\alpha_t,\P^0_{(X_t^{\xi,\alpha},\alpha_t)})\big]  \Big] \\
& =  \; \E \Big[  \sum_{t\in\N} \beta^t \tilde f(\P^0_{X_t^{\xi,\alpha}},\balpha_t) \Big] \; = \; \widetilde V^{\balpha}(\mu), \;\; \mbox{ with } \; \mu = \Lc(\xi).  \label{VtildeV} 
\end{align}
It follows that $V(\xi)$ $\leq$ $\widetilde V(\mu)$, for $\mu = \Lc(\xi)$.  Our goal is  to prove the equality, which implies in particular that $V$ is law-invariant, and to obtain as a by-product the corresponding Bellman fixed point equation that characterizes analytically the solution to the CMKV-MDP.

\subsection{Bellman fixed point on $\Pc(\Xc)$}

We derive and study the Bellman equation corresponding to the general lifted MDP \eqref{dynmugen}-\eqref{deftildeV} on $\Pc(\Xc)$. 

By defining this MDP on the canonical space $(E^0)^{\N}$, we identify $\eps^0$ with the canonical identity function in $(E^0)^{\N}$, and $\eps^0_t$ with the $t$-th projection in $(E^0)^{\N}$. We also denote by $\theta:(E^0)^{\N}\rightarrow (E^0)^{\N}$ the shifting operator, 
defined by $\theta ((e^0_t)_{t\in\N})$ $=$ $(e^0_{t+1})_{t\in\N}$. 
Via this identification, an open-loop  control $\bnu$ $\in$ $\bAc$ is a sequence $(\bnu_t)_t$  where $\bnu_t$ is a measurable function from  $(E^0)^t$ into $\bA$,  with the convention that $\bnu_0$ is simply a constant in 
$\bA$.
Given  $\bnu$ $\in$ $\bAc$, and $e^0$ $\in$ $E^0$, we define 
$\vec{\bnu}^{e^0}$ $:=$ $(\vec{\bnu}_ t^{e^0})_t$ $\in$ $\bAc$, where $\vec{\bnu}_t^{e^0}(.)$ $:=$ $\bnu_{t+1}(e^0,.)$, $t\in\N$. Given $\mu$ $\in$ $\Pc(\Xc)$, and 
$\bnu$ $\in$ $\bAc$, we denote by $(\mu_t^{\mu,\bnu})_t$ the solution to \eqref{dynmugen} on the canonical space, which satisfies the flow property
\begin{align}
\big(\mu_{t+1}^{\mu,\bnu},\bnu_{t+1}\big) & \equiv \; \big(\mu_t^{\mu_1^{\mu,\bnu},\vec{\bnu}^{\eps^0_1}(\theta(\eps^0))}, \vec{\bnu}_t^{\eps^0_1}(\theta(\eps^0))\big), \quad t \in \N.   
\end{align}
where $\equiv$ denotes the equality between functions on the canonical space. Given that 
$\eps^0_1\independent \theta (\eps^0)\overset{d}{=}\eps^0$, we obtain that the expected gain of this MDP in \eqref{deftildeV} satisfies the relation
\begin{align} \label{tildeVnu} 
\widetilde V^{\bnu}(\mu) &= \; \tilde f(\mu,\bnu_0)  + \beta \E \Big[ \widetilde V^{\vec{\bnu}^{\eps_1^0}}(\mu_1^{\mu,\bnu})   \Big].  
\end{align}
Let us denote by $L^\infty(\Pc(\Xc))$ the set of bounded real-valued functions on $\Pc(\Xc)$, 
and by $L^\infty_m(\Pc(\Xc))$ the subset of measurable functions in $L^\infty(\Pc(\Xc))$. 
We then introduce the Bellman ``operator'' $\Tc:L^\infty_m(\Pc(\Xc))\rightarrow L^\infty(\Pc(\Xc))$ defined 
for any $W$ $\in$ $L^\infty_m(\Pc(\Xc))$ by: 
\begin{align} \label{defTc} 
[\Tc W](\mu) & := \; \sup_{\ba \in \bA} \Big\{  \tilde f(\mu,\ba) + \beta \E \big[ W\big( \tilde F(\mu,\ba,\eps_1^0)\big) \big] \Big\},  \quad \mu \in \Pc(\Xc). \end{align}
Notice that the $\sup$ can a priori lead to a non measurable function $\Tc W$.
This Bellman operator is consistent with the lifted MDP derived in Section \ref{seclift}, with cha\-racteristics $(\Pc(\Xc),\hat A(\Xc),\hat F,\hat f,\beta)$, although this MDP is not always well-posed. Indeed, its corresponding Bellman operator is well-defined as it only involves the random variable  $\eps_1^0$ at time $1$, hence only requires the measurability of $e^0$ $\mapsto$ 
$\hat F(\mu,\hat a,e^0)$, for any $(\mu,\hat a)$ $\in$ $\P(\Xc)\times\hat A(\Xc)$ (which holds true), and it turns out that it coincides with $\Tc$.

\begin{Proposition} \label{equivTc} 
For any $W$ $\in$ $L^\infty_m(\Pc(\Xc))$, and $\mu$ $\in$ $\Pc(\Xc)$, we have
\begin{align} \label{Belhat} 
[\Tc W](\mu) \; = \;  \sup_{\hat a \in \hat A(\Xc)}  [\hat\Tc^{\hat a} W](\mu)  \; = \;  \sup_{\mra \in  L^0(\Xc\times [0,1]; A)} [\T^\mra W](\mu),
\end{align}
where $\hat\Tc^{\hat a}$ and $\T^\mra$ are the operators defined on $L^\infty(\Pc(\Xc))$ by 
\begin{align}
[\hat\Tc^{\hat a} W](\mu) & = \;   \hat f(\mu,\hat a) + \beta \E \big[ W\big( \hat F(\mu,\hat a,\eps_1^0)\big) \big],     \\
[\T^{\mra} W](\mu) & = \; \E \Big[ f(\xi,\mra(\xi,U),\Lc(\xi,\mra(\xi,U))) + \beta W\big(  \P^0_{F(\xi,\mra(\xi,U),\Lc(\xi,\mra(\xi,U)),\eps_1,\eps_1^0)}   \big) \Big],  \label{defTa} 
\end{align} 
for any $(\xi,U)$ $\sim$  $\mu\otimes\Uc([0,1])$ (it is clear that the right-hand side in \eqref{defTa} does not depend on the choice of such $(\xi,U)$).  Moreover, we have 
\begin{align}
[\Tc W](\mu) & = \; \sup_{\alpha_0 \in L^0(\Omega; A)} \E \Big[ f(\xi,\alpha_0,\Lc(\xi,\alpha_0)) + \beta W\big(  \P^0_{F(\xi,\alpha_0,\Lc(\xi,\alpha_0),\eps_1,\eps_1^0)}   \big) \Big].  \label{BelF} 
\end{align}
\end{Proposition}
{\bf Proof.} Fix $W$ $\in$ $L^\infty_m(\Pc(\Xc))$, and $\mu$ $\in$ $\Pc(\Xc)$. 
Let $\ba$ be arbitrary in $\bA$.  Since $\bp(\mu,\ba)$ has first marginal equal to $\mu$, there exists by assertion  3 in Lemma \ref{condlaw}  a probability kernel 
$\hat a$ $\in$ $\hat A(\Xc)$ such that $\bp(\mu,\ba)$ $=$ $\mu\cdot\hat a$. Therefore, $\tilde F(\mu,\ba,e^0)$ $=$ $\hat F(\mu,\hat a,e^0)$, $\tilde f(\mu,\ba)$ $=$ $\hat f(\mu,\hat a)$, which implies that $[\Tc W](\mu)$ $\leq$ 
$ \sup_{\hat a \in \hat A(\Xc)} [\hat\Tc^{\hat a} W](\mu)$ $=:$ $\T^1$.

Let us consider  the operator $\Rc$ defined by 
 \begin{equation*}
 \begin{array}{rcl}
 \Rc :  L^0(\Xc\times [0,1];A)  & \longrightarrow & \hat A(\Xc)  \\
 \mra   & \longmapsto &   \hat a: \;\; \hat a(x) = \Lc\big( \mra(x,U) \big), \;\; x \in \Xc,  \; U \sim \Uc([0,1]), 
 \end{array}
 \end{equation*}
 and notice that it is surjective from $L^0(\Xc\times [0,1];A)$ into $\hat A(\Xc)$, by Lemma 2.22 in \cite{Kallenberg}.  
 By noting that for any  $\mra$ $\in$ $L^0(\Xc\times [0,1];A)$, and  $(\xi,U)$ $\sim$  $\mu\otimes\Uc([0,1])$, we have $\Lc\big( \xi,\mra(\xi,U) \big)$ $=$ $\mu\cdot \Rc(\mra)$, it follows that 
$[ \mathbb{T}^\mra W](\mu)$ $=$  $[ \hat\Tc^{\Rc(\mra)} W](\mu)$. Since $\Rc$ is surjective, this yields $\T^1$ $=$ $\sup_{\mra\in L^0(\Xc\times [0,1];A)}[\mathbb{T}^\mra W](\mu)$ $=:$ $\T^2$. 
 
Denote by $\T^3$ the right-hand-side in \eqref{BelF}.  It is clear that $\mathbb{T}^2$ $\leq$ $\mathbb{T}^3$. Conversely, let  $\alpha_0$ $\in$ $L^0(\Omega;A)$. 
We then set $\ba$ $=$ $\Lc(\xi,\alpha_0)$ $\in$ $\Pc(\Xc\times A)$, and notice that the first marginal of $\ba$ is $\mu$. 
Thus, $\bp(\mu,\ba)$ $=$ $\Lc(\xi,\alpha_0)$, and so
\begin{align}
\tilde f(\mu,\ba) & = \;  \int_{\Xc\times A}  f(x,a,\bp(\mu,\ba))\bp(\mu,\ba)(\d x,\d a) \; = \; \E \big[ f(\xi,\alpha_0,\Lc(\xi,\alpha_0)) \big] \\
\tilde F(\mu,\ba,\eps_1^0) & = \;  F(\cdot, \cdot,\bp(\mu,\ba),\cdot,\eps_1^0)\star\big(\bp(\mu,\ba)\otimes\lambda_{\eps_{}}\big) \; = \; \P^0_{F(\xi,\alpha_0,\Lc(\xi,\alpha_0),\eps_1,\eps_1^0)}.   
\end{align}
We deduce that $\mathbb{T}^3$ $\leq$ $[\Tc W](\mu)$, which gives finally the equalities \eqref{Belhat} and \eqref{BelF}.  
\ep

\vspace{3mm}

We state the basic properties of the Bellman operator $\Tc$.

\begin{Proposition}\label{lemTc}\label{proregul} 
Assume that  $({\bf H_{lip}})$ holds true.
(i)  The operator $\Tc$ is monotone increasing: for $W_1,W_2\in L_m^\infty(\Pc(\Xc))$, if $W_1$ $\leq$ $W_2$,  then $\Tc W_1$ $\leq$ $\Tc W_2$.
(ii) Furthermore, it is contracting on $L^\infty_m(\Pc(\Xc))$ with Lipschitz factor $\beta$, and admits a unique fixed point 
in $L^\infty_m(\Pc(\Xc))$, denoted by $V^\star$, hence solution to: 
\begin{align}
V^\star & = \; \Tc V^\star. 
\end{align}
(iii) $V^\star$ is $\gamma$-H\"older, with $\gamma= \min \left(1, \frac{\vert \ln \beta\vert}{\ln (2K_F)}\right)$, i.e. there exists  some positive constant $K_\star$ (depending only on  $K_F$, $K_f$, $\beta$, and explicit in the proof),  
such that
\begin{align}
\big\vert V^\star(\mu) - V^\star(\mu') \big\vert & \leq \;  K_\star\Wc(\mu,\mu')^\gamma, \quad \forall  \mu,\mu' \in \Pc(\Xc).   
\end{align}
\end{Proposition}
\noindent {\bf Proof.}  (i) The monotonicity of $\Tc$ is shown by standard arguments. 

\noindent(ii) The $\beta$-contraction property of $\Tc$ is also obtained by standard arguments. Let us now prove by induction that the iterative sequence $V_{n+1}$ $=$ $\Tc V_n$, with $V_0\equiv 0$ is well defined and such that 
\begin{align}\label{prop-n}
\vert V_{n}(\mu)-V_{n}(\mu')\vert\leq 2K_f \sum_{t=0}^\infty \beta^t \min((2K_F)^t \Wc (\mu,\mu'),\Delta_\Xc)
\end{align}
for all $n\in \N$. The property is obviously satisfied for $n=0$. Assume that the property holds true for a fixed $n\in \N$, and let us prove it for $n+1$. First of all, the inequality \eqref{prop-n} implies that $V_n$ is continuous, and thus $V_n\in L^\infty_m(\Pc(\Xc)$. Therefore, $V_{n+1}=\Tc V_n$ is well defined. Fix $\mu,\mu'\in \Pc(\Xc)$. In order to use the expression \eqref{BelF} of the Bellman operator $\Tc$, we consider an optimal coupling $(\xi,\xi')$ of $\mu$ and $\mu'$, i.e. $\xi\sim \mu$, $\xi'\sim \mu'$, and $\E [d(\xi,\xi')]=\Wc(\mu,\mu')$, and fix an $A$-valued random variable $\alpha_0$. Let us start with two preliminary estimations: under (${\bf H_{lip}}$), we have
\begin{align}
\E \left[\vert f(\xi, \alpha_0, \Lc(\xi,\alpha_0))-f(\xi', \alpha_0, \Lc(\xi',\alpha_0))\vert\right]&\leq K_f(\E [d(\xi,\xi')]+\Wc(\Lc(\xi,\alpha_0),\Lc(\xi',\alpha_0)))\nonumber \\
&\leq K_f(\E [d(\xi,\xi')]+\E [d((\xi,\alpha_0),(\xi',\alpha_0))]) \nonumber \\
&\leq 2 K_f\E [d(\xi,\xi')]=2 K_f\Wc(\mu,\mu'). \label{f-estim}
\end{align}
Similarly, for $e^0\in E^0$, we have
\begin{align}\label{F-estim}
\E [d(F(\xi,\alpha_0,\Lc(\xi,\alpha_0),\eps^1_1,e^0),F(\xi',\alpha_0,\Lc(\xi,\alpha_0),\eps^1_1,e^0))]&\leq 2 K_F\Wc(\mu,\mu').
\end{align}
Now, we prove the hereditary property. The definition of $\Tc$ and $V_{n+1}$ combined with \eqref{f-estim} and the induction hypothesis, imply that
\beqs
\vert V_{n+1}(\mu)-V_{n+1}(\mu')\vert 
&\leq & 2K_f\Wc(\mu,\mu')+\beta \E [ 2K_f\sum \beta^t \min((2K_F^t \Wc (\mu_1,\mu'_1),\Delta_\Xc)]
\enqs
where $\mu_1=\Lc^0(F(\xi,\alpha_0,\Lc(\xi,\alpha_0),\eps^1_1,\eps^0_1))$ and $\mu'_1=\Lc^0(F(\xi',\alpha_0,\Lc(\xi',\alpha_0),\eps^1_1,\eps^0_1))$. By Jensen's inequality and \eqref{F-estim}, we have
\beqs
\vert V_{n+1}(\mu)-V_{n+1}(\mu')\vert &\leq & 2K_f\min(\Wc(\mu,\mu'),\Delta_\Xc)+\beta   2K_f\sum \beta^t \min((2K_F^t \E\Wc (\mu_1,\mu'_1),\Delta_\Xc)\\
&\leq & 2K_f\min(\Wc(\mu,\mu'),\Delta_\Xc)+\beta   2K_f\sum \beta^t \min((2K_F^t 2K_F \Wc(\mu,\mu'),\Delta_\Xc)\\
&\leq & 2K_f\sum \beta^t \min((2K_F^t \Wc(\mu,\mu'),\Delta_\Xc). 
\enqs
This concludes the induction and proves that $V_n$ is well defined and satisfies the inequality \eqref{prop-n} for all $n\in \N$. As $\Tc$ is $\beta$-contracting, a standard argument from the proof of the Banach fixed point theorem shows that $(V_n)_n$ is a Cauchy sequence in the complete metric space $L^\infty_m(\Pc(\Xc))$, and therefore admits a limit $V^\star\in L^\infty_m(\Pc(\Xc))$. Notice that
\beqs 
V^\star(\mu)=\lim_n V_{n+1}(\mu)=\lim_n \Tc V_{n}(\mu)=\Tc V^\star
\enqs
by continuity of the contracting operator $\Tc$.

\noindent (iii) By sending  $n$ to infinity in \eqref{prop-n},  we obtain
\beqs
\vert V(\mu)-V(\mu')\vert\leq 2K_f \sum_{t=0}^\infty \beta^t \min((2K_F)^t \Wc (\mu,\mu'),\Delta_\Xc).
\enqs
Then, a similar derivation as in the end of proof of  Theorem \ref{theoL1chaos} shows the required $\gamma$-H\"older property of $V^\star$. 
\ep

\begin{Remark} 
{\rm  In the proof of Proposition \ref{lemTc}, one could also have proved that the set $\Sc$ of functions $W:\Pc(\Xc)\rightarrow \R$ such that 
\begin{align}
\vert W(\mu)-W(\mu')\vert\leq 2K_f \sum_{t=0}^\infty \beta^t \min((2K_F)^t \Wc (\mu,\mu'),\Delta_\Xc)
\end{align}
for all $\mu,\mu'\in\Pc(\Xc)$ is a complete metric space, as it is a closed set of the complete metric space $L^\infty_m(\Pc(\Xc))$, and is stabilized by the contracting operator $\Tc$ (which is essentially proved by replacing $V_n$ by $W$ in the proof). One could then have invoked the Banach fixed point theorem on this set $\Sc$, implying the existence and uniqueness of the fixed point $V^\star$. Notice that this argument would not work if we considered, instead of $\Sc$, the set of $\gamma$-Hölder continuous functions. Indeed, while it is true that such set is stabilized by $\Tc$ (it essentially follows from \eqref{f-estim} and \eqref{F-estim}), the set of $\gamma$-Hölder continuous functions is not closed in $L^\infty_m(\Pc(\Xc))$ (and thus not a complete metric space): there might indeed exist a converging sequence of $\gamma$-Hölder continuous functions with multiplicative factors (in the Hölder property) tending toward infinity, such that the limit function is not $\gamma$-Hölder anymore.
}
\ep
\end{Remark}

\vspace{3mm}

As a consequence of  Proposition \ref{lemTc}, we can easily show the following relation between the value function $\tilde V$ of the general lifted MDP, and the fixed point $V^\star$ of the Bellman operator.

\begin{Lemma} \label{lemtildeVstar}
 For all $\mu$ $\in$ $\Pc(\Xc)$,
 we have $\tilde{V}(\mu)\leq V^\star(\mu)$.
 \end{Lemma}
 \noindent{\bf Proof.}
From \eqref{tildeVnu}, we have
\beqs 
& & \inf_{\mu \in \Pc(\Xc)} \big\{ V^\star(\mu)-\widetilde V^{\bnu}(\mu) \big\}  \\
&\geq& \inf_{\mu\in\Pc(\Xc)}  \left\{ \Tc V^\star(\mu)-\left(\tilde f(\mu,\bnu_0)  + \beta \E \Big[ V^{\star}(\mu_1^{\mu,\bnu})   \Big] \right) +\beta \E \Big[ V^{\star}(\mu_1^{\mu,\bnu})-\widetilde V^{\vec{\bnu}^{\eps_1^0}}(\mu_1^{\mu,\bnu})\Big]\right\} \\
&\geq& \beta \E \Big[ V^{\star}(\mu_1^{\mu,\bnu})-\widetilde V^{\vec{\bnu}^{\eps_1^0}}(\mu_1^{\mu,\bnu})\Big]\geq \beta \inf_{\mu\in\Pc(\Xc)} \big\{ V^\star(\mu)-\widetilde V^{\bnu}(\mu) \big\}. 
\enqs
This shows that $\inf_{\mu\in\Pc(\Xc)}(V^\star(\mu)-\widetilde V^{\bnu}(\mu))\geq 0$, hence
\beqs 
\widetilde V^{\bnu}(\mu)\leq V^\star(\mu)\quad \forall \mu\in\Pc(\Xc).
\enqs 
Taking the sup over $\bnu\in\bAc$, we obtain the required result. 
\ep


\subsection{Building $\epsilon$-optimal randomized feedback controls}

We  aim to prove rigorously  the equality $\tilde V$ $=$ $V^\star$, i.e., the value function $\tilde V$ of the general lifted MDP satisfies the Bellman fixed point equation: $\tilde V$ $=$ $\Tc\tilde V$, 
and also to show the existence of an $\epsilon$-optimal control for $\tilde V$.  
Notice that it cannot be obtained directly from classical theory of MDP as we consider here open-loop controls $\bnu$ $\in$ $\bAc$ while MDP usually deals with feedback controls on 
finite-dimensional spaces.  Anyway, following the standard notation in MDP theory with state space $\Pc(\Xc)$ and action space $\bA$, and in connection with the Bellman operator in \eqref{defTc}, we introduce, for $\bpi$ $\in$ $L^0(\Pc(\Xc);\bA)$ (the set of measurable functions from $\Pc(\Xc)$ into $\bA$) called (measurable) feedback policy, the so-called $\bpi$-Bellman 
operator $\Tc^{\bpi}$  on $L^\infty(\Pc(\Xc))$, defined for $W \in L^\infty(\Pc(\Xc))$  by
\begin{align} \label{defTpi} 
[\Tc^{\bpi} W](\mu) &=   \tilde f(\mu,\bpi(\mu)) + \beta \E \big[ W\big( \tilde F(\mu,\bpi(\mu),\eps_1^0)\big) \big], \quad \mu \in \Pc(\Xc). 
\end{align}
 
As for the Bellman operator $\Tc$, we have the basic properties on the operator $\Tc^{\bpi}$. 

\begin{Lemma} \label{lemTcpi} 
Fix $\bpi$ $\in$ $L^0(\Pc(\Xc);\bA)$. 
\begin{itemize}
\item[(i)] The operator $\Tc^{\bpi}$ is contracting on $L^\infty(\Pc(\Xc))$ with Lipschitz factor $\beta$, and admits a unique fixed point denoted $\tilde V^{\bpi}$.  
\item[(ii)] Furthermore, it is monotone increasing: for $W_1,W_2\in L^\infty(\Pc(\Xc))$, if $W_1$ $\leq$ $W_2$,  then $\Tc^{\bpi} W_1$ $\leq$ $\Tc^{\bpi} W_2$.
\end{itemize}
\end{Lemma}

\begin{Remark} \label{remVpi} 
{\rm  It is well-known from MDP theory that  the fixed point $\tilde V^{\bpi}$ to the operator $\Tc^{\bpi}$ is equal to
\begin{align}
\tilde V^{\bpi}(\mu)  &= \; \E \Big[ \sum_{t\in \N} \tilde f(\mu_t,\bpi(\mu_t)) \Big], 
\end{align}
where $(\mu_t)$ is the MDP in \eqref{dynmugen} with the feedback and stationary control $\bnu^{\bpi}$ $=$ $(\bnu^{\bpi}_t)_t$  $\in$ $\bAc$ defined by $\bnu^{\bpi}_t$ $=$ $\bpi(\mu_t)$, $t$ $\in$ $\N$. In the sequel, we shall then identify by misuse of notation $\tilde V^{\bpi}$ and  $\tilde V^{\bnu^{\bpi}}$ as defined in \eqref{deftildeV}.   
}
\ep
\end{Remark}

\vspace{2mm}

Our ultimate goal being to solve the  CMKV-MDP, we introduce a subclass of feedback policies for the lifted MDP.  

\begin{Definition}
{\bf (Lifted randomized feedback policy)} 

\noindent A feedback policy $\bpi\in L^0(\Pc(\Xc);\bA)$ is a lifted randomized feedback policy if there exists a measurable function $\mfa$  $\in$ 
$L^0(\Pc(\Xc)\times \Xc\times [0,1];A)$, called randomized feedback policy,  such that $\big(\xi,\mfa(\mu,\xi,U)\big)$ $\sim$ $\bpi(\mu)$, for all $\mu\in\Pc(\Xc)$, with $(\xi,U)\sim\mu\otimes \Uc([0,1])$. 
\end{Definition}

\begin{Remark} \label{remlift} 
{\rm Given $\mfa$ $\in$ $L^0(\Pc(\Xc)\times \Xc\times [0,1];A)$, denote by $\bpi^{\mfa}$ $\in$ $L^0(\Pc(\Xc);\bA)$ the associated lifted  randomized feedback policy, i.e., 
$\bpi^\mfa(\mu)$ $=$ $\Lc\big(\xi,\mfa(\mu,\xi,U)\big)$,  for $\mu\in\Pc(\Xc)$, and $(\xi,U)\sim\mu\otimes \Uc([0,1])$.  By definition of the $\pi$-Bellman operator 
$\Tc^{\bpi^{\mfa}}$ in \eqref{defTpi}, and observing that $\bp(\mu,\bpi^{\mfa}(\mu))$ $=$ $\bpi^{\mfa}(\mu)$ $=$ $\Lc\big(\xi,\mra^\mu(\xi,U)\big)$, where we set $\mra^\mu$ $=$ $\mfa(\mu,\cdot,\cdot)$ $\in$ 
$L^0(\Xc\times [0,1]:A)$, we see (recalling the notation in \eqref{defTa}) that for all $W$ $\in$ $L^\infty(\Pc(\Xc))$, 
\begin{align} \label{lienT} 
[\Tc^{\bpi^{\mfa}} W](\mu) & = [\T^{\mra^\mu} W](\mu), \quad \mu \in \Pc(\Xc). 
\end{align}
On the other hand, let $\xi$ $\in$ $L^0(\Gc;\Xc)$ be some initial state satisfying the randomization hypothesis {\bf Rand$(\xi,\Gc)$}, and  
denote by $\alpha^\mfa$ $\in$ $\Ac$ the randomized feedback stationary 
control defined by $\alpha_t^\mfa$ $=$ $\mfa(\P_{X_t}^0,X_t,U_t)$, where $X$ $=$ $X^{\xi,\alpha^\mfa}$ is the state process in \eqref{dynXMKV2} of the CMKV-MDP, 
and $(U_t)_t$ is an i.i.d. sequence of uniform $\Gc$-measurable random variables independent of $\xi$.  
By construction,  the associated lifted control  $\boldsymbol{\alpha^\mfa}$ $=$ $\bLi^0_\xi(\alpha^\mfa)$  
satisfies  $\boldsymbol{\alpha_t^\mfa}$ $=$ $\P^0_{(X_t,\alpha_t^\mfa)}$ $=$ $\bpi^\mfa(\mu_t)$, where $\mu_t$ $=$ $\P_{X_t}^0$, $t$ $\in$ $\N$. 
Denoting by $V^\mfa$ $:=$ $V^{\alpha^\mfa}$ the associated expected gain of the CMKV-MDP, and recalling Remark \ref{remVpi}, we see from \eqref{VtildeV} that $V^\mfa(\xi)$ $=$ 
$\tilde V^{\bnu^{\bpi^\mfa}}(\mu)$ $=$ $\tilde V^{\bpi^\mfa}(\mu)$,  where  $\mu$ $=$ $\Lc(\xi)$.   
}
\ep
\end{Remark}

\vspace{2mm}

We show a verification type result for the general lifted MDP, and as a byproduct for the CMKV-MDP,  by means of the Bellman operator.

\begin{Proposition}\label{lemverif} 
{\bf (Verification result)} 

\noindent Fix $\epsilon\geq0$, and suppose that  there exists an $\epsilon$-optimal feedback policy $\bpi_\epsilon$ $\in$ $L^0(\Pc(\Xc);\bA)$ 
for $V^\star$ in the sense that
\begin{align}
V^\star &\leq \;  \Tc^{\bpi_\epsilon} V^\star  +\epsilon.   
\end{align}
Then, $\bnu^{\pi_\epsilon}$ $\in$ $\bAc$  is $\frac{\epsilon}{1-\beta}$-optimal for $\tilde V$, i.e.,  $\tilde V^{\bpi_\epsilon}$ $\geq$ $\tilde V - \frac{\epsilon}{1-\beta}$, 
and we have $\tilde V$ $\geq$ $V^\star - \frac{\epsilon}{1-\beta}$.  Furthermore, if $\bpi_\epsilon$ is a lifted randomized feedback policy, i.e., $\bpi_\eps$ $=$ 
$\bpi^{\mfa_\epsilon}$, for some $\mfa_\eps$ $\in$ $L^0(\Pc(\Xc)\times \Xc\times [0,1];A)$, then under {\bf Rand$(\xi,\Gc)$}, $\alpha^{\mfa_\epsilon}$ $\in$ $\Ac$ is an $\frac{\epsilon}{1-\beta}$-optimal control for $V(\xi)$, i.e., 
$V^{\mfa_\epsilon}(\xi)$ $\geq$ $V(\xi) - \frac{\epsilon}{1-\beta}$, and we have  $V(\xi)$ $\geq$ $V^\star(\mu) - \frac{\epsilon}{1-\beta}$, for $\mu$ $=$ $\Lc(\xi)$. 
\end{Proposition}
{\bf Proof.} Since $\tilde V^{\bpi_\epsilon}$ $=$ $\Tc^{\bpi_\epsilon} \tilde V^{\bpi_\epsilon}$, and recalling from Lemma \ref{lemtildeVstar} that $V^\star$ $\geq$ $\tilde V$ $\geq$ 
$\tilde V^{\bpi_\epsilon}$, we have for all $\mu$ $\in$ $\Pc(\Xc)$,
\begin{align}
 \Big| (V^\star - \tilde V^{\bpi_\epsilon})(\mu) \Big|  &\leq \; \Big|  \Tc^{\bpi_\epsilon} (V^\star- \tilde V^{\bpi_\epsilon})(\mu)+\epsilon \Big| \; \leq \;  \beta \Vert V^\star-\tilde V^{\bpi_\epsilon}\Vert +\epsilon,
\end{align} 
where we used the $\beta$-contraction property of $\Tc^{\bpi_\epsilon}$ in Lemma \ref{lemTcpi}.   We deduce that  
$\Vert V^\star -\tilde V^{\bpi_\epsilon} \Vert$ $\leq$ $\frac{\epsilon}{1-\beta}$,  and then, $\tilde V$ $\geq$  $\tilde{V}^{\bpi_\epsilon}$ $\geq$ 
$V^\star-\frac{\epsilon}{1-\beta}$, which combined with $ V^\star\geq \tilde V$, shows the first assertion.  Moreover, if $\bpi_\eps$ $=$ $\bpi^{\mfa_\epsilon}$ is a lifted randomized feedback policy, then by Remark \ref{remlift}, and under 
{\bf Rand$(\xi,\Gc)$}, we have $V^{\mfa_\epsilon}(\xi)$ $=$ $\tilde V^{\bpi_\epsilon}(\mu)$. Recalling that $V(\xi)$ $\leq$ $\tilde V(\mu)$, and together with the first assertion, this proves the required result.  
\ep


\begin{Remark} \label{remDPP} 
{\rm  If we can find for any $\epsilon$ $>$ $0$, an $\epsilon$-optimal lifted randomized feedback policy for $V^\star$, then according to Proposition \ref{lemverif}, and under {\bf Rand$(\xi,\Gc)$}, 
one could restrict to randomized feedback policies in the computation of the  optimal value $V(\xi)$ of the CMKV-MDP,  i.e.,  $V(\xi)$ $=$ $\sup_{\mfa \in L^0(\Pc(\Xc)\times \Xc\times [0,1];A)}  V^{\mfa}(\xi)$. 
Moreover,  this would prove that $V(\xi)$ $=$ $\tilde V(\mu)$ $=$ $V^\star(\mu)$, hence $V$ is law-invariant, and satisfies the Bellman fixed equation. 

Notice that instead of proving directly the dynamic programming Bellman equation for $V$, we start from the fixed point solution $V^\star$  to the Bellman equation, and show via a verification result that $V$ is indeed equal to  $V^\star$, hence satisfies the Bellman equation. 

By the formulation \eqref{Belhat} of the Bellman operator in Proposition \ref{equivTc}, and the fixed point equation satisfied by $V^\star$, we know that for all $\epsilon$ $>$ $0$, and $\mu$ $\in$ $\Pc(\Xc)$, there exists 
$\mra_\epsilon^\mu$ $\in$ $L^0(\Xc\times [0,1];A)$ such that 
\begin{align} \label{VepsY} 
V^\star(\mu) & \leq [\T^{\mra_\epsilon^\mu} V^\star ](\mu) + \epsilon.
\end{align}
The crucial issue is to prove that the mapping $(\mu,x,u)$ $\mapsto$ $\mfa_\epsilon(\mu,x,u)$ $:=$ $\mra_\epsilon^\mu(x,u)$ is measu\-rable so that it defines a randomized feedback policy 
$\mfa_\epsilon$ $\in$ $L^0(\Pc(\Xc)\times \Xc\times [0,1];A)$, and an associated lifted randomized feedback policy $\bpi^{\mfa_\epsilon}$. Recalling the relation \eqref{lienT}, 
this would then show that $\bpi^{\mfa_\epsilon}$ is a $\epsilon$-optimal lifted randomized feedback policy for $V^\star$, and we could apply the verification result. 
}
\ep
\end{Remark}

\vspace{1mm}

We now address the measurability issue for proving the existence of an $\epsilon$-optimal randomized feedback policy for $V^\star$. The basic idea is to construct as in \eqref{VepsY} an 
$\epsilon$-optimal $\mra_\epsilon^\mu$ $\in$ $L^0(\Xc\times [0,1];A)$ for $V^\star(\mu)$ when $\mu$ lies in a suitable finite grid of $\Pc(\Xc)$,  
and then ``patchs" things together to obtain an $\eps$-optimal randomized feedback policy.  This is made possible under some uniform continuity 
property of $V^\star$.

 \vspace{1mm}

\vspace{5mm}

The next result provides a suitable discretization of the set of probability measures.

\begin{Lemma} \label{lemquantif} {\bf (Quantization of $\Pc(\Xc)$)}

\noindent Fix $\eta$ $>$ $0$. 
Then for each finite  $\eta/2$-covering $\Xc_\eta$ of $\Xc$, one can construct a finite subset  $\Mc_\eta$ of $\Pc(\Xc)$, of size $N_\eta$ $=$ $n_\eta^{\#\Xc_\eta - 1}$, where $n_\eta$ is a grid size of $[0,1]$, that is an $\eta$-covering of $\Pc(\Xc)$.
\end{Lemma}
{\bf Proof.}
As $\Xc$ is compact, there exists a finite subset $\Xc_\eta\subset \Xc$ such that $d(x,x_\eta)\leq \eta/2$ for all $x\in\Xc$, where $x_\eta$ denotes the projection of $x$ on 
$\Xc_\eta$. Given  $\mu\in\Pc(\Xc)$, and $\xi$ $\sim$ $\mu$, we denote by $\xi_\eta$ the quantization, i.e., the projection of $\xi$ on $\Xc_\eta$,  and by 
$\mu_\eta$ the discrete law of $\xi_\eta$. Thus, $\E[d(\xi,\xi_\eta)]$ $\leq$ $\eta/2$, and therefore $\Wc(\mu,\mu_\eta)\leq \eta/2$.   The probability measure $\mu_\eta$ lies 
in $\Pc(\Xc_\eta)$, which  is identified with the simplex of $[0,1]^{\#\Xc_\eta}$.   We then use another grid $G_\eta$ $=$ 
$\{\frac{i}{n_\eta}: i=0,\ldots,n_\eta\}$ of $[0,1]$, and project its weights $\mu_\eta(y)$ $\in$ $[0,1]$, $y \in \Xc_\eta$, on $G_\eta$, in order to obtain another discrete probability measure $\mu_{\eta,n_\eta}$. 
From the dual Kantorovich representation of Wasserstein distance, it is easy to see  that for $n_\eta$ large enough, $\Wc(\mu_\eta, \mu_{\eta, n_\eta})$ $\leq$ $\eta/2$, and so $\Wc(\mu, \mu_{\eta, n_\eta})$ $\leq$ $\eta$. 
We conclude the proof by noting that   $\mu_{\eta,n_\eta}$  belongs to the set $\Mc_\eta$ of probability measures on $\Xc_\eta$  with weights valued in the finite grid $G_\eta$, hence $\Mc_\eta$ is a finite set of $\Pc(\Xc_\eta)$, 
of cardinal $N_\eta$ $=$ $n_\eta^{\#\Xc_\eta - 1}$.  
\ep

\begin{Remark}
{\rm
Lemma \ref{lemquantif} is actually a simple consequence of Prokhorov's theorem, but the above proof has some advantages:
\begin{itemize}
    \item it provides an explicit construction of the quantization grid $\Mc_\eta$,
    \item it explicitly gives the size of the grid as a function of $\eta$, which is particularly useful for computing the time/space complexity of algorithms,
    \item this special grid simultaneously quantizes $\Pc(\Xc)$ and $\Xc$, in the sense that the measures from $\Mc_\eta$ are all supported on the finite set $\Xc_\eta$. This is also useful for algorithms because for $\mu\in \Mc_\eta$, $ \hat a \in \hat A(\Xc)$, and $W\in L^\infty(\R)$, the expression $\hat{\Tc}^{\hat a} W(\mu)$ only depends upon $(\hat a(x))_{x\in\Xc_\eta}$. Therefore, in the Bellman fixed point equation, the computation of an $\epsilon$-argmax over the set $L
   ^0(\Xc,\Pc(A))$ is reduced to a computation of an $\epsilon$-argmax over the set $\Pc(A)^{\Xc_\eta}$, which is more tractable for a computer. 
\end{itemize} 
}
\end{Remark}

\vspace{5mm}

We can conclude this paragraph by showing the existence of an $\epsilon$-optimal lifted randomized feedback policy for the general lifted MDP on $\Pc(\Xc)$, and obtain as a by-product the corresponding Bellman fixed point equation for its value function and for the optimal value of the CMKV-MDP under randomization hypothesis.

\begin{Theorem} \label{theomainDPP} 
Assume that  $({\bf H_{lip}})$  holds true.   Then, for all $\epsilon>0$,  there exists a lifted randomized feedback policy $\bpi^{\mfa_\epsilon}$, for some $\mfa_\epsilon$ $\in$ $L^0(\Pc(\Xc)\times\Xc\times [0,1];A)$, that is 
$\epsilon$-optimal for $V^\star$.   Consequently, under  {\bf Rand$(\xi,\Gc)$}, the randomized feedback stationary control $\alpha^{\mfa_\epsilon}$ $\in$ $\Ac$ is $\frac{\epsilon}{1-\beta}$-optimal for $V(\xi)$, and we have 
$V(\xi)$ $=$ $\tilde V(\mu)$ $=$ $V^\star(\mu)$, for $\mu$ $=$ $\Lc(\xi)$, which thus satisfies the Bellman fixed point equation.   
\end{Theorem}
{\bf Proof.}  Fix $\epsilon$ $>$ $0$, and given $\eta$ $>$ $0$, consider a quantizing grid $\Mc_\eta$ $=$ $\{\mu^1,\ldots,\mu^{N_\eta}\}$ $\subset$ $\Pc(\Xc)$ as in Lemma \ref{lemquantif}, and an associated partition  $C_\eta^i$, $i$ $=$ $1,\ldots,N_\eta$, of $\Pc(\Xc)$, satisfying 
\begin{align}
C_\eta^i & \subset \; B_\eta(\mu^i) := \Big\{ \mu \in \Pc(\Xc): \Wc(\mu,\mu^i) \leq \eta \Big\}, \quad i=1,\ldots,N_\eta. 
\end{align}
For any $\mu^i$, $i$ $=$ $1,\ldots,N_\eta$, and  by \eqref{VepsY}, there exists $\mra_\epsilon^i$ $\in$ $L^0(\Xc\times [0,1];A)$ such that 
\begin{align} \label{VepsYi} 
V^\star(\mu^i) & \leq [\T^{\mra_\epsilon^i} V^\star](\mu^i) + \frac{\epsilon}{3}.  
\end{align}
From the partition $C^i_\eta$, $i$ $=$ $1,\ldots,N_\eta$ of $\Pc(\Xc)$, associated to $\Mc_\eta$,  
we  construct the function $\mfa_\epsilon$ $:$ $\Pc(\Xc)\times\Xc\times [0,1]$ $\rightarrow$ $A$ as follows. Let $h_1$, $h_2$ be two measurable functions from $[0,1]$ into $[0,1]$, 
such that if $U\sim\Uc([0,1])$, then $(h_1(U),h_2(U))\sim \Uc([0,1])^{\otimes 2}$.  We then define
\beqs
\mfa_\epsilon(\mu,x,u) &=& \mra_\epsilon^i\big(\zeta(\mu,\mu^i,x,h_1(u)),h_2(u)\big), \;\; \mbox{ when } \mu \in  C_\eta^i, \; i =1,\ldots,N_\eta, \; x \in \Xc, u \in [0,1], 
\enqs
where $\zeta$ is the measurable coupling function defined in Lemma \ref{lemcoupling1}. Such function  $\mfa_\eps$  is clearly measurable, i.e., 
$\mfa_\eps$ $\in$ $L^0(\Pc(\Xc)\times\Xc\times [0,1];A)$, and we  denote by $\bpi_\eps$ $=$ $\bpi^{\mfa_\eps}$ the associated lifted randomized feedback policy, which satisfies 
\begin{align} \label{TcT}
[\Tc^{\bpi_\eps} V^\star](\mu^i) &= \; [\T^{\mra_\epsilon^i} V^\star](\mu^i), \quad i=1,\ldots,N_\eta, 
\end{align}
by \eqref{lienT}. 
Let us now check that such $\bpi_\epsilon$ yields an $\epsilon$-optimal randomized feedback policy for $\eta$ small enough.  
For $\mu$ $\in$ $\Pc(\Xc)$, with $(\xi,U)$ $\sim$ $\mu\otimes\Uc([0,1])$, we set $U_1$ $:=$ $h_1(U)$, $U_2$ $:=$ $h_2(U)$,  
and define $\mu_\eta$ $=$ $\mu^i$,  when $\mu$ $\in$ $C_\eta^i$, $i=1,\ldots,N_\eta$,  and $\xi_\eta$ $:=$ $\zeta(\mu,\mu_\eta,\xi,U_1)$. 
Observe by Lemma \ref{lemquantif} that $\Wc(\mu,\mu_\eta)$ $\leq$ $\eta$, and  by Lemma \ref{lemcoupling1} that $(\xi_\eta,U_2)$ $\sim$ $\mu_\eta\otimes\Uc([0,1])$. 
We then write for any $\mu$ $\in$ $\Pc(\Xc)$, 
 \begin{align}
[\Tc^{\bpi_\epsilon} V^\star](\mu) - V^\star(\mu) &= \; 
\Big( [\Tc^{\bpi_\epsilon}  V^\star](\mu)  - [\Tc^{\bpi_\epsilon} V^\star](\mu_\eta) \Big) + \Big([\Tc^{\bpi_\epsilon} V^\star](\mu_\eta) - V^\star(\mu_\eta)\Big)  \\
&  \quad\quad    + \big(V^\star(\mu_\eta) -  V^\star(\mu)\big)  \\
& \geq  \Big( [\Tc^{\bpi_\epsilon} V^\star](\mu)-[\Tc^{\bpi_\epsilon} V^\star](\mu_\eta) \Big) - \frac{\epsilon}{3} -  \frac{\epsilon}{3},  \label{interTC} 
\end{align}
where we used \eqref{VepsYi}-\eqref{TcT} and the fact that  $|V^\star(\mu_\eta) -  V^\star(\mu)|$ $\leq$ $\epsilon/3$ for $\eta$ small enough by uniform continuity of $V^\star$ in Proposition \ref{proregul}. 
Moreover, by observing that $\mfa_\epsilon(\mu,\xi,U)$ $=$ $\mfa_\epsilon(\mu_\eta,\xi_\eta,U_2)$ $=:$ $\alpha_0$, so that $\bpi_\epsilon(\mu)$ $=$ $\Lc(\xi,\alpha_0)$, $\bpi_\epsilon(\mu_\eta)$ $=$ $\Lc(\xi_\eta,\alpha_0)$, 
we have 
\begin{align}
[\Tc^{\bpi_\epsilon} V^\star](\mu) &= \; \E \Big[ f(Y) + \beta V^\star(\P^0_{F(Y,\eps_1,\eps_1^0}) \Big], \\
[\Tc^{\bpi_\epsilon} V^\star](\mu_\eta) &= \; \E \Big[ f(Y_\eta) + \beta V^\star(\P^0_{F(Y_\eta,\eps_1,\eps_1^0}) \Big], 
\end{align}
where $Y$ $=$ $(\xi,\alpha_0,\bpi_\epsilon(\mu))$, and $Y_\eta$ $=$ $(\xi_\eta,\alpha_0,\bpi_\epsilon(\mu_\eta))$. Under  $({\bf H_{lip}})$, by using  the 
$\gamma$-H\"older property of $V^\star$ with constant $K_\star$ in Proposition \ref{proregul}, and by definition of the Wasserstein distance (recall that $\xi$ $\sim$ $\mu$, $\xi_\eta$ 
$\sim$ $\mu_\eta$), we then get 
\beqs
\big \vert [\Tc^{\bpi_\epsilon} V^\star](\mu)- [\Tc^{\bpi_\epsilon} V^\star](\mu _\eta)\big\vert &\leq & 2K \E\big[ d(\xi,\xi_\eta) \big] \\
& &   + \beta K_\star \E \Big[ 
\E\big[ d\big(F(\xi,\alpha_0,\bpi_\epsilon(\mu),\eps_1,e), F(\xi_\eta,\alpha_0,\bpi_\epsilon(\mu_\eta),\eps_1,e)\big)^\gamma\big]_{e:=\eps^0_1} \Big]  \\
&\leq & 2K \E\big[ d(\xi,\xi_\eta) \big] \\
& &   + \beta K_\star \E \Big[ 
\E\big[ d\big(F(\xi,\alpha_0,\bpi_\epsilon(\mu),\eps_1,e), F(\xi_\eta,\alpha_0,\bpi_\epsilon(\mu_\eta),\eps_1,e)\big)\big]_{e:=\eps^0_1} \Big]^\gamma  \\
& \leq  & C\E\big[ d(\xi,\xi_\eta) \big]^\gamma.  
\enqs
for some constant $C$ independent from $\E\big[ d(\xi,\xi_\eta) \big]$. Now, by the coupling Lemma \ref{lemcoupling1}, one can choose $\eta$ small enough so that  $C\E\big[ d(\xi,\xi_\eta) \big]^\gamma\leq \frac{\epsilon}{3}$. Therefore,  $\big\vert [\Tc^{\bpi_\epsilon} V^\star](\mu)-[\Tc^{\bpi_\epsilon} V^\star](\mu _\eta)\big\vert$ $\leq$ $\epsilon/3$, and, plugging into \eqref{interTC},  
we obtain $\Tc^{\bpi_\epsilon} V^\star(\mu) - V^\star(\mu)$ $\geq$ $-\epsilon$, for all $\mu\in\Pc(\Xc)$, which means  that $\bpi_\epsilon$ is $\epsilon$-optimal for $V^\star$. 
The rest of the assertions in the Theorem follows  from the verification result in Proposition \ref{lemverif}.  
\ep

\begin{Remark}
{\rm
We stress the importance of the coupling Lemma in the construction of $\epsilon$-optimal control in Theorem \ref{theomainDPP}. 
Indeed, as we do not make any regularity assumption on $F$ and $f$  with respect to  the ``control arguments", the only way  to make  $[\Tc^{\bpi_\epsilon} V^\star](\mu)$ and 
$[\Tc^{\bpi_\epsilon} V^\star](\mu_\eta)$ close to each other  is to couple terms to have the same control  in $F$ and $f$. This is achieved by turning $\mu$ into $\mu_\eta$, $\xi$ into $\xi_\eta$ and set 
$\alpha_0$ $=$ $\mfa_\epsilon(\mu,\xi,U)$ $=$ $\mfa_\epsilon(\mu_\eta,\xi_\eta,U_2)$. Turning $\mu$ into $\mu_\eta$ is a simple quantization, but turning $\xi$ into $\xi_\eta$ is obtained  
thanks to the coupling Lemma.
}
\ep
\end{Remark}

\begin{Remark}
{\rm
Theorem \ref{theomainDPP}, although applying to a more general case than the results from Section \ref{seclift}, provides a weaker result. Indeed, it does not state that any control $\bnu$ for the lifted MDP $\tilde{V}(\mu)$ can be represented, i.e.,  associated to a control $\alpha$ for $V(\xi)$ such that $\balpha_t:=\P^0_{(X_t,\alpha_t)}=\bnu_t$ for all $t\in\N$. This theorem only implies that one can {\em restrict}  the optimization to representable controls without changing the optimal value. Consequently, contrarily to Theorem \ref{theonocommon} and Theorem \ref{theodiscret}, here one cannot  conclude that an optimal control for $V(\xi)$ exists iff an optimal control for $\tilde{V}(\mu)$ exists. More precisely, it is possible that an optimal control $\bnu$ for $\tilde{V}(\mu)$ exists but cannot be associated to a control $\alpha$ for $V(\xi)$ such that $\balpha=\bnu$, and thus the existence of an optimal control for $\tilde{V}(\mu)$ does not guarantee the existence of an optimal control for $V(\xi)$.
}
\ep
\end{Remark}

\subsection{Relaxing the randomization hypothesis} 


We now relax the randomization hypothesis, and make the weaker assumption (recall Remark \ref{remrand}) that the initial information  filtration  $\Gc$ is  rich enough.

 \vspace{1mm}
 
 We need to state some uniform continuity property on the value function $V$ of the CMKV-MDP. 

\begin{Lemma} \label{unifV} 
Assume that  $({\bf H_{lip}})$  holds true.  Then, there exists a constant $K$ (only depending upon $K_F, \beta$, and $\Delta_\Xc$) such  that
\begin{align}
\sup_{\alpha\in\Ac} \big| V^\alpha(\xi) - V^\alpha(\xi') \big| & \leq \; K_\gamma 
 \big( \E\big[ d(\xi,\xi') \big] \big)^\gamma, \quad \forall \xi,\xi' \in  L^0(\Gc;\Xc),
\end{align}
where $\gamma=\min\big(1,\frac{\vert \ln(\beta)\vert}{(\ln 2K)_+}  \big)$.
Consequently, $V$ is $\gamma$-H\"older  on $L^0(\Gc;\Xc)$ endowed with the $L^1$-distance. 
\end{Lemma}
{\bf Proof.}
Fix $\xi, \xi' \in L^0(\Omega, \Xc)$, and consider an arbitrary $\alpha$ $=$ $\alpha^\pi$ $\in$ $\Ac$ associated to an open-loop policy $\pi\in\Pi_{OL}$. 
By Proposition \ref{condlaw-adapted}, there exists a measurable function $f_{t,\pi}\in L^0(\Xc\times G\times E^t\times (E^0)^t;\Xc)$, 
s.t. $X^{\xi,\alpha}_t= f_{t,\pi}(\xi,\Gamma, (\eps_s)_{s\leq t}, (\eps^0_s)_{s\leq t}) $, thus $\P^0_{X^{\xi,\alpha}_t}=\Lc(f_{t,\pi}(\xi,\Gamma, (\eps_s)_{s\leq t}, (e^0_s)_{s\leq t}))_{e^0_s=\eps^0_s,s\leq t}$. 
We thus have
\beqs
\Wc(\P^0_{(X^{\xi,\alpha}_t,\alpha_t)}, \P^0_{(X^{\xi',\alpha}_t,\alpha_t)}) \leq \E \Big[d\Big(f_{t,\pi}(\xi,\Gamma, (\eps_s)_{s\leq t}, (e^0_s)_{s\leq t}), f_{t,\pi}(\xi',\Gamma, (\eps_s)_{s\leq t}, (e^0_s)_{s\leq t})\Big)\Big]_{e^0_s=\eps^0_s,s\leq t},
\enqs
and so 
\begin{align}\label{wassdist}
\E \Big[ \Wc(\P^0_{(X^{\xi,\alpha}_t,\alpha_t)}, \P^0_{(X^{\xi',\alpha}_t,\alpha_t)}) \Big]  &\leq \;  \E \big[ d(X^{\xi,\alpha}_t, X^{\xi',\alpha}_t) \big],	
\end{align}
Under the Lipschitz condition on $f$ in  $({\bf H_{lip}})$, we then have
\begin{align} \label{Vlip}
\vert V^\alpha(\xi)-V^\alpha(\xi')\vert 
& \leq \;  2K_f \sum_{t=0}^\infty \beta^t\E \big[ d(X^{\xi,\alpha}_t, X^{\xi',\alpha}_t)\big].
\end{align}
By conditioning, and from the transition dynamics of the state process, we see that  for $t\in\N$, 
$\E\big[d(X^{\xi,\alpha}_{t+1}, X^{\xi',\alpha}_{t+1})\big]$ $=$ $\E\big[\Delta(\alpha^i_t,X^{\xi,\alpha}_{t},\P^0_{(X^{\xi,\alpha}_{t},\alpha_t)},  X^{\xi',\alpha}_{t},\P^0_{(X^{\xi',\alpha}_{t},\alpha_t)}, \eps^0_{t+1})\big]$,  where 
 \beqs
 \Delta(a,x,\nu,x',\nu',e^0_{t+1}) &=&  \E \big[d (F(x, a, \nu,\eps_{t+1},e^0_{t+1}), F(x', a, \nu',\eps_{t+1},e^0_{t+1}))\big].
 \enqs
 By the Lipschitz condition on $F$ in $({\bf H_{lip}})$, we thus have
 \beqs
 \E\big[ d(X^{\xi,\alpha}_{t+1}, X^{\xi',\alpha}_{t+1})\big] & \leq &  K_F\E [d(X^{\xi,\alpha}_{t}, X^{\xi',\alpha}_{t})+ \Wc(\P^0_{(X^{\xi,\alpha}_t,\alpha_t)}, \P^0_{(X^{\xi',\alpha}_t,\alpha_t)})] \; \leq \; 2K_F\E\big[ d(X^{\xi,\alpha}_{t}, X^{\xi',\alpha}_{t}) \big], 
 \enqs
 where the last inequality comes from \eqref{wassdist}.  Denoting by $\delta_t(\xi,\xi')$ $:=$ $\sup_{\alpha\in\Ac} \E\big[ d(X_t^{\xi,\alpha},X_t^{\xi',\alpha})\big]$, and noting that $\delta_0(\xi,\xi')$ $=$ $\E d(\xi,\xi')$, 
 it follows by induction that there exists $K\in\R$ such that 
 \begin{align}
 \delta_t(\xi,\xi') & \leq \; s_t(\E d(\xi, \xi')), \quad s_t(m) := m (2K)^t, \; m \geq 0, \; t \in \N.  
 \end{align}
 By the same arguments as in Theorem \ref{theoL1chaos},  and choosing $\gamma$ as in the assertion of the theorem,  we obtain that there exists $K$ (only depending upon $K_F, \beta$, and $\Delta_\Xc$) such that
 \begin{align}
 \sum_{t=0}^\infty \beta^t  \delta_t(\xi,\xi')  & \leq \; 
 K   \big( \E\big[ d(\xi,\xi') \big] \big)^\gamma, 
 \quad \xi,\xi' \in  L^0(\Gc;\Xc), 
 \end{align}
 and we conclude with \eqref{Vlip}. 
\ep

\vspace{3mm}



\begin{Theorem} \label{theomainrelax}
Assume that $\Gc$ is rich enough and $({\bf H_{lip}})$  holds true.  Then, for any $\xi$ $\in$ $L^0(\Gc;\Xc)$, $V(\xi)$ $=$ $\tilde {V}(\mu)$, where $\mu$ $=$ $\Lc(\xi)$. 
Consequently,  $V$ is law-invariant, identified with $\tilde V$, and satisfies the Bellman fixed point equation $\tilde V$ $=$ $\Tc\tilde V$.  Moreover, for all $\epsilon$ $>$ $0$, there exists an $\epsilon$-optimal randomized feedback control for $V(\xi)$. 
\end{Theorem}
{\bf Proof.}  As $\Xc$ is compact, there exists a finite subset $\Xc_\eta\subset \Xc$ such that $d(x,x_\eta)\leq \eta$ for all $x\in\Xc$, where 
$x_\eta$ denotes the projection of $x$ on $\Xc_\eta$.  
Fix $\xi$ $\in$ $L^0(\Gc;\Xc)$, and set $\mu$ $=$ $\Lc(\xi)$. Let us then consider a random variable $\xi'\sim\mu$  defined on another probability universe along with an independent uniform law $U'$. 
We set $\Gamma':=(\xi',U')$, and $\Gc'$ $=$ $\sigma(\Gamma')$.  By construction the randomization hypothesis  {\bf Rand$(\xi',\Gc')$} holds true, and we then have $V(\xi')=\tilde{V}(\mu)$ from Theorem \ref{theomainDPP}. 
Consider now the quantized random variables $\xi_\eta$ and $\xi'_\eta$, which have the same law, and satisfy respectively the randomization hypothesis {\bf Rand$(\xi_\eta,\Gc)$} and {\bf Rand$(\xi'_\eta,\Gc')$} from Lemma \ref{lemrandom}.  From Theorem \ref{theomainDPP}, 
we deduce that  $V(\xi_\eta)$ $=$ $\tilde{V}(\Lc(\xi_\eta))$ $=$ $V(\xi'_\eta)$. By uniform continuity of $V$, it follows by sending  $\eta$ to zero, that  
$V(\xi)$ $=$ $V(\xi')$, and thus $V(\xi)=\tilde{V}(\mu)$, which proves the required result.

Finally, the existence of an $\epsilon$-optimal control for $V(\xi)$ is obtained as follows. From the uniform continuity of $V$ in Lemma \ref{unifV}, there exists $\eta$ small enough so that  
$|V(\xi) - V(\xi_\eta)|$ $\leq$ $\epsilon/2$. We then build according to Theorem \ref{theomainDPP}  an $\epsilon/2$-optimal control for $V(\xi_\eta)$, which yields an $\epsilon$-optimal (randomized feedback stationary) control for $V(\xi)$. 
\ep

\vspace{2mm}


\begin{Remark} \label{remTca} 
{\rm From Theorems \ref{theomainDPP} and \ref{theomainrelax}, under the condition that $\Gc$ is rich enough and $({\bf H_{lip}})$  holds true, the value function $V$ of the CMKV-MDP is law-invariant,  and the supremum in the Bellman fixed point equation for $V$ $\equiv$ $\tilde V$  with the operator $\Tc$ can be restricted to lifted randomized feedback policies, i.e., 
\begin{align} \label{BelV2}
V &= \;  \Tc V \; = \; \sup_{\mfa \in L^0(\Pc(\Xc)\times\Xc\times [0,1];A)} \Tc^\mfa V
\end{align} 
where we set $\Tc^\mfa$ $:=$ $\Tc^{\bpi^\mfa}$ equal to 
\beqs
[\Tc^\mfa W](\mu) & = &  \E \Big[ f(Y^\mfa(\mu,\xi,U)) + \beta W(\P^0_{F(Y^\mfa(\mu,\xi,U),\eps_1,\eps_1^0}) \Big],
\enqs
with  $Y^\mfa(\mu,x,u)$ $:=$ $(x,\mfa(\mu,x,u),\bpi^\mfa(\mu))$,  and  $(\xi,U)$ $\sim$ $\mu\otimes\Uc([0,1])$. Notice that this Bellman fixed point equation is not the same as the Bellman fixed point equation obtained by optimizing over feedback controls only (not randomized nor open-loop). Let us call $V_f$ the associated optimal value. Then it is known that
\begin{align} \label{BelV2}
V_f &= \;  \Tc_f V_f \; = \; \sup_{\mfa_f \in L^0(\Pc(\Xc)\times\Xc;A)} \Tc^{\mfa_f} V.
\end{align} 
In other words, in the feedback case, the sup in the Bellman fixed point equation is only taken over (non-randomized) feedback policies.
}
\ep
\end{Remark}

\subsection{Feedback vs open-loop vs randomized controls}

In this paper, we have  mentioned different types of controls: open-loop controls, feedback controls, and randomized feedback controls.  
To fix ideas, let us consider a CMKV-MDP where {\bf Rand$(\xi,\Gc)$} does not initially hold true, and let us address three problems:
\begin{itemize}
    \item {\bf Feedback problem}: Optimizing over stationary feedback controls. We note $V_f$ the corresponding optimal value.
    \item {\bf Open-loop problem}: Optimizing over open-loop controls. We note $V_{ol}$ the corresponding optimal value.
    \item {\bf Randomized feedback problem}: Enlarging $\Gc$ with an independent uniform random variable $U$, i.e. $\tilde{\Gc}:=\sigma(\Gc,U)$, so that {\bf Rand$(\xi,\tilde{\Gc})$} holds true, and then optimizing over stationary randomized feedback controls. We note $V_r$ the corresponding optimal value.
\end{itemize}

We shall compare these problems by relying on the results of this paper as well as the examples presented at the end of this section. We will focus on the following questions:

\begin{itemize}
    \item[(1)]{\bf When is Bellman fixed point equation satisfied?} Theorem \ref{theomainDPP} shows that $V_r$ satisfies a Bellman equation, and we have already mentioned in Remark \ref{remTca} that $V_f$ does satisfy a Bellman fixed point equation as well (although not the same). Regarding the open-loop problem, we know that if $\Gc$ is rich enough, $V_{ol}$ satisfies the same Bellman equation as $V_r$. However, when $\Gc$ is not rich enough, Example \ref{example-OL-disc} below  provides a counter-example, where the Bellman fixed point equation is not satisfied. 
    \item[(2)] {\bf When do the optimal values coincide?} Theorem \ref{theomainDPP} shows that $V_r$ is the same as the optimal value when the optimization is performed over open-loop controls (with initial information $\tilde{\Gc}$), and thus $V_r(\xi)\geq V_{ol}(\xi)$ (because $\Gc\subset\tilde{\Gc}$). Also, we clearly have $V_{ol}(\xi)\geq V_f(\xi)$. The problem is now to figure out if the inequalities can be strict. If $\Gc$ is not rich enough, Examples \ref{example-F-disc}, \ref{example-OL-disc}, and \ref{example-R-disc} below  illustrate that one can have $V_r(\xi)>V_{ol}(\xi)>V_f(\xi)$. When $\Gc$ is rich enough, Theorem \ref{theomainrelax} shows that $V_r(\xi)=V_{ol}(\xi)$. However, regarding the comparison with $V_f$, it is remarkable that Examples \ref{example-F-cont} and \ref{example-R-cont-1} below illustrate that the inequality $V_{r}=V_{ol}>V_f$ can be strict even when $\Gc$ is rich enough, and actually even when $\xi$ is atomless.
    \item[(3)] {\bf When does an optimal control exist?} Given that in the comparison between the randomized feedback problem and the open-loop problem, we have seen that the optimal values are equal as soon as $\Gc$ is rich enough, it is natural to wonder if there is still a qualitative difference between both cases. It is notable that Examples \ref{example-R-cont-2} and \ref{example-OL-cont} below  illustrate that an optimal control can exist in the randomized feedback case but not in the open-loop case. Furthermore, the optimal control from the randomized feedback case can have a much simpler and more regular form than any $\epsilon$-optimal control from the open-loop case.
\end{itemize}

We conclude this section by a set of simple examples illustrating the subtle differences between each type of problem. 
\begin{Example}[Feedback problem, $\Gc$ not rich enough]\label{example-F-disc}
Let us take an example similar to Example 3.1 in \cite{guetal19}.   
Consider $\Xc=\{-1,1\}=A$, $\eps_1\sim\Bc(1/2)$, $F(x,a,\nu,e,e^0)= ax$, $f(x,a,\nu)=-\Wc(\text{pr}_1\star \nu, \Bc(1/2))$. In other words, the reward is maximal and equal to $0$ when the law of the state is a Bernoulli$(1/2)$ on $\Xc$, and minimal equal to $-1/2$ when the law of the state is a Dirac ($\delta_{-1}$ or $\delta_{1}$). Assume that $\Gamma=1$ a.s., so that  that $\Gc$ is the trivial $\sigma$-algebra. In this case, $\xi=:x$ is necessarily deterministic. We perform the optimization over feedback controls. It is clear that the law of $X_t$ will always be a Dirac, and thus the gain will be $V_f(\xi)=\sum_{t=0}^\infty \beta^t (- \frac{1}{2})=-\frac{1}{2(1-\beta)}$ which is the worst possible gain.
\end{Example}

\begin{Example}[Open-loop problem, $\Gc$ not rich enough]\label{example-OL-disc}
 Let us consider the same problem as in Example \ref{example-F-disc}, but optimize over open-loop controls instead. $\Gc$ being the trivial $\sigma$-algebra, $\xi=:x$ and $\alpha_0$ are then necessarily deterministic, and thus $X^{\xi,\alpha}_1=\alpha_0\xi$ has to be deterministic as well, which yields rewards at $t=0,1$ both equal to $-\frac{1}{2}$. 
 By choosing a control $\alpha_0=1$, $\alpha_1=\eps_1$ and $\alpha_t=1$ afterwards, we have $\P^0_{X^{x,\alpha}_t}=\delta_x$ for $t=0,1$,  and $\P^0_{X^{x,\alpha}_t}=\Bc(1/2)$ afterwards. 
This control is clearly optimal,  and the associated gain is $V_{ol}(\xi)=-\frac{1+\beta}{2}$. Notice that in this example, where $\Gc$ is not rich enough, $V$ does not satisfy a Bellman equation. Indeed, if it did, we would have 
$V_{ol}(x)= \underset{a\in A}{\sup} (-1/2 +\beta V_{ol}(ax))$, which is equivalent to $-\frac{1+\beta}{2}=-\frac{1}{2} -\beta\frac{1+\beta}{2}$,  and this is clearly false.
\end{Example}

\begin{Example}[Randomized feedback problem, $\Gc$ not rich enough]\label{example-R-disc}
 Let us consider the same problem as in Example  \ref{example-F-disc}, but let us enlarge $\Gc$ with an independent uniform random variable $U$, i.e. $\tilde{\Gc}:=\sigma(\Gc,U)$, so that {\bf Rand$(\xi,\tilde{\Gc})$} holds true, and then optimizing over stationary randomized feedback controls. The enlarged initial information allows to set $\alpha_0=sgn(U-\frac{1}{2})$ and $\alpha_t=1$ for $t\in\N_\star$. It is clear that the strategy is optimal and leads to a gain $V_r(\xi)=-\frac{1}{2}$.
\end{Example}

\begin{Example}[Feedback problem, $\Gc$ rich enough, $\xi$ atomless]\label{example-F-cont}
Consider $\Xc=\{-1,1\}\cup [2,3]$, $A=\{-1,1\}$, $\eps_1\sim\Bc(1/2)$, $F(x,a,\nu,e,e^0)= a x {\bf 1}_{x\in \{-1,1\}}+{\bf 1}_{x\in [2,3]}$, $f(x,a,\nu)=-\Wc(\text{pr}_1\star \nu, \Bc(1/2))$. 
Assume that $\Gamma\sim \Uc([0,1]^2)$ and $\xi=2+\Gamma$ a.s.. In particular, $\Gc$ is rich enough and $\xi$ is atomless.  This is the same example as Example \ref{example-R-disc} except for the choice of $\Gamma$ and $\xi$, and for the fact that we have extended $F$ on $\Xc=\{-1,1\}\cup [2,3]$.
Assume that we optimize over feedback controls. Given the definition of $F$, it is clear that at time $t=0$, the control has no impact, and that $X_1=1$ a.s.. From $t=1$, we thus fall back into the problem from Example \ref{example-F-disc} and thus the optimal value is $V_f(\xi)=-\Wc(\Uc([2,3]),\Bc(1/2))-\beta\frac{1+\beta}{2}$.
\end{Example}

\begin{Example}[Open-loop problem, $\Gc$ rich enough, $\xi$ atomless]\label{example-R-cont-1}
Let us consider the same problem as in Example \ref{example-F-cont}, but let us optimize over open-loop controls instead. Again, at time $t=0$, the control has no impact, and $X_1=1$ a.s.. From $t=1$, we fall back into the problem from Example \ref{example-OL-disc} except that $\Gc$ is now rich enough. By Theorem \ref{theomainrelax}, the optimal value from time $t=1$ is thus the same as for the randomized feedback problem from Example \ref{example-R-disc}, and thus the optimal value is $V_{ol}(\xi)=-\Wc(\Uc([2,3]),\Bc(1/2))-\beta\frac{1}{2}$.
\end{Example}

\begin{Example}[Randomized feedback problem, existence of an optimal control]\label{example-R-cont-2}
Assume that $\Xc=[-1,1]$, $A=\{-1,1\}$, $\eps_1\sim \Uc([-1,1])$, $F(x,a,\nu,e,e^0)= ax$, $f(x,a,\nu)=-\Wc(\text{pr}_1\star \nu, \Uc([-1,1]))$,  $\Gamma=(U,W)\sim \Uc([0,1]^2)$, 
and $\xi= W$. Then {\bf Rand$(\xi,\Gc)$} holds true. Similarly to Example \ref{example-R-disc}, we set $\alpha_0=sgn(U-\frac{1}{2})$ and $\alpha_t=1$ for $t\in\N_\star$. It is straightforward to show that the strategy is optimal and also leads to a gain $V_r(\xi)=-\frac{1}{2}$.
\end{Example}

\begin{Example}[Open-loop problem, no optimal control]\label{example-OL-cont}
Let us take the same problem as in Example \ref{example-R-cont-2}, but with $\Gamma=W$. 
Then {\bf Rand$(\xi,\Gc)$} does not hold true anymore. Indeed, any $\Gc$-measurable variable independent from $\xi$ is necessarily deterministic and therefore cannot be uniform on $[0,1]$.
By Theorem \ref{theomainrelax}, as $\Gc$ is rich enough, we know that we will also have $V_{ol}(\xi)=-\frac{1}{2}$. However, without the additional uniform variable $U$, the existence of an optimal control does not hold true anymore. Indeed, from Example \ref{example-R-cont-2}, an optimal strategy would imply $X_1\sim\Uc([-1,1])$. 
On the other hand,  since $\alpha_0$ is  $\Gc$-measurable, it has to take the (non-randomized) form $\alpha_0=h(\xi)$ for $h\in L^0(\Xc,A)$.  The latter  implies  that the support of the distribution of  $X_1=h(\xi)\xi$ is of Lebesgue measure $1$, which is in contradiction with the $\Uc([-1,1])$-distribution of $X_1$. 
Moreover, from the proof of Theorem \ref{theomainrelax}, it is clear that building an approximate optimal control would require to introduce a quantized version $\xi_\eta$ of $\xi$ for $\eta$ small, and then to use Lemma \ref{lemrandom} to build a uniform variable $U_\eta\independent \xi_\eta$. Therefore, the approximate optimal control would be much more complex than the simple optimal control from Example \ref{example-R-cont-2}. 
\end{Example}

\subsection{Computing value function and  $\epsilon$-optimal strategies in CMKV-MDP}

Having established the correspondence of our CMKV-MDP with lifted MDP on $\Pc(\Xc)$, 
and the associated Bellman fixed point equation, we can (up to a simple discretization of the state space in the Bellman fixed point equation) design two methods for computing the value function and optimal strategies: 

\vspace{1mm}

\noindent {\bf (a) Value iteration.} We approximate the value function $V$ $=$ $\tilde V$ $=$ $V^\star$ by iteration from the Bellman operator: $V_{n+1}$ $=$ $\Tc V_n$, and at iteration $N$, we compute an approximate optimal  randomized feedback policy $\mfa_N$ by  (recall Remark \ref{remTca}) 
\begin{align}
\mfa_N & \in \; {\rm arg}\max_{\mfa\in L^0(\Pc(\Xc)\times\Xc\times [0,1];A)} \Tc^{\mfa}  V_N. 
\end{align}
From $\mfa_N$, we then construct an approximate randomized feedback stationary control $\alpha^{\mfa_N}$ according to the procedure described in Remark \ref{remlift}.  

\vspace{1mm}

\noindent {\bf (b) Policy iteration.} Starting from some initial randomized feedback policy $\mfa_0$ $\in$ $L^0(\Pc(\Xc)\times\Xc\times [0,1];A)$, we iterate according to: 
\begin{itemize}
\item Policy evaluation: we compute the expected gain $\tilde V^{\bpi^{\mfa_0}}$ of the lifted MDP
\item Greedy strategy: we compute 
\begin{align}
\mfa_{k+1}  & \in \; {\rm arg}\max_{\mfa\in L^0(\Pc(\Xc)\times\Xc\times [0,1];A)} \Tc^{\mfa}  \tilde V^{\bpi^{\mfa_k}}.  
\end{align}
\end{itemize}
We stop at iteration $K$ to obtain $\mfa_K$, and then construct an approximate randomized feedback control $\alpha^{\mfa_K}$ according to the procedure described in Remark \ref{remlift}.

\vspace{3mm} 
 
\noindent  {\bf Practical computation.}  Since a randomized feedback control $\alpha$  is a measurable function $\mfa$  of  $(\P^0_{X_t^{\xi,\alpha}}, X_t^{\xi,\alpha}, U_t)$, we would need to compute and store the 
(conditional) law of the state process, which is infeasible in practice when  $\Xc$ is a continuous space. 
In this case,  to circumvent this issue, a natural idea is to discretize the compact space $\Xc$ by considering 
a finite subset $\Xc_\eta$ $=$ $\{x^1,\ldots,x^{N_\eta}\}$ $\subset$ $\Xc$ associated with  a partition  $B_\eta^i$, $i$ $=$ $1,\ldots,N_\eta$,  of $\Xc$,  satisfying: 
$B_\eta^i$ $\subset$ $\big\{ x \in \Xc: d(x,x^i) \leq \eta \big\}$,  $i=1,\ldots,N_\eta$, with $\eta$ $>$ $0$. For any $x$ $\in$ $\Xc$, we denote by $[x]_\eta$ (or simply $x_\eta$) its projection on $\Xc_\eta$, defined by: 
$x_\eta$ $=$ $x^i$,   for $x \in B_\eta^i$, $i=1,\ldots,N_\eta$.



\begin{Definition}[Discretized CMKV-MDP] \label{defdiscret} 
Fix $\eta>0$. Given $\xi\in L^0(\Gc;\Xc_\eta)$, and a control $\alpha\in\Ac$, we denote by $X^{\eta,\xi,\alpha}$ the McKean-Vlasov MDP on $\Xc_\eta$ given by
\begin{align} \label{dynXMKV2dis} 
 X^{\eta,\xi,\alpha}_{t+1} & =   \big[ F(X^{\eta,\xi,\alpha}_t, \alpha^{}_t, \P^0_{(X^{\eta,\xi,\alpha}_t,\alpha_t)},\eps_{t+1}, \eps^0_{t+1}) \big]_\eta,  \quad 	t\in\N, \;\; X_0^{\eta,\xi,\alpha}  = \xi, 
\end{align}
i.e., obtained by projecting the state on $\Xc_\eta$ after each application of the transition function $F$. The associated expected gain $V_\eta^\alpha$ is defined by
\begin{align} \label{defValphadis}
V_\eta^\alpha(\xi) &= \E \Big[ \sum_{t=0}^\infty  \beta^t f\big(X^{\eta,\xi,\alpha}_t,\alpha_t,\P^0_{(X^{\eta,\xi,\alpha}_t,\alpha_t)}\big) \Big].  
\end{align}
\end{Definition}

Notice that the (conditional) law of the discretized CMKV-MDP on $\Xc_\eta$ is now valued in a finite-dimensional space (the simplex of $[0,1]^{N_\eta}$), which makes the computation of the associated 
randomized feedback control  accessible, although computationally challenging due to the high-dimensionality (and beyond the scope of this paper). The next result states that an $\epsilon$-optimal randomized feedback control in the initial CMKV-MDP can be approximated by a randomized feedback control in the discretized CMKV-MDP.

\begin{Proposition} \label{theodisc}
Assume that $\Gc$ is rich enough and $({\bf H_{lip}})$  holds true.  Fix $\xi\in L^0(\Gc;\Xc)$. Given $\eta>0$, let us define $\xi_\eta$ the projection of $\xi$ on $\Xc_\eta$. As {\bf Rand$(\xi_\eta,\Gc)$} holds true, let us consider 
an i.i.d. sequence $(U_{\eta,t})_{t\in\N}$ of $\Gc$-measurable uniform variables independent of $\xi_\eta$. For $\epsilon>0$, let $\mfa_\epsilon$ be a randomized feedback policy that is $\epsilon$-optimal for the Bellman fixed point equation satisfied by $V$. Finally, let $\alpha^{\eta,\epsilon}$ be the randomized feedback control in the discretized CMKV-MDP 
recursively defined by $\alpha_t^{\eta,\epsilon}$ $=$ 
$\mfa_\epsilon(\P^0_{X^{\eta,\epsilon}_{t}},X^{\eta,\epsilon}_{t},U_{\eta,t})$, $t$ $\in$ $\N$,  where we set $X_t^{\eta,\epsilon}$ $:=$ $X^{\eta,\xi_\eta,\alpha^{\epsilon,\eta}}_{t}$. 
Then the control $\alpha^{\eta,\epsilon}$ is $\Oc(\eta^\gamma+\epsilon)$-optimal for the CMKV-MDP $X$ with initial state $\xi$, where $\gamma =\min \big(1, \frac{\vert \ln \beta\vert}{(\ln 2K)_+}\big)$. 
\end{Proposition}
{\bf Proof.}  {\it Step 1.} Let us show that
\begin{align}\label{distdisc}
\underset{\alpha\in\Ac}{\sup}\sum_{t=0}^\infty \beta^t \E \big[ d(X^{\xi,\alpha}_t, X^{\eta,\xi_\eta,\alpha}_t) \big] & \leq  \; C \eta^\gamma, 
\end{align} 
for some constant $C$ that depends only on $K$, $\beta$ and $\gamma$. 
Indeed, notice by definition of the projection on $\Xc_\eta$, and by a simple conditioning argument  that for all $\alpha$ $\in$ $\Ac$, and $t$ $\in$ $\N$, 
\begin{align}
\E \big[  d(X^{\xi,\alpha}_{t+1}, X^{\eta,\xi_\eta,\alpha}_{t+1}) \big] & \leq \; \eta + 
\E \big[ \Delta\big( X^{\xi,\alpha}_{t}, X^{\eta, \xi_\eta,\alpha}_{t},\alpha_t, \P^0_{(X^{\xi,\alpha}_t,\alpha_t)}, \P^0_{(X^{\eta,\xi_\eta,\alpha}_t,\alpha_t)},\eps^0_{t+1} \big) \big], 
\end{align}
where 
\begin{align}
\Delta(x,x',a,\nu,\nu',e^0) &= \; \E[d(F(x,a,\nu,\eps_{t+1},e^0), F(x',a,\nu',\eps_{t+1},e^0))]. 
\end{align} 
Under {$({\bf H_{lip}})$, we then get
\begin{align}
\E \big[  d(X^{\xi,\alpha}_{t+1}, X^{\eta,\xi_\eta,\alpha}_{t+1}) \big] & \leq \; \eta + 
K  \E\Big[d(X^{\xi,\alpha}_{t}, X^{\eta, \xi_\eta,\alpha}_{t})+\Wc(\P^0_{X^{\xi,\alpha}_t}, \P^0_{X^{\eta,\xi_\eta,\alpha}_t})\Big] \\
& \leq \; \eta + 2K \E\big[d(X^{\xi,\alpha}_{t}, X^{\eta, \xi_\eta,\alpha}_{t}) \big],
\end{align} 
by the same argument as in \eqref{wassdist}.  Hence,  the sequence $(\E\big[d(X^{\xi,\alpha}_{t}, X^{\eta, \xi_\eta,\alpha}_{t})\big] )_{t\in\N}$ satisfies the same type of induction inequality as in \eqref{estid} in 
Theorem  \ref{theoL1chaos} with $\eta$ instead of $M_N$, and thus the same derivation leads to  the required result \eqref{distdisc}. From the Lipschitz condition on $f$, we deduce by the same 
arguments as in \eqref{wassdist} in  Lemma \ref{unifV} that 
\begin{align} \label{estiVetauni}
\sup_{\alpha\in\Ac} \big| V^\alpha(\xi_\eta) -V_\eta^\alpha(\xi_\eta) \big|
&= \; \Oc(\eta^\gamma). 
\end{align}

\vspace{2mm}

\noindent {\it Step 2.} Denote by $\mu$ $=$ $\Lc(\xi)$, and $\mu_\eta$ $=$ $\Lc(\xi_\eta)$, and observe that $\Wc(\mu,\mu_\eta)$ $\leq$ $\E[d(\xi,\xi_\eta)]$ $\leq$ $\eta$. We write 
\beqs
V^{\alpha^{\eta,\epsilon}}(\xi)-V(\xi) &=& \big[ V^{\alpha^{\eta,\epsilon}}(\xi) - V^{\alpha^{\eta,\epsilon}}(\xi_\eta) \big] + \big[ V^{\alpha^{\eta,\epsilon}}(\xi_\eta)-V^{\alpha^{\eta,\epsilon}}_\eta(\xi_\eta) \big] \\
& & \; + \;  \big[ V^{\alpha^{\eta,\epsilon}}_\eta(\xi_\eta)-V(\xi_\eta) \big] + \big[ V(\xi_\eta)-V(\xi) \big]  \; =: \;  I_1 + I_2 + I_3 + I_4. 
\enqs
The first and last terms $I_1$ and $I_4$  are smaller than $\Oc(\eta^\gamma)$ by the $\gamma$-H\"older property of $V^\alpha$ and $V$ in Lemma \ref{unifV}. By  \eqref{estiVetauni}, the second term $I_2$ is  of  order  
$\Oc(\eta^\gamma) $ as well for $\eta$ small enough. Regarding the third term $I_3$,  notice that by definition, $V^{\alpha^{\eta,\epsilon}}_\eta(\xi_\eta)$ corresponds to the gain associated to the randomized feedback policy $\mfa_\epsilon$ for the discretized CMKV-MDP. 
Denote by $\bpi_{\epsilon}$ the lifted randomized feedback policy associated to $\mfa_{\epsilon}$, and recall by Remark \ref{remlift} the identification with the lifted MDP:  
$V^{\alpha^{\eta,\epsilon}}_\eta(\xi')$ $=$ $\tilde V_\eta^{\bpi_{\epsilon}}(\mu')$, $\mu'$ $=$ $\Lc(\xi')$, where  $\tilde V_\eta^{\bpi_{\epsilon}}$ is the expected gain of the lifted MDP associated to the discretized CMKV-MDP, hence fixed point of the operator
\begin{align}
[\Tc_\eta^{\mfa_{\epsilon}} W](\mu') &= \; \E \Big[ f(Y^{\mfa_\epsilon}(\mu',\xi',U))  + \beta W\big(\P^0_{\big[F(Y^{\mfa_\epsilon}(\mu',\xi',U),\eps_1,\eps_1^0)\big]_\eta}\big)   \Big], 
\end{align}
$Y^\mfa(\mu,x,u)$ $=$ $(x,\mfa(\mu,x,u),\bpi^\mfa(\mu))$ and $(\xi',U)\sim\mu' \otimes\Uc([0,1])$.  Recalling that $V(\xi')$ $=$ $\tilde V(\mu')$, $\mu'$ $=$ $\Lc(\xi')$,  
with $\tilde V$ fixed point to the Bellman operator 
$\Tc$, it follows that 
\beqs
I_3 \; = \;  \tilde V_\eta^{\bpi_{\epsilon}}(\mu_\eta) - \tilde V(\mu_\eta) & = &  
\Big(  [\Tc_\eta^{\mfa_{\epsilon}} \tilde V_\eta^{\bpi_{\epsilon}}](\mu_\eta) - [\Tc_\eta^{\mfa_{\epsilon}} \tilde V](\mu_\eta)  \Big) + 
\Big( [\Tc_\eta^{\mfa_{\epsilon}} \tilde V](\mu_\eta) -  [\Tc^{\mfa_{\epsilon}} \tilde V](\mu_\eta) \Big) \\
& & \; + \; \Big( [\Tc^{\mfa_{\epsilon}} \tilde V](\mu_\eta) -   \tilde V(\mu_\eta) \Big) \; =: \; I_3^1 + I_3^2 + I_3^3.   
\enqs
By definition of $\mfa_\epsilon$,  we have $|I_3^3|$ $\leq$ $\epsilon$. 
For $I^2_3$ notice that the only difference between the operators $\Tc_\eta^{\mfa_{\epsilon}}$ and  $\Tc^{\mfa_{\epsilon}}$
is that  $F$ is projected on $\Xc_\eta$. Thus, 
\beqs
\Big| [\Tc_\eta^{\mfa_{\epsilon}} \tilde V](\mu_\eta) -  [\Tc^{\mfa_{\epsilon}} \tilde V](\mu_\eta) \Big|
 &\leq &  \beta \E \Big[ \big\vert \tilde V \big(\P^0_{[F(Y_\eta,\eps_1,\eps_1^0)]_\eta}\big)    - \tilde V\big(\P^0_{F(Y_\eta,\eps_1,\eps_1^0)}\big)   \big\vert \Big], 
\enqs  
where $Y_\eta$ $=$ $(\xi_\eta,\mfa_\epsilon(\mu,\xi_\eta,U),\bpi_\epsilon(\mu_\eta))$.  It is clear by definition of the Wasserstein distance and the projection on $\Xc_\eta$ that 
\beqs
\Wc\big(\P^0_{[F(Y_\eta,\eps_1,\eps_1^0)]_\eta},  \P^0_{F(Y_\eta,\eps_1,\eps_1^0)} \big) & \leq & \E^0[d(F(Y_\eta,\eps_1,\eps^0_1), [F(Y_\eta,\eps_1,\eps^0_1)]_\eta)] \; \leq \; \eta. 
\enqs
From the $\gamma$-H\"older property of $\tilde V$ in Proposition \ref{proregul}, we deduce that $I_3^2$ $=$  $\Oc(\eta^\gamma)$.
Finally, for $I_3^1$, since $\Tc_\eta^{\mfa_{\epsilon}}$ is a  $\beta$-contracting operator on $(L^\infty(\Mc_\eta),\|\cdot\|_{\eta,\infty})$, we have  
\begin{align}
\big|[\Tc_\eta^{\mfa_{\epsilon}} \tilde V_\eta^{\bpi_{\epsilon}}](\mu_\eta) - [\Tc_\eta^{\mfa_{\epsilon}} \tilde V](\mu_\eta)\big| & \leq \; 
\beta\Vert \tilde V^{\bpi_\epsilon}_\eta -\tilde V\Vert_{\eta,\infty}, 
\end{align}
and thus $\vert \tilde V_\eta^{\bpi_{\epsilon}}(\mu_\eta) - \tilde V(\mu_\eta) \vert$ $=$ $\vert I_3\vert$ 
$\leq$ $\vert I_3^1\vert + \vert I_3^2\vert + \vert I_3^3\vert$ $\leq$ $\beta\Vert \tilde V^{\bpi_\epsilon}_\eta -\tilde V\Vert_{\eta,\infty} +\Oc(\eta^\gamma+\epsilon)$. Taking the $\sup$ over $\mu_\eta\in\Mc_\eta$ on the left, we obtain  
that $\Vert \tilde V^{\bpi_\epsilon}_\eta -\tilde V\Vert_{\eta,\infty} \leq \frac{1}{1-\beta}\Oc(\eta^\gamma +\epsilon)=\Oc(\eta^\gamma+\epsilon) $, and we conclude that 
$\vert I_3\vert\leq \Vert \tilde V^{\bpi_\epsilon}_\eta -\tilde V\Vert_{\eta,\infty} \leq  \Oc(\eta^\gamma + \epsilon)$, which ends the proof. 
\ep

\begin{Remark}
\rm{ Back to the $N$-agent MDP problem with open-loop controls, recall from Section \ref{secchaos}, that it suffices to find an $\epsilon$-optimal open-loop policy $\pi^\epsilon\in\Pi_{OL}$ for the CMKV-MDP, 
as it will automatically be $\Oc(\epsilon)$-optimal for the $N$-agent MDP with $N$ large enough. For instance, the construction of an $\epsilon$-optimal control $\alpha^\epsilon$ given by 
Proposition \ref{theodisc} can be associated to an $\epsilon$-optimal open-loop policy $\pi^\epsilon$ such that $\alpha^\epsilon=\alpha^{\pi^\epsilon}$ (where $\pi_t^{\epsilon}$ is a measurable function of  
$(\Gamma,(\eps_s)_{s\leq s\leq t},(\eps^0_s)_{s\leq s\leq t})$). The processes $\Oc(\epsilon)$-optimal for the $i$-th agent  $\alpha^{\epsilon,i}_t$ is then the result of the same construction but with  
$(\Gamma^i, \eps^i, \eps^0)$ instead of $(\Gamma, \eps, \eps^0)$, i.e. replacing $\xi$ by $\xi^i$, $U_{\eta,t}$ by $U_{\eta,t}^i$, and $(\Gamma, \eps, \eps^0)$ by $(\Gamma^i, \eps^i, \eps^0)$ in  
Proposition \ref{theodisc}. Notice that this construction never requires the access to  the  individual's states $X^{i,N}$.
 }
\end{Remark}

\begin{Remark}[$Q$ function] \label{remQ} 
{\rm  In view of the Bellman fixed point equation satisfied by the value function $V$ of the CMKV-MDP in terms of randomized feedback policies, 
let us introduce the corresponding state-action value function $Q$ defined on $\Pc(\Xc)\times \hat A(\Xc)$ by 
\begin{align}
Q(\mu,\hat a) &=  [\hat\Tc^{\hat a} V](\mu) \; = \;  \hat f(\mu,\hat a) + \beta \E \big[ V\big( \hat F(\mu,\hat a,\eps_1^0)\big) \big],
\end{align}
From Proposition \ref{equivTc}, and since $V$ $=$ $\Tc V$,  we recover  the standard connection between the value function and the state-action value function, namely $V(\mu)$ $=$ 
$\sup_{\hat a \in \hat A(\Xc)} Q(\mu,\hat a)$, from which we obtain the Bellman equation for the $Q$ function:
\begin{align} \label{belQ}
Q(\mu,\hat a) &=  \hat f(\mu,\hat a) + \beta \E \Big[ \sup_{\hat a' \in \hat A(\Xc)} Q\big( \mu_1^{\hat a},\hat a'\big) \Big],
\end{align}
 where we set $\mu_1^{\hat a}$ $=$ $\hat F(\mu,\hat a,\eps_1^0)$. Notice that this $Q$-Bellman equation extends the  equation in \cite{guetal19} (see their Theorem 3.1) derived in the 
 no common noise case and when there is no mean-field dependence with respect to the law of the control. The Bellman equation \eqref{belQ} is the starting point in a model-free framework when the state transition function is unknown (in other words in  the context of reinforcement learning) for the design of $Q$-learning algorithms in order to estimate  the $Q$-value function by 
 $Q_n$, and then to compute a relaxed control by
\begin{align}
\hat a_n^\mu & \in \; {\rm arg}\max_{\hat a\in \hat A(\Xc)} Q_n(\mu,\hat a), \quad \mu \in \Pc(\Xc). 
\end{align} 
From Lemma 2.22 \cite{Kallenberg}, one can associate to such probability kernel $\hat a_n^\mu$,  a function $\mfa_n$ $:$ 
$\Pc(\Xc)\times\Xc\times [0,1]$ $\rightarrow$ $A$, such that 
$\Lc(\mfa_n(\mu,x,U))$ $=$ $\hat a_n^\mu(x)$, $\mu$ $\in$ $\Pc(\Xc)$, $x$ $\in$ $\Xc$,  where $U$ is an uniform random variable.  
In practice, one has to discretize the state space $\Xc$ as in Definition \ref{defdiscret}, and then to quantize the space $\Pc(\Xc)$ as in Lemma \ref{lemquantif} in order to reduce the 
learning problem to a finite-dimensional problem for the computation of an approximate optimal randomized feedback policy $\mfa_n$ for the CMKV-MDP. 
}
\ep
\end{Remark}

\begin{Remark}[About the practical implementation of an RL algorithm] \label{remRL } 
{\rm  The individualized and open-loop form of controls has been motivated in the introduction from the point of view of applications. Let us now discuss some implications regarding the practical implementation of algorithms. As described in the introduction, the algorithms for computing the actions will be running in parallel on each individual's device (phone or computer), and in order to reduce the cost of data transfer, they will only access the data available on their respective devices. For instance, on device $i$, only the individual's own data $\eps^i$ and the public data $\eps^0$ will be accessible. These data are used to compute the McKean-Vlasov process $X^i$ associated to the individual. 
A randomized feedback policy is then applied to $X^i$ to compute the action $\alpha^i$. 

An important consequence is that each algorithm  must be able to compute the McKean-Vlasov process. Therefore:
(1) Either the dynamic $F$ and reward $f$ are explicitly modeled. 
(2) Either the population's behavior is simulated without explicitly describing $F$ and $f$.
The first case naturally leads to MDP algorithms (value/policy iteration), while the second case is better suited for RL algorithms. We refer to the updated book \cite{sutbar} for a general overview of RL algorithms, and to the more recent papers \cite{carlautan19} and \cite{guoetal19} in the framework of mean-field control/games.
}
\end{Remark}

\section{Conclusion} \label{secconclu} 

We have developed a theory for mean-field Markov decision processes with common noise and open-loop controls, called CMKV-MDP, for general state space and action space. 
Such problem is motivated and shown to be the asymptotic problem of a large population of cooperative agents under mean-field interaction controlled by a social planner/influencer, and we provide a rate of convergence of the 
$N$-agent model to the CMKV-MDP.  We prove the correspondence of CMKV-MDP with a general lifted MDP on the space of probability measures, and emphasize the role of relaxed control, which is crucial to characterize the solution via the Bellman fixed point equation. Appro\-ximate randomized feedback controls are obtained from the Bellman equation in a model-based framework, and future work under 
investigation will develop algorithms in a model-free framework, in other words in the context  of reinforcement learning with many interacting and cooperative agents.

\appendix

\section{Some useful results  on conditional law} \label{annexe_conditional}

\setcounter{Lemma}{0}
\setcounter{Definition}{0}

\begin{Lemma}\label{wassmes}
Let $(S,\Sc)$,  $(T,\Tc)$, and $(U,\Uc)$ be three measurable spaces, and $F\in L^0( (S,\Sc)\times (T,\Tc); (U,\Uc)) $ be a measurable function, then the function $\hat{F}:(\Pc(S),\Cc(S))\times (T,\Tc)\rightarrow (\Pc(U),\Cc(U))$ given by $\hat{F}(\mu, x):= F(\cdot,x)\star \mu$ is measurable.
\end{Lemma}
{\bf Proof.} This follows from the measurability of the maps: 
\begin{itemize}
	\item  $x\in(S,\Sc)\mapsto \delta_x\in(\Pc(S),\Cc(S))$, 
	\item  $(\mu,\nu)\in(\Pc(S),\Cc(S))\times (\Pc(T),\Cc(T))\mapsto \mu\otimes \nu\in(\Pc(S\times T),\Cc(S\times T)) $, 
	\item  $\mu\in(\Pc(S),\Cc(S))\mapsto F\star \mu \in (\Pc(T),\Cc(T)) $, 
	\end{itemize}
and the measurability of the  composition $(\mu,x)\mapsto(\mu,\delta_x)\mapsto\mu\otimes \delta_x\mapsto F\star (\mu\otimes \delta_x) = \hat{F}(\mu,x)$.
\ep

\vspace{3mm}


\begin{Lemma}[Conditional law]\label{condlaw}
Let $(S,\Sc)$ and   $(T,\Tc)$ be two measurable spaces. 
\begin{enumerate}
\item If $(S,\Sc)$ is a Borel space, there exists a conditional law of $Y$ knowing $X$.
\item If $Y=\varphi(X,Z)$ where $Z\independent X$ is a random variable  valued in a measurable space $V$ and $\varphi:S\times V\rightarrow T$
is a measurable function, then  $\Lc( \varphi(x,Z))\mid_{x=X}$ is a conditional law of $Y$ knowing $X$. 
In the case $S=S_1\times S_2$, $X=(X_1,X_2)$, and $Y=\varphi(X_1,Z)$, then $\P^X_Y$ $=$  $\Lc( \varphi(x_1,Z))\mid_{x_1=X_1}$,  and thus $\P^X_Y$ is $\sigma(X_1)$-measurable in 
$(\Pc(T),\Cc(T))$.
\item For any probability kernel $\nu$ from $S$ to $\Pc(T)$,  there exists a measurable function $\phi$ $:$ 
$S\times [0,1]\rightarrow T$ s.t. $\nu(s)$ $=$ $\Lc( \phi(s,U))$, for all $s\in S$, where $U$ is a uniform random variable. 
\end{enumerate}
\end{Lemma}
\noindent{\bf Proof.}
The first assertion is stated in  Theorem 6.3 in \cite{Kallenberg}, and the second one follows  from Fubini's theorem.  The third assertion is a consequence of the two others. 
\ep



\begin{Proposition}\label{condlaw-adapted}
Given an open-loop control $\alpha$ $\in$ $\Ac$, and an initial condition $\xi$ $\in$ $L^0(\Xc;\Gc)$, the solution  $X^{\xi,\alpha}$ to the conditional McKean-Vlasov equation is 
such that: for all  $t\in\N$, $X^{\xi,\alpha}_t$ is $\sigma(\xi,\Gamma, (\eps)_{s\leq t},(\eps^0_s)_{s\leq t})$-measurable, 
and $\P^0_{(X^{\xi,\alpha}_t ,\alpha_t)}$ is $\Fc^0_t$-measurable.  
\end{Proposition}
{\bf Proof.}
We prove the result by induction on $t$. It is clear for $t=0$. Assuming that it holds true for some $t\in\N$, we write
\beqs
X^{\xi,\alpha}_{t+1} =  F(X^{\xi,\alpha}_t, \alpha_t, \P^0_{(X^{\xi,\alpha}_t,\alpha_t)},\eps_{t+1}, \eps^0_{t+1}),  \quad 	t\in\N.
\enqs
By induction hypothesis,  there is  a measurable function $f_{t+1}$ $:$ $\Xc\times G\times E^{t+1}\times (E^0)^{t+1}$ $\rightarrow$ $\Xc$ s.t. 
$X^{\xi,\alpha}_{t+1}$ $=$ $f_{t+1}(\xi, \Gamma, (\eps_s)_{s\leq t+1}, (\eps^0_s)_{s\leq t+1})$, and thus $X^{\xi,\alpha}_{t+1}$ is $\sigma(\xi,\Gamma, (\eps)_{s\leq t+1},(\eps^0_s)_{s\leq t+1})$-measurable and  $\P^0_{(X^{\xi,\alpha}_{t+1},\alpha^{}_{t+1})}$ is $\sigma(\eps^0_s,s\leq t+1)$-measurable by Lemma \ref{condlaw}. 
 \ep

\section{Wasserstein convergence of the empirical law}

\setcounter{Lemma}{0}
\setcounter{Definition}{0}


\vspace{2mm}

\begin{Proposition}[Conditional Wasserstein convergence of the empirical measure]\label{annexe_conditional_wassertein_convergence}
Let $E,F$ be two measurable spaces and $G$ a compact Polish space. Let $X$ be an $E$-valued random variable independent from  a family of i.i.d. $F$-valued variables $(U_i)_{i\in\N}$,  and a measurable function $f:E\times F\rightarrow G$. Then 
\beqs
\Wc\Big(\frac{1}{N}\sum_{i=1}^N\delta_{f(X,U_i)},\P^X_{f(X,U_1)}\Big) & \overset{a.s.}{\underset{N\rightarrow \infty}{\longrightarrow}}& 0. 
\enqs
\end{Proposition}
\noindent{\bf Proof.}
It suffices to observe that the probability of this event is one by conditioning w.r.t. $X$ and use the analog non-conditional result, which follows from the fact that 
Wasserstein distance metrizes weak convergence (as $G$ is compact), and the fact that empirical measure converges weakly.  
\ep

\section{Proof of coupling results} \label{coupling}

\setcounter{Lemma}{0}
\setcounter{Definition}{0}

 \begin{Lemma} \label{lemc2}
 Let  $U,V$ be two independent  uniform variables, and $F$ a distribution function on $\R$. We have
 \beqs 
 \Big(F^{-1}(U), F(F^{-1}(U))-U\Big)\overset{d}{=} (F^{-1}(U), V \Delta F (F^{-1}(U)),
 \enqs
 where we denote $\Delta F$ $:=$ $F-F_-$. 
 \end{Lemma}
 {\bf Proof.}
 Notice that $F(F^{-1}(U))-U$ is the position (from top to bottom) of $U$ in the set $\{u\in [0,1], F^{-1}(u)=F^{-1}(U)\}$ and is thus smaller than $\Delta F(F^{-1}(U))$. Now, given a measurable function $f\in L^0(A\times [0,1];\R)$, we have 
 \begin{align}\label{eqcoup}
&\E \Big[ f\big(F^{-1}(U), F(F^{-1}(U))-U) \Big]\\ &  =  \E \big[  f( F^{-1}(U), 0){\bf 1}_{\Delta F(F^{-1}(U))=0} \big] + \E \big[  f( F^{-1}(U), F(F^{-1}(U))-U)){\bf 1}_{\Delta F(F^{-1}(U))>0} \big]. 
 \end{align}
 The second term can be decomposed as 
  \beqs
 \sum_{\Delta F(c)>0}\E \Big[ f\big( c, F(c) - U\big){\bf 1}_{F^{-1}(U)=c} \Big] &=&  \sum_{\Delta F(c)>0}\int_0^1 f\left( c, \Delta F(c)u)\right)\Delta F(c)du. 
 \enqs
  where the equality comes from a change of variable. Summing over $\Delta F(c)>0$, we obtain $\E\big[f\big( F^{-1}(U), V \Delta F(F^{-1}(U))\big){\bf 1}_{\Delta F (F^{-1}(U))>0}\big]$, and combined with \eqref{eqcoup}, we get 
 $$\E \Big[ f\big(F^{-1}(U), F(F^{-1}(U))-U) \Big] =\E\big[f\big( F^{-1}(U), V \Delta F(F^{-1}(U))\big)\big],$$
 which proves the result. 
 \ep

 \vspace{13mm}

 \begin{Lemma}\label{annexe_embedding}
 Let $\Xc$ be a compact Polish space, then there exists an embedding $\phi\in L^0(\Xc,\R)$ such that
 \begin{enumerate}
 \item $\phi$ and $\phi^{-1}$ are uniformly continuous,
 \item for any probability measure $\mu\in\Pc(\Xc)$, we have $\text{Im}\left( F^{-1}_{\phi\star\mu}\right)\subset \text{Im} (\phi)$. In particular, $\phi^{-1}\circ F^{-1}_{\phi\star\mu}$ is well posed.
 \end{enumerate} 
 \end{Lemma}
 \noindent{\bf Proof.}
  \noindent 1. Without loss of generality, we assume that $\Xc$ is bounded by $1$. Fix a countable dense family $(x_n)_{n\in\N}$ in $\Xc$. We define the map $\phi_1 :x\in\Xc \overset{}{\mapsto} (d(x,x_n))_{n\in\N}\in [0,1]^\N$. Let us endow $[0,1] ^\N$ with the metric 
 $\mfd((u_n)_{n\in\N},(v_n)_{n\in\N}):=\sum_{n\geq 0}\frac{1}{2^n}\vert u_n -v_n\vert$. $\phi_1$ is clearly injective and uniformly continuous (even Lipschitz). The compactness  of $\Xc$ implies that its inverse $\phi_1^{-1}$ is uniformly continuous as well. Let us now consider 
 $\phi_2: ([0,1]^\N ,\mfd)\mapsto [0,1]$ where $\phi_2((u_n)_{n\in \N})$ essentially groups the decimals of the real numbers $u_n$, $n\in\N$, in a single real number. More precisely, let $\iota:\N\rightarrow \N^2$ be a surjection, then we define the $k$-th decimal of $\phi_2((u_n)_{n\in \N}) $ as the $(\iota (k))_2$-th decimal of $u_{(\iota (k))_1}$ (with the convention that for a number with two possible decimal representations, we choose the one that ends with $000...$). $\phi_2$ is clearly injective, uniformly continuous, as well as its inverse $\phi_2^{-1}$. Thus, $\phi:=\phi_2\circ \phi_1$ defines an embedding of $\Xc$ into $\R$, such that $\phi$ and $\phi^{-1}$ are uniformly continuous.

\noindent  2. $F^{-1}_{\phi\star \mu}$ being caglad, and $\text{Im}(\phi)$ being closed (by compactness  of $\Xc$), it is enough to prove that $F^{-1}_{\phi\star \mu}(u)\in \text{Im}\phi$ for almost every $u\in[0,1]$ (in the Lebesgue sense). However, given a uniform variable $U$, we have $F^{-1}_{\phi\star \mu}(U)\sim \phi\star \mu$, and thus
\beqs
\P(F^{-1}_{\phi\star \mu}(U)\in\text{Im}(\phi))=\P_{Y\sim \mu}(\phi(Y)\in\text{Im}(\phi))=1.
\enqs
 \ep

\vspace{5mm}

\noindent {\bf Proof of Lemma \ref{lemcoupling1} }

 \noindent (1) We first consider the case where $\Xc$ $\subset$ $\R$.  Let us call $F_\mu$ the distribution function of $\mu$ $\in$ $\Pc(\Xc)$, and $F^{-1}_\mu$ its generalized inverse. 
 Let us define the function $\zeta$ $:$ $\Pc(\Xc)\times\Pc(\Xc)\times\Xc\times [0,1]$ $\rightarrow$ $\Xc$ by 
 \beqs
 \zeta(\mu,\mu', x,u) & := &  F_{\mu'}^{-1}\big(F_\mu(x)- u \Delta F_\mu (x)\big),
 \enqs
which is measurable by noting that the measurability in $\mu,\mu'$ comes from the continuity of
 \beqs
 \Pc(\Xc) &\rightarrow& L_{caglad}^1([0,1],\Xc) \\
 \mu &\mapsto& F_\mu^{-1}.
 \enqs
 By construction, we then have for any $\xi$ $\sim$ $\mu$, and $U,V$ two independent  uniform variables,  independent of $\xi$
 \beqs
 (\xi, \zeta(\mu,\mu',\xi,V)) &=& (\xi, F_{\mu'}^{-1}\big(F_\mu(\xi)- V \Delta F_\mu (\xi)\big))\\
 &\overset{d}{=}& (F_\mu^{-1}(U), F_{\mu'}^{-1}\big(F_\mu(F_\mu^{-1}(U))- V \Delta F_\mu (F_\mu^{-1}(U))\big))\\
  &=& (F_\mu^{-1}(U), F_{\mu'}^{-1}\big(F_\mu(F_\mu^{-1}(U))- V \Delta F_\mu (F_\mu^{-1}(U))\big))\\
   &\overset{d}{=}& \big(F_\mu^{-1}(U), F_{\mu'}^{-1}(U)\big),   
 \enqs
 where the last equality holds by Lemma \ref{lemc2}. It is well-known (see e.g. Theorem 3.1.2 in \cite{rachev}) that $\big(F_\mu^{-1}(U), F_{\mu'}^{-1}(U)\big)$ is an optimal coupling for $(\mu,\mu')$, and so 
 $\Wc(\mu,\mu')$ $=$ $\E\big[d(\xi,\zeta(\mu,\mu',\xi,V))\big]$.  
 
 \vspace{2mm}

 \noindent (2) Let us now consider the case of a general compact Polish space $\Xc$.  Denoting  by $\zeta_\R$ the "$\zeta$" from the case "$\Xc\subset\R$", 
 and considering  an embedding $\phi\in L^0(\Xc,\R)$ as in Lemma \ref{annexe_embedding},  let us define
 \beqs
 {\zeta (\mu,\mu',x,u):=\phi^{-1}(\zeta_\R (\phi\star \mu , \phi\star \mu', \phi(x),u))},
 \enqs
which is well posed by definition of $\zeta_\R$ and Lemma \ref{annexe_embedding}. Now, fix $\xi\sim  \mu$, $U$ a uniform variable independent of $\xi$, and define $\xi':=\zeta(\mu,\mu',\xi,U)$. By definition of $\zeta$, its clear that $\xi'\sim \mu'$, and
 \begin{align}\label{eqphi}
\E \big[ d(\phi(\xi),\phi(\xi')) \big] = \Wc(\phi\star\mu, \phi\star\mu').  	
 \end{align}
Fix $\epsilon>0$.  We are looking for $\eta,\delta>0$ such that  
\beqs
\Wc(\mu,\mu')<\eta \; \Rightarrow \;  \Wc(\phi\star\mu, \phi\star\mu') <\delta \; \Leftrightarrow \;  \E\big[d(\phi(\xi),\phi(\xi'))\big] <\delta \; \Rightarrow \;  \E[d(\xi,\xi')] <\epsilon. 
\enqs
 Let us first show that there exists $\delta>0$ such that ${\E [d(\phi(\xi),\phi(\xi'))]<\delta\Rightarrow \E[d(\xi,\xi')]<\epsilon}$. 
 Fix $\gamma>0$ such that $d(x,x')<\gamma\Rightarrow d(\phi^{-1}(x),\phi^{-1}(x'))<\frac{\epsilon}{2}$. Denoting by $\Delta_\Xc$ the diameter of $\Xc$, we then  have 
\beqs
\E[d(\xi,\xi')] \leq \E [d(\xi,\xi'){\bf 1}_{d(\phi(\xi),\phi(\xi'))<\gamma}]+ \frac{\Delta_\Xc}{\gamma}\E\big[d(\phi(\xi),\phi(\xi'))\big] \leq \frac{\epsilon}{2}+\frac{\Delta_\Xc}{\gamma} \E\big[d(\phi(\xi),\phi(\xi'))\big],
\enqs
so that we can choose $\delta= \frac{\gamma}{\Delta_\Xc }\frac{\epsilon}{2}$. On the other hand, by uniform continuity of $\phi$ and by definition of the Wasserstein metric, there exists $\eta>0$ such that 
$d(\mu,\mu')<\eta\Rightarrow \Wc(\phi\star\mu, \phi\star\mu') <\delta$. From \eqref{eqphi}, we thus conclude that $d(\mu,\mu')<\eta \Rightarrow \E[d(\xi,\xi')]<\epsilon$. 
 \ep

\vspace{9mm}

\small

\bibliographystyle{plain}

\bibliography{references}




  




\end{document}